\documentclass[journal, onecolumn]{IEEEtran}
\usepackage{amsthm}

\usepackage{tikz}
\usetikzlibrary{arrows.meta, positioning, shapes,snakes, calc}
\usepackage{subfig}
\usepackage{amssymb}
\usepackage{mathtools}
\usepackage{cases}
\usepackage{autobreak}
\usepackage{mathrsfs}  
\usepackage{booktabs}
\usepackage{multirow}
\usepackage{color}
\usepackage{hyperref}
\usepackage{cite}
\usepackage{amsmath}
\usepackage{extarrows}
\usepackage{bm}
\usepackage{comment}
\usepackage{algorithmic}
\usepackage{algorithm}
\usepackage{caption}
\usepackage{graphicx}
\usepackage{pifont}
\usepackage{bbm}
\usepackage{dsfont}
\usepackage{diagbox}
\newcommand{\BA}{{\bold{A}}}
\newcommand{\BB}{{\bold{B}}}
\newcommand{\BC}{{\bold{C}}}
\newcommand{\BD}{{\bold{D}}}
\newcommand{\BE}{{\bold{E}}}
\newcommand{\BF}{{\bold{F}}}
\newcommand{\BG}{{\bold{G}}}
\newcommand{\BH}{{\bold{H}}}
\newcommand{\BI}{{\bold{I}}}
\newcommand{\BJ}{{\bold{J}}}
\newcommand{\BK}{{\bold{K}}}

\newcommand{\BM}{{\bold{M}}}

\newcommand{\BQ}{{\bold{Q}}}
\newcommand{\BR}{{\bold{R}}}
\newcommand{\BS}{{\bold{S}}}
\newcommand{\BT}{{\bold{T}}}
\newcommand{\BU}{{\bold{U}}}

\newcommand{\BX}{{\bold{X}}}
\newcommand{\BY}{{\bold{Y}}}
\newcommand{\BZ}{{\bold{Z}}}

\newcommand{\Ba}{{\bold{a}}}
\newcommand{\Bb}{{\bold{b}}}

\newcommand{\Bd}{{\bold{d}}}

\newcommand{\Bg}{{\bold{g}}}
\newcommand{\Bh}{{\bold{h}}}

\newcommand{\Bn}{{\bold{n}}}

\newcommand{\Bu}{{\bold{u}}}
\newcommand{\Bv}{{\bold{v}}}

\newcommand{\Bx}{{\bold{x}}}
\newcommand{\By}{{\bold{y}}}
\newcommand{\Bz}{{\bold{z}}}

\newcommand{\BGamma}{\bold{\Gamma}}
\newcommand{\BPsi}{\bold{\Psi}}
\newcommand{\BOmega}{\bold{\Omega}}
\newcommand{\BTheta}{\bold{\Theta}}
\newcommand{\BSigma}{\bold{\Sigma}}

\newcommand{\BPi}{\bold{\Pi}}

\newcommand{\BPhi}{\boldsymbol{\Phi}}

\newcommand{\Bmu}{\boldsymbol{\mu}}

\newcommand{\Bxi}{\boldsymbol{\xi}}

\newcommand{\Bkappa}{\boldsymbol{\kappa}}

\newcommand{\Bdelta}{\boldsymbol{\delta}}

\newcommand{\SI}{\mathscr{I}}

\newcommand{\Bone}{\bold{1}}
\newcommand{\Bzero}{\bold{0}}

\newcommand*{\rom}[1]{\expandafter\@slowromancap\romannumeral #1@}

\newcommand{\E}{\mathbb{E}}
\newcommand{\PP}{\mathbb{P}}
\newcommand{\dd}{\mathrm{d}}

\newcommand{\C}{\mathbb{C}}
\newcommand{\R}{\mathbb{R}}
\newcommand{\Ind}{\mathbbm{1}}

\DeclareMathOperator{\Tr}{Tr}
\DeclareMathOperator{\dist}{dist}

\DeclareMathOperator{\diag}{diag}
\DeclarePairedDelimiter\abs{\lvert}{\rvert}%
\DeclarePairedDelimiter\norm{\lVert}{\rVert}%
\makeatletter
\let\oldabs\abs
\def\abs{\@ifstar{\oldabs}{\oldabs*}}
\let\oldnorm\norm
\def\norm{\@ifstar{\oldnorm}{\oldnorm*}}
\makeatother
\newcommand\numberthis{\addtocounter{equation}{1}\tag{\theequation}}
\newtheorem{remark}{Remark}
\newtheorem{theorem}{Theorem}
\newtheorem{lemma}{Lemma}
\newtheorem{proposition}{Proposition}
\newtheorem{assumption}{Assumption}
\newtheorem{corollary}{Corollary}

\ifCLASSINFOpdf
\else
\fi

\title{No Eigenvalues Outside the Limiting Support of
Generally Correlated and Noncentral Sample Covariance Matrices}

\author{\IEEEauthorblockN {Zeyan Zhuang, \textit{Graduate Student Member}, \textit{IEEE}, Xin Zhang, \textit{Member}, \textit{IEEE}, Dongfang Xu, \textit{Member}, \textit{IEEE}, and Shenghui Song, \textit{Senior Member}, \textit{IEEE}}

\thanks{
Zeyan Zhuang and Shenghui Song are with the Department of Electronic and Computer Engineering,
The Hong Kong University of Science and Technology, Hong Kong (e-mail: zzhuangac@connect.ust.hk; eeshsong@ust.hk).
\par
Xin Zhang is with the School of Cyber
Science and Technology, Beihang University, Beijing 100191, China (e-mail: zhangxinn@buaa.edu.cn).
\par
Dongfang Xu is with the Division of Integrative Systems and Design,
The Hong Kong University of Science and Technology, Hong Kong (e-mail: eedxu@ust.hk).
}
}
\begin{document}
\maketitle
\begin{abstract}
Spectral properties of random matrices play an important role in statistics, machine learning, communications, and many other areas. Engaging results regarding the convergence of the empirical spectral distribution (ESD) and the ``no-eigenvalue'' property have been obtained for random matrices with different correlation structures. However, the related spectral analysis for generally correlated and noncentral random matrices is still incomplete, and this paper aims to fill this research gap. Specifically, we consider matrices whose columns are independent but with non-zero means and non-identical correlations. Under high-dimensional asymptotics where both the number of rows and columns grow simultaneously to infinity, we first establish the almost sure convergence of the ESD for the concerned random matrices to a deterministic limit, assuming mild conditions. Furthermore, we prove that with probability 1, no eigenvalues will appear in any closed interval outside the support of the limiting distribution for matrices with sufficiently large dimensions. The above results can be applied to different areas such as statistics, wireless communications, and signal processing. In this paper, we apply the derived results to two communication scenarios: 1) We determine the limiting performance of the signal-to-interference-plus-noise ratio for multi-user multiple-input multiple-output (MIMO) systems with linear minimum mean-square error receivers; and 2) We establish the invertibility of zero-forcing precoding matrices in downlink MIMO systems, providing theoretical guarantees.
\end{abstract}
\section{Introduction}
The covariance structure of random matrices is fundamental in many areas, including multivariate statistics \cite{yao2015sample}, machine learning \cite{couillet_liao_2022}, and wireless communications \cite{couillet2011random}, particularly regarding the properties of eigenvalue distribution.
Consider a $p$-by-$n$ data matrix $\BSigma = [\Bxi_1, \Bxi_2, ..., \Bxi_n]$ comprising $n$ independent and $p$-dimensional observations.
The sample covariance matrix is defined as $
    \BS = \sum_{j=1}^n \Bxi_j\Bxi_j^H = \BSigma \BSigma^H
$.
Classical statistical frameworks typically assume data abundance, where the dimension $p$ is fixed while the sample size $n \to \infty$. If the elements of $\BSigma$ are independent and identically distributed (i.i.d.) with $\E[\Bxi_j]_i = 0$ and $\E|\sqrt{n}[\Bxi_j]_i|^2 = 1$,  the strong law of large numbers ensures that $\BS$ converges to the identity matrix $\BI_p$ almost surely.
\par
However, the classical paradigm breaks down in the modern ``big data'' era, where high-dimensional data is prevalent. A key question arises: how do the eigenvalues of $\BS$ behave when the dimension $p$ and sample size $n$ increase simultaneously? Seminal work \cite{marchenko1967distribution} showed that as $p$, $n \to \infty$ with $p / n \to c \in (0, \infty)$, the empirical spectral distribution (ESD) of $\BS$ converges almost surely to a deterministic limit known as the Marčenko-Pastur (MP) distribution. Specifically, when $c \in (0, 1)$, the MP distribution is supported over $[(1 - \sqrt{c})^2, (1 + \sqrt{c})^2]$. The analysis for the limiting spectral behavior has been extended to various statistical models \cite{SILVERSTEIN1995Empirical, SILVERSTEIN1995strong, Zhanglixin2006spectral, MeiTianxing2023singular,  Dozier2007information, hachem2007deterministic,  hachem2006varianceprofile, zhou2023limiting}. An important class is the noncentral model \cite{hachem2006varianceprofile, Dozier2007information, hachem2007deterministic, zhou2023limiting}, which generalizes the centered model. The noncentral model is also known as the information-plus-noise model, which has significant applications in image processing and matrix denoising \cite{couillet_liao_2022}. 
\par
While the above-mentioned limiting spectral distributions indicate the proportion of eigenvalues within a specific region, they can not guarantee that extreme eigenvalues (the smallest and largest) are within the limiting support. Given the importance of these extreme eigenvalues, their asymptotic behavior has been studied with compelling results. For i.i.d. cases, Yin \textit{et al.} \cite{yin1988limit} established that if $\mathbb{E}|\sqrt{n}[\Bxi_j]_i|^4 < \infty$, the largest eigenvalue $\lambda_{\max}(\mathbf{S})$  almost surely converges to $(1 + \sqrt{c})^2$, the right endpoint of the MP support. The almost sure convergence of the smallest eigenvalue to the left endpoint was established in \cite{Bai1993smalleig}. More generally, in \cite{bai1998no}, Bai \textit{et al.} demonstrated that almost surely, no eigenvalues of $\BS$ will appear in any interval outside the limiting support when ${ \{ \sqrt{n} \Bxi_j\} }$ are centered with identical covariance. This no-eigenvalue property is vital for studying spiked eigenvalues \cite{BAIK2006spiked} and linear spectral statistics (LSS) \cite{Bai2004CLTaop} in large random matrices.
\par
% An important class of random matrix models is the information-plus-noise model \cite{hachem2006varianceprofile, hachem2013bilinear}, which is crucial for data analysis. Notably, this model generalizes the centered case. The limiting spectrum of the noncentral model has been studied in \cite{girko2012theory, Dozier2007analysis, hachem2006varianceprofile,  hachem2007deterministic, zhou2023limiting}, and the non-eigenvalues property was recently investigated in \cite{Bai2025noeigen}. However, these works considered the case of homogeneous noise, i.e., the correlation $ n\E \Bxi_j \Bxi_j^H$ does not vary with index $j$. To the best of our knowledge, the research on the analysis of the spectrum of the general correlated noncentral model is still in its early stage.
%  \par
However, the related investigations for generally correlated and noncentral models are quite limited. In this paper, we study the spectral behavior of the covariance matrix for such random matrices, whose columns are noncentral and have non-identical correlations. In particular, we establish the convergence of the ESD and prove that, with probability 1, no eigenvalues lie in any specific closed interval outside the limiting support.
\subsection{Related Works}
In this section, we review the related works on the convergence of the ESD and the no-eigenvalue property.  
\subsubsection{Convergence of the ESD}  The study of the ESD for $\BS$ has a long history, accompanied by increasing complexity in the correlation and mean structures. The earliest results can be traced back to \cite{marchenko1967distribution}, where Marčenko and Pastur investigated the case with
$
 \Bxi_j = \sqrt{\tau}_j \Bx_j
$,
where $ {\tau_j} $ is a sequence of non-negative numbers satisfying specific ergodicity conditions, and $\{ {\sqrt{n}\Bx_j} \}$ is a sequence of i.i.d. and centered random vectors with $\mathrm{Cov}(\sqrt{n} \Bx_j) = \BI_p$ and finite forth moment, leading to the celebrated MP law. Subsequently,  Silverstein \textit{et al.} \cite{SILVERSTEIN1995Empirical, SILVERSTEIN1995strong} studied the ESD for the case with
$
\Bxi_j = \BT^{1/2} \Bx_j
$, where $\BT$ is Hermitian non-negative. In particular, it was established that, when the spectral distribution of $\BT$ converges and the elements of $\Bx_j$ are i.i.d. with zero mean and
variance $1 / n$, the ESD of $\BS$ converges weakly to the generalized MP distribution. Motivated by the advancement in multivariate statistics,  Zhang \cite{Zhanglixin2006spectral} generalized these results to the separable correlated model, which takes the form
$
 \Bxi_j = \sqrt{\tau_j} \BT^{1/2} \Bx_j , 
$
where $\BT$ is Hermitian non-negative, $\{\tau_j\}$ is a sequence of non-negative numbers, and $\{ \sqrt{n}\Bx_j \}$ are random vectors with i.i.d. standardized elements and finite fourth moments. Recently, these findings were extended by Mei \textit{et al.} in \cite{MeiTianxing2023singular} to consider  a generally correlated model
$ \Bxi_j = \BB_{j} \Bx_j  $,
where $\{\BB_j\}$ is a sequence of deterministic matrices, and $\{\Bx_j\}$ is a sequence of independent random vectors that satisfy mild concentration conditions.
\par
The aforementioned works considered centered random matrices. In the following, we review the related works for noncentral models. In \cite{Dozier2007information}, Dozier \textit{et al.} studied the matrix 
$ \BS = (\BA + \sigma \BX)(\BA + \sigma \BX)^H$,
where $\BA$ and $\BX$ are independent, $\BX$ has i.i.d. standardized entries, and the ESD of $\BA\BA^H$ converges to a non-random distribution. It was shown in \cite{Dozier2007information} that the limiting spectral distribution (LSD) of $\BS$ exists and its Stieltjes transform is uniquely defined by the LSD of $\BA\BA^H$. In \cite{zhou2023limiting}, the convergence of the ESD for the covariance
$ \BS = (\BA + \BT^{1/2} \BX)(\BA + \BT^{1/2}\BX)^H$, where $\BT$ is Hermitian non-negative, was established under the condition that $\BA\BA^H$  and  $\BT$ are commutative and their eigenvalue distributions converge. In \cite{hachem2006varianceprofile}, Hachem \textit{et al.} studied the convergence of the ESD for the variance profile model
$ \BS = (\BY + \BA)(\BY +  \BA)^H$, where $\BA$ is a rectangular diagonal matrix. The entries of $\BY$ are given by $[\BY]_{i, j} = \sigma(i / p, j /n)X_{ij} / \sqrt{n}$, where $\sigma: [0, 1]^2 \to \mathbb{R}$ denotes the variance profile function and $X_{ij}$ is standardized with finite $4 + \varepsilon$ moment. However, as emphasized in \cite{hachem2007deterministic},  when $\BA$ is a general matrix, the ESD of $\BS$ may not converge even if the eigenvalue distribution of $\BA\BA^H$ converges. To tackle that issue, an alternative approach is to investigate the convergence properties of the resolvent of $\BS$, which characterizes the convergence of linear spectral statistics (LSS) even if the limit of ESD does not exist. In \cite{hachem2007deterministic}, the deterministic approximation of the trace form for the resolvent of the noncentral covariance model $ \BS = (\BY + \BA)(\BY + \BA)^H $ was established,
where the rows and columns of $\BA$ have bounded Euclidean norms. This result generalizes that in \cite[Chapter 7]{girko2012theory}, where Girko established an approximation rule for the resolvent under the condition that the entries of $\Bxi_j$ are independent, despite differing variances. However, in practice, the correlations of $\Bxi_j$ can be more complex, necessitating the consideration of a generally correlated model.

\subsubsection{No-Eigenvalue Property} Following the pioneering work in \cite{bai1998no}, Dozier \textit{et al.} extended the no-eigenvalue property to the separable correlated model in \cite{PAUL2009Noeigen}. 
 In \cite{Yiming2023random}, a similar phenomenon was observed for the noncentral model when the mean matrix $\BA$ has finite rank. This result was further extended to the moderate low-rank case in \cite{Xiaoyu2025spiked}, specifically when the rank of $\BA$ is of the order $O(n^{1/3})$. In \cite{Bai2012Noeigen}, Bai \textit{et al.} demonstrated that for the noncentral random matrix with isotropic correlations, i.e., $\mathrm{Cov}(\sqrt{n} \Bxi_j) = \BI_p$ and general $\BA$, the no-eigenvalue property holds. This result was later generalized in \cite{Bai2025noeigen} to the case with $ \mathrm{Cov}(\sqrt{n} \Bxi_j) = \mathbf{T}$ where $\BT$ is Hermitian non-negative.
In \cite{AblaNoEigen}, Kammoun \textit{et al.} proved that no eigenvalues exist outside the support of the LSD for generally correlated and centered Gaussian random matrices. In \cite{Yin2022strong}, Yin showed that the no-eigenvalue phenomenon occurs when $\Bxi_j = \BB \Bx_j$ and the correlation matrix $\mathbf{B}$ can take any form. However, the research on the no-eigenvalue property for noncentral models is still in its infancy. In particular, the results for the variance profile model, which serves as a special case of the generally correlated model, have not been established yet.
\subsection{Contributions}
The contributions of this work are listed as follows.
\par
\subsubsection{Convergence of the Resolvent} For generally correlated and noncentral covariance matrices, we derive the deterministic approximation of the resolvent. In particular, we prove that the bilinear and trace forms of the resolvent converge to those of a deterministic matrix, which is characterized by a set of fixed-point equations. Based on the convergence of the resolvent, we obtain the functional equation for the LSD and demonstrate that the ESD will converge to the LSD when the correlation matrices and the mean satisfy commutative conditions.
 
\subsubsection{Analysis of the Fundamental Equations} We establish the analytical properties of the solutions to the fixed-point equations. In particular, the solutions belong to a class of functions that are Stieltjes transforms of non-negative measures. We prove that under certain mild conditions, the support of these measures associated with the solutions is a subset of the support of the LSD.
 
\subsubsection{No-Eigenvalue Property} We prove that, almost surely, no eigenvalues of the generally correlated and noncentral matrices will appear in any specific interval outside the support of the LSD. Based on this result, we then show that the largest eigenvalue of the covariance matrix is almost surely finite when the support of the LSD is bounded.

\subsubsection{Applications} We apply the theoretical results to communication systems and obtain new insights. With the deterministic approximation of the resolvent, we demonstrate the almost sure convergence of the signal-to-interference-plus-noise ratio (SINR) for uplink multi-user multiple-input multiple-output (MIMO) systems with linear minimum mean-square error (LMMSE) receivers. Additionally, based on the no-eigenvalue property, we prove that for downlink multi-user MIMO  communications, the Gram matrix of the channel is almost surely invertible, ensuring that the zero-forcing precoding is valid.
\subsection{Paper Outline}
The paper is organized as follows. In Section \ref{Sec_Problem_statement}, we introduce the statistical model and assumptions. In Section \ref{Sec_Main_results}, we present the main results: 1) almost sure convergence of the trace/bilinear forms for the resolvent, and 2) no eigenvalues outside the limiting support. In Section \ref{Sec_applications}, we apply these results to MIMO systems, analyzing LMMSE receivers for uplink multi-user MIMO systems and zero-forcing precoding for downlink multi-user MIMO systems. Numerical simulations are provided in Section \ref{Sec_Simulations} to validate the accuracy of the theoretical results and Section \ref{Sec_Conclusion} concludes this paper. The key mathematical tools are provided in Appendix \ref{App_mathmatical_tools}.
 \subsection{Notations} 
 We adopt the following notations throughout the paper.
 \par
\subsubsection{Notations for Matrices and Vectors} We use bold 
 upper-case letters and bold lower-case letters to represent matrices and vectors, respectively. The $(i, j)$-th element of matrix $\mathbf{A}$ and the $i$-th element of vector $\mathbf{a}$ are denoted by $[\mathbf{A}]_{i, j}$ and  $[\mathbf{a}]_i$, respectively. $\mathbf{A}^T$ and $\mathbf{A}^H$ represent the transpose or conjugate transpose of $\mathbf{A}$, and $\mathbb{C}^p$ and $\mathbb{C}^{p \times n}$ denote the space of $p$-dimensional complex vectors and the $p$-by-$n$ complex matrix space, respectively.  The conjugate transpose operator and transpose operator are denoted by $(\cdot)^T$ and $(\cdot)^H$, respectively. $\Tr \BA$ refers to the trace of $\BA$.
\par
\subsubsection{Norms for Matrices and Vectors} When dealing with matrices, $\norm{\cdot}$ denotes the spectral norm, $\norm{\cdot}_{\infty}$ refers to the max-row norm, and $\norm{\cdot}_1$ represents the max-column norm. In particular, $\norm{\BA}_{\infty} = \max_{1 \leq i \leq n} \sum_{j=1}^n |[\BA]_{i, j}|$ and $\norm{\BA}_{\infty} = \max_{1 \leq j \leq n} \sum_{i=1}^n |[\BA]_{i, j}|$ if $\BA$ is a $n$-by-$n$ matrix. $\rho(\BA)$ represents the spectral radius of matrix $\BA$. In the case of vectors, $\norm{\cdot}$ denotes the Euclidean norm, and $\norm{\cdot}_{\infty}$ represents the $\ell_{\infty}$-norm.
\par
\subsubsection{Constants and Functions} We use $\Bzero_n$ and $\Bzero_{p\times n}$ to denote the $n$-dimensional all-zero vector and $p$-by-$n$ all-zero matrix, respectively, while $\Bone_n$ and $\Bone_{p\times n}$ represent the corresponding all-one vector and matrix, respectively.  $\BI_n$ represents the identity matrix of size $n$. We denote $[n] = \{1, \ldots, n\}$ as the set of positive integers, and the complex unit is denoted by $\jmath$. $\Ind_{\{\cdot\}}$ or $\Ind\{\cdot\}$ represents the indicator function.
$\Im(a)$ and $\Re(a)$ denote the imaginary and real part of complex number $a$, respectively.
The set $\R_* = \R - \{0 \}$ denotes the real axis except $0$. $\R_+$ and $\C^+$ represent the sets $\{x \in \R: x \geq 0\}$ and $\{z \in \mathbb{C}: \Im(z) > 0\}$, respectively. The support of measure $\mu$ is denoted by $\mathsf{Supp}(\mu)$.
\par
\subsubsection{Probability Measures} The notation $\PP(\cdot)$ denotes the probability measure and $\E (\cdot)$ represents the expectation. The notations $\xrightarrow{a.s.}$ and $\Rightarrow$ indicate convergence almost surely and convergence in distribution, respectively. The event $\{E_n, \text{i.o.}\}$ represents that the event $E_n$ occurs infinitely often, i.e., $\{E_n,  \text{i.o.}\} = \cap_{m=1}^{\infty} \cup_{n=m}^{\infty} E_n$. 
\par
\subsubsection{Asymptotic Notations} A constant is denoted by $K$ and $K_{\Xi}$ represents the constant related to the parameter(s) $\Xi$. The values of the constants may vary from line to line.  We use the asymptotical notation $a \lesssim b$ and $a \lesssim_{\Xi} b$ if and only if $\limsup |a / b| \leq K$ (almost surely) and $\limsup |a / b| \leq K_{\Xi}$ (almost surely), respectively.

\section{Problem Statement}
\label{Sec_Problem_statement}
In this section, we introduce the statistical model and the associated problem formulation. Consider the random matrices 
\begin{equation}
    \BSigma = \BA + \BY \in \mathbb{C}^{p \times n},
\end{equation}
where $\BA$ is deterministic and $\BY$ is random. In many scenarios, $\BA$ represents the deterministic signal component, while $\BY$ denotes the random noise. In this paper, we consider $\BY$ with the following general correlation structure
\begin{equation}
    \BY = \frac{1}{\sqrt{n}} \begin{bmatrix}
        \BB_1 \Bx_1 &\BB_2 \Bx_2& \ldots & \BB_n \Bx_n
    \end{bmatrix}, \label{Eq_random_matrix_Model}
\end{equation}
where $\Bx_j \in \C^{d_j}, j \in [n]$ are i.i.d. random vectors with zero mean and $\BB_j \in \mathbb{C}^{p \times d_j}$ denotes the correlation of the $j$-th column of $\BY$. The objective of this paper is to investigate the random spectrum of $\BSigma$ in a high-dimensional setting. To this end, we make the following assumptions.
\begin{assumption} (On the Asymptotic Regime)
\label{Assumpt_1}
The dimensions $n$, $p$, and $\{d_j\}_{1 \leq j \leq n}$ approach infinity with
\begin{align}
    0 < \liminf_{n \to \infty} \min_{1 \leq j \leq n} \frac{d_j}{n} \leq \limsup_{n \to \infty} \max_{1 \leq j \leq n} \frac{d_j}{n} < \infty, \text{ and } 0 < \liminf_{n \to \infty} \frac{p}{n} \leq \limsup_{n \to \infty} \frac{p}{n} < \infty, \label{Eq_Assumpt_1_Asym}
\end{align}
\end{assumption}
This assumption is common in the study of high-dimensional random matrices \cite{bai2010spectral,  AblaNoEigen}. In fact, $p = p(n)$ and $d_j = d_j(n)$ for $j \in [n]$ can be viewed as sequences indexed by $n$, i.e., the numbers $p$ and $d_j$ can be viewed as sampled from $\{ p(n) \}_{n \geq 1}$ and $\{ d_j(n) \}_{n \geq j}$, respectively. For notational simplicity, we use $n \to \infty$ to denote the double asymptotic regime in \eqref{Eq_Assumpt_1_Asym}. 
\begin{assumption}(On the Randomness)
\label{Assumpt_2}
Let $\{ X_{ij}\}_{i, j \geq 1}$ be a double array of i.i.d. complex random entries with $\E X_{11} = 0$ and $\E \abs{X_{11}}^2 = 1$. Moreover, there exists $\varepsilon > 0$ suth that $\E \abs{X_{11}}^{4 + \varepsilon} < \infty$. The elements of the random vectors $\{\Bx_j\}$ in \eqref{Eq_random_matrix_Model} come from the double array with $[\Bx_j]_i = X_{ij}$.
\end{assumption}
Here, the number $\varepsilon$ could be arbitrary and we assume $\varepsilon \leq 1 / 4$ without loss of generality. Strictly speaking, the vector ${\Bx}_j$ should be indexed by $n$, i.e., $\mathbf{x}_j = \mathbf{x}_{nj}$. To simplify the notation, we omit $n$ when there is no ambiguity. Note that this omission has been applied to $\BSigma = \BSigma_n$, $\BY = \BY_n$, $\BA = \BA_n$, and $\{\BB_j\} = \{\BB_{nj}\}$.
% \begin{assumption} (On the Correlations)
% \label{Assumpt_3}
% Define $\BOmega_j = \BB_j\BB_j^H$, then 
%   \begin{equation}
%     \omega^{-} \leq \liminf_{n \to \infty} \min_{1 \leq j \leq n}\lambda_{\min}(\BOmega_j)  \leq \limsup_{n \to \infty} \max_{1 \leq j \leq n}\lambda_{\max}(\BOmega_j) \leq \omega^+,
%   \end{equation}
% where $\omega^-$ and $\omega^+$ are positive numbers and $\lambda_{\min}(\BOmega_j)$ and  $\lambda_{\max}(\BOmega_j)$ denote the smallest and largest
% eigenvalues of $\BOmega_j$, respectively.
% \end{assumption}
\begin{assumption} (On the Correlations)
\label{Assumpt_3}
Define $\BOmega_j = \BB_j\BB_j^H$, then we have 
  \begin{equation}
    0 < \liminf_{n \to \infty} \min_{1 \leq j \leq n} \frac{1}{p} \Tr \BOmega_j \leq \limsup_{n \to \infty} \max_{ 1\leq j \leq n}\norm{\BOmega_j} < \infty. \label{Eq_Assumpt_3}
  \end{equation}
\end{assumption}
 Assumption \ref{Assumpt_3} indicates that the largest eigenvalue of the correlation matrix $\BOmega_j$ is asymptotically bounded, and the extreme low-rank case, such as the case with finite rank, will not occur for ${\BOmega_j}$. This assumption is stronger than that in \cite[Assumption A-2]{AblaNoEigen}, which assumes that the smallest eigenvalue of $\BOmega_j$ is bounded away from 0.
\begin{assumption} (On the Mean)
\label{Assumpt_4}
     The spectral norm of $\BA$ is uniformly bounded with 
     \begin{equation}
         \limsup_{n \to \infty}  ~ \norm{\BA}  < \infty.  \label{Eq_Assumpt_4}
     \end{equation}
\end{assumption}
This assumption ensures that the random matrix is ``non-trivial''. Intuitively, if some singular values of $\BA$ diverge to infinity, it implies that the subspace of the signals associated with these singular values have infinite energy, which makes estimation and detection trivial. Compared to \cite[Chapter 7]{girko2012theory}, which assumes $\limsup_n \|\BA\|_1 + \|\BA\|_{\infty} < \infty$, this assumption is more general given  $\|\BA\| \leq \sqrt{\|\BA\|_1 \|\BA\|_{\infty}} $.
% \begin{table*}[!htbp]
%     \centering
%     \label{Table_Comlex}
% \begin{tabular}{|c|c|c|c|c|c|c|c|c| }
% \toprule
% \diagbox{$\BA$}{$\BB_j$}  & $\BI_p$ & $\tau_j \BI_p$& $\BT^{1 / 2}$ & $b_j \BA^{1/2}$ & $\diag(\sigma^2(i / p, j / n))$ &  General $\BB_j$ \\
% \midrule
% $\Bzero_{p \times n}$ & i.i.d. model & Marchenko-Pastur & Sample Covariance & Bi correlated Model & $\diag(\sigma^2(i / p, j / n))$ &  \\
% \bottomrule
% \end{tabular}
% \end{table*} 
\par
Define the sample covariance of $\BSigma$ as $\BS = \BSigma \BSigma^H$ and denote the $p$ eigenvalues of $\BS$ as $\lambda_1 \geq \lambda_2 \geq \ldots \geq \lambda_p \geq 0$. The ESD of $\BS$ is given by
\begin{equation}
    F^{\BS, n}(x) = \frac{1}{p} \sum_{j=1}^p \mathbbm{1}_{\left\{ \lambda_j \leq x \right\}}. \label{EQ_ESD_F_S_n}
\end{equation}
In the following, we will investigate the weak convergence of the random distribution $F^{\BS, n}$ and the distribution for the eigenvalues of $\BS$ outside the limiting support of $F^{\BS, n}$. We discuss the generality and practical applications of the concerned model \eqref{Eq_random_matrix_Model} in the following remark.
\begin{remark}
The random matrix model in \eqref{Eq_random_matrix_Model} is very general. In fact, by appropriately specifying the correlation matrices $\{\BB_j\}_{j \leq n}$ and the mean $\BA$, \eqref{Eq_random_matrix_Model} can degenerate to a series of previously studied models \cite{marchenko1967distribution, SILVERSTEIN1995Empirical, SILVERSTEIN1995strong , hachem2006varianceprofile, Dozier2007information, hachem2007deterministic, MeiTianxing2023singular, Zhanglixin2006spectral, zhou2023limiting}. For instance, when $ \BB_j = \sqrt{\tau_j} \BT^{1/2}$, $j \in [n]$, with $\BT$ being Hermitian non-negative and $\BA = \Bzero_{p \times n}$, \eqref{Eq_random_matrix_Model} reduces to the separable correlated  model \cite{Zhanglixin2006spectral}.  With the general correlation structure, this model can be applied to the analysis of multi-user MIMO communications with linear transceivers \cite{WagnerMISO2012, Kammoun2009clt, Kammoun2019RZF}. Furthermore, this model can also be applied to multivariate time series analysis \cite{HaoyangLiu2015lineartimeseries}, random or non-random signal testing with heterogeneous noise \cite{Yiming2023random}, mean and covariance estimation for high-dimensional diffusion processes with anisotropic co-volatility \cite{Zhengxinghua2011diffusion}, etc. In this paper, we will apply the theoretical results to two communication systems in Section \ref{Sec_applications} to obtain new physical insights.
\end{remark}
% We note that the model in \eqref{Eq_random_matrix_Model} is very general. In fact, by setting $\BB_j$, $j \in [n]$ to
% \begin{equation}
%     \sqrt{\tau_j} \BI_p, ~~  \BT^{1 / 2}, ~~ \sqrt{\tau_j}\BT^{1 / 2}, ~~ \diag(\sigma_{ij}; 1 \leq i \leq p),  ~~ \text{general}, 
% \end{equation}
% and $\BA$ to $\Bzero_{p \times n}$, pseudo-diagonal, or general, \eqref{Eq_random_matrix_Model} can degenerate to a series of the previously studied models \cite{SILVERSTEIN1995Empirical, SILVERSTEIN1995strong , hachem2006varianceprofile, Dozier2007analysis, hachem2007deterministic, MeiTianxing2023singular, Zhanglixin2006spectral, zhou2023limiting}. 
% As a result, the main results of this paper, i.e., Theorem \ref{Thm_conver_ESD} and Theorem \ref{Thm_No_Eigenvalues}, provide a unified generalization and can be applied to these models. 
\section{Main Results}
\label{Sec_Main_results}
To study the spectrum of $\BS$, a very effective method is to investigate the linear functionals of the resolvent for $\BS$, which is defined as
\begin{equation}
 \BQ(z) = \left( \BS - z\BI_p \right)^{-1}, ~~ z \in \mathbb{C}^+. \label{Eq_Def_Resolvent}
\end{equation}
We note that the definition of $z$ can be extended to points in the whole complex plane except for the eigenvalues of $\BS$. However, the analytic properties of the resolvent in $\mathbb{C}^+$ are generally sufficient for the study of the spectrum. The connection between the resolvent $\BQ(z)$ and the ESD $F^{\BS, n}$ is achieved by the Stieltjes transform, defined for the non-negative measure $\nu$ over $\R^+$ by
\begin{equation}
    m_{\nu}(z) := \int_{\R} \frac{\nu(\dd \lambda)}{\lambda - z}, ~~ z \in \mathbb{C}^+.
\end{equation}
Similar to the characteristic function, the Stieltjes transform is invertible in the sense that the measure $\nu$ can be recovered by the inversion formula \cite[Proposition 2.2]{hachem2007deterministic}. 
Then, by the definition \eqref{Eq_Def_Resolvent}, the Stieltjes transform of $F^{\BS, n}$ is given by 
\begin{equation}
    m_{\BS, n}(z) = \int_{\R^+} \frac{F^{\BS, n}(\dd \lambda)}{\lambda - z}= \frac{1}{p} \Tr  \BQ(z).
\end{equation}
Therefore, it suffices to establish the convergence of $\BQ(z)$ to study the convergence of ESD and the empirical measure outside the limiting support.
For ease of presentation, we will first introduce the fundamental fixed-point system of $2n$ equations that will be used to approximate the resolvent:
% \begin{subequations}
% \begin{align}
%   \delta_j(z) &= \frac{1}{n} \Tr \left[ \BOmega_j\left(-z \left(\BI_p + \sum_{l=1}^n \frac{\BOmega_l \widetilde{\delta}_l(z)}{n} \right) + \BA(\BI_n + \BD_{\delta}(z))^{-1}\BA^H \right)^{-1}\right], \\
%   \widetilde{\delta}_j(z) &= \left[ \left( -z \left(\BI_n + \BD_{\delta}(z) \right) + \BA^H \left(\BI_p + \sum_{l=1}^n \frac{\BOmega_l \widetilde{\delta}_l(z)}{n}\right)^{-1}\BA \right)^{-1} \right]_{j, j}, ~~ \forall j \in [n].
%   % \delta_l(z) &= \frac{\Tr \BOmega_l \BTheta(z)}{n},~~ \BD_{\delta}(z) = \diag\left(\delta_1(z), \ldots, \delta_n(z) \right), ~~ \widetilde{\delta}_l(z) = [\widetilde{\BTheta}(z)]_{l,l},   
% \end{align}
% \label{Eq_DE} 
% \end{subequations}
% where $\BD_{\delta}(z) = \diag(\delta_1(z), \delta_2(z), \ldots, \delta_n(z))$. 
% \begin{subequations}
% \begin{align}
%   \BTheta(z) &= \left[-z \left(\BI_p + \sum_{l=1}^n \frac{\BOmega_l \widetilde{\delta}_l(z)}{n} \right) + \BA(\BI_n + \BD_{\delta}(z))^{-1}\BA^H \right]^{-1}, \\
%   \widetilde{\BTheta}(z) &= \left[ -z \left(\BI_n + \BD_{\delta}(z) \right) + \BA^H \left(\BI_p + \sum_{l=1}^n \frac{\BOmega_l \widetilde{\delta}_l(z)}{n}\right)^{-1}\BA \right]^{-1},  \\
%   \delta_l(z) &= \frac{\Tr \BOmega_l \BTheta(z)}{n},~~ \BD_{\delta}(z) = \diag\left(\delta_1(z), \ldots, \delta_n(z) \right), ~~ \widetilde{\delta}_l(z) = [\widetilde{\BTheta}(z)]_{l,l},   
% \end{align}
% \label{Eq_DE} 
% \end{subequations}
\begin{equation}
    \left\{\begin{aligned}
        &\delta_i(z) = \frac{\Tr \BOmega_i \BTheta(z)}{n} , &  i \in [n], \\
        &\widetilde{\delta}_j(z) = \left[\widetilde{\BTheta}(z)\right]_{j, j}, & j \in [n],
    \end{aligned}\right. \label{Eq_DE}
\end{equation}
where
\begin{align*}
    \BF(z) &= \left[-z\left(\BI_p + \sum_{j=1}^n \frac{\BOmega_j\widetilde{\delta}_j(z)}{n} \right) \right]^{-1}, \\
    \widetilde{\BF}(z) &= \diag \left(\frac{-1}{z(1 + \delta_i(z))} ; 1 \leq i \leq n \right), \\
    \BTheta(z) &= \left( \BF^{-1}(z) - z \BA \widetilde{\BF}(z) \BA^H \right)^{-1}, \numberthis \label{Eq_DE_1}\\
    \widetilde{\BTheta}(z) &= \left( \widetilde{\BF}^{-1}(z) - z \BA^H {\BF}(z) \BA \right)^{-1}.  
\end{align*}
% \begin{subequations}
% \begin{align}
%   \BTheta(z) &= \left[-z \left(\BI_p + \sum_{l=1}^n \frac{\BOmega_l \widetilde{\delta}_l(z)}{n} \right) + \BA(\BI_n + \BD_{\delta}(z))^{-1}\BA^H \right]^{-1}, \\
%   \widetilde{\BTheta}(z) &= \left[ -z \left(\BI_n + \BD_{\delta}(z) \right) + \BA^H \left(\BI_p + \sum_{l=1}^n \frac{\BOmega_l \widetilde{\delta}_l(z)}{n}\right)^{-1}\BA \right]^{-1},  \\
%   \delta_l(z) &= \frac{\Tr \BOmega_l \BTheta(z)}{n},~~ \BD_{\delta}(z) = \diag\left(\delta_1(z), \ldots, \delta_n(z) \right), ~~ \widetilde{\delta}_l(z) = [\widetilde{\BTheta}(z)]_{l,l},   
% \end{align}
% \label{Eq_DE} 
% \end{subequations}
The details of the approximation rules are postponed to Theorem \ref{Thm_conver_ESD}. Define $\mathcal{S}$ as the set of Stieltjes transforms of non-negative measures over $\R^+$.
% \begin{equation}
%     \mathcal{S} = \left\{ m: m(z) = \int_{\R^+} \frac{\mu(\dd \lambda)}{\lambda - z}, z \in \mathbb{C}^+  \right\}
% \end{equation}
The analytic properties of the class $\mathcal{S}$ are provided in \cite[Proposition 2.2]{hachem2007deterministic}. The following proposition demonstrates the existence and uniqueness of the solutions to the aforementioned fixed-point equations in $\mathcal{S}$.
\begin{proposition}
\label{Prop_exists_unique}
    There exist unique solutions $(\delta_1, \ldots, \delta_n, \widetilde{\delta}_1, \ldots, \widetilde{\delta}_n)$ for \eqref{Eq_DE} with $\delta_j \in \mathcal{S}$ and $\widetilde{\delta}_j \in \mathcal{S}$, $\forall j \in [n]$.
\end{proposition}
\textit{Proof:} The proof of Proposition \ref{Prop_exists_unique} follows the same logic as that in \cite[Lemma 2.4]{hachem2007deterministic}. Here, we provide a proof sketch. The existence of the solution is established by using Picard’s approximation method. More specifically, we rewrite \eqref{Eq_DE} as $(\Bdelta, \widetilde{\Bdelta})$ = $f(\Bdelta, \widetilde{\Bdelta})$, where $\Bdelta = (\delta_1, ..., \delta_n)$ and $\widetilde{\Bdelta} = (\widetilde{\delta}_1, ..., \widetilde{\delta}_n)$. Then, we construct the initial point $(\Bdelta^0, \widetilde{\Bdelta}^0) \in \mathcal{S}^{2n}$, and then proceed inductively to show $(\Bdelta^{k+1}, \widetilde{\Bdelta}^{k+1}) = f(\Bdelta^{k}, \widetilde{\Bdelta}^{k}) \in \mathcal{S}^{2n}$ for each $k$. Furthermore, it can be shown that each element of $(\Bdelta^k, \widetilde{\Bdelta}^k)$ over a certain unbounded domain $\mathbb{D} \subset \mathbb{C}^+$ forms a Cauchy sequence for the sup norm, i.e., $\norm{f}_{\infty} = \sup_{z \in \mathbb{D}} |f(z)|$. Consequently, we denote  $(\Bdelta, \widetilde{\Bdelta}) = \lim_{k \to \infty}(\Bdelta^k, \widetilde{\Bdelta}^k)$ as the corresponding limit over $\mathbb{D}$. Since each element of $\Bdelta^k$ and $\widetilde{\Bdelta}^k$ is bounded on any compact subset of $\mathbb{C}^+$, the normal family theorem \cite{rudin1987real} and the properties of the Stieltjes transform \cite[Proposition 2.2]{hachem2007deterministic} imply that $(\Bdelta, \widetilde{\Bdelta}) \in \mathcal{S}^{2n}$. This establishes the existence of the solutions. The uniqueness assertion follows by assuming the existence of two distinct solutions, taking their difference, and deriving a contradiction. \qed
\begin{remark}
\label{Rem_iteration}
The proof demonstrates that the numerical solution for 
$\delta_j(z), \widetilde{\delta}_j(z)$ with a given $z \in \mathbb{C}^+$ can be obtained by using the classical fixed-point iteration algorithm. We note that for any $m \in \mathcal{S}$, the corresponding measure $\mu$ is supported on $\mathbb{R}^+$. Therefore, the Stieltjes transform $m(z)$ can be analytically extended to $z \in (-\infty, 0)$. Given the solutions $(\Bdelta, \widetilde{\Bdelta}) \in \mathcal{S}^{2n}$, the same iterative algorithm can be used to evaluate the values of $(\Bdelta(z), \widetilde{\Bdelta}(z))$ over $z \in (-\infty, 0)$.
\end{remark}
In the following, when we refer to the fixed-point system of equations \eqref{Eq_DE} and \eqref{Eq_DE_1}, we assume the unique solution described in Proposition \ref{Prop_exists_unique} without further distinction.
\par
\subsection{Weak Convergence of $F^{\BS, n}$}
It is well established that the vague convergence \cite{bai2010spectral} of a family of measures is equivalent to the pointwise convergence of their corresponding Stieltjes transforms on $\mathbb{C}^+$. The following theorem demonstrates the deterministic approximations for the trace and bilinear forms of $\BQ(z)$.
\begin{theorem}
\label{Thm_conver_ESD}
Assume Assumptions \ref{Assumpt_1}-\ref{Assumpt_4} hold. Then, for given $z \in \mathbb{C^+}$, we have, as $n \to \infty$, 
\begin{equation}
    \frac{1}{p} \Tr\BC\left[\BQ(z) - \BTheta(z)\right]\xrightarrow{a.s.} 0, ~~  \Bu^H \left[\BQ(z) - \BTheta(z)\right] \Bv \xrightarrow{a.s.} 0, \label{Eq_Thm1_Conv_Resolvent}
\end{equation}
where $\BC \in \mathbb{C}^{p \times p}$ is a deterministic matrix with bounded spectral norm, i.e., $\limsup_n \norm{\BC} < \infty$, and $\Bu \in \mathbb{C}^p$ and $\Bv \in \mathbb{C}^p$  are deterministic vectors with bounded Euclidean norm, i.e., $\limsup_n\max\{ \norm{\Bu}, \norm{\Bv} \}  < \infty$.
\end{theorem}
\textit{Proof:} The proof of Theorem \ref{Thm_conver_ESD} is given in Appendix \ref{App_Converge_resolvent}. \qed
\par
Denote $m_n(z) = \Tr \BTheta(z) / p$. It can be verified that $m_n(z) \in \mathcal{S}$ with the integration representation $m_n(z) = \int_{\mathbb{R}^+} \frac{F^n(\dd \lambda)}{\lambda - z}$. By letting $\BC = \BI_p$ in Theorem \ref{Thm_conver_ESD}, we have $m_{\BS, n}(z) - m_n(z) \to 0$ almost surely. Due to the complex correlation structure and the existence of the mean, the Stieltjes transform $m_n(z)$ may not converge as $n \to \infty$. As a result, the limiting distribution may not exist. However, we have the following ``convergence in distribution" property. Specifically, following Lemma \ref{Lemm_Bound_on_qterms} and the same method as in the proof of \cite[Corollary 2.7]{hachem2007deterministic}, we have 
\begin{equation}
    \int_{0}^{\infty} f(\lambda) F^{\BS, n}(\dd \lambda)- \int_{0}^{\infty} f(\lambda) F^n(\dd \lambda)\xrightarrow{a.s.} 0, \label{Eq_ESD_LSD_convergence}
\end{equation}
for any bounded continuous function $f: \mathbb{R}^+ \to \mathbb{R}$.  Without causing any ambiguity, we call $F^n$ the LSD in this work. 
\begin{remark}
To prove the convergence of \eqref{Eq_ESD_LSD_convergence}, only the convergence of the trace form for the resolvent is needed. Write the eigen-decomposition of $\BS$ as $\BU_S \diag({\lambda_i ; 1 \leq i \leq p}) \BU_S^H$ and denote $ f(\BS) = \BU_S \diag({f(\lambda_i) ;1 \leq  i \leq p}) \BU_S^H$. \eqref{Eq_ESD_LSD_convergence} implies that
$
   \Tr f(\BS) / p = \sum_{j=1}^p [f(\BS)]_{j, j} / p
$
, i.e., the average of the diagonal elements of $f(\BS)$, will converge. However, this result does not imply the convergence of the specific elements $[f(\BS)]_{i, j}$, which requires the convergence of the bilinear form of the resolvent \cite{nagel2021functional}. The bilinear forms of the resolvent are also useful for the study of the eigenvectors for the covariance matrix $\BS$ \cite{hachem2013bilinear, BAIK2006spiked, couillet_liao_2022, Xia2013AoS}. 
\end{remark}
\begin{remark}
    The convergence in \eqref{Eq_Thm1_Conv_Resolvent} also holds for $z \in (-\infty, 0)$.
\end{remark}
 Next, we show that under specific conditions on the statistics of $\BA$ and ${\BOmega_j}$, the limit of $F^n$ exists. The details are provided in the following corollary.
\begin{corollary}
\label{Coro_variance profile}
 Assume $\lim_{n \to \infty} p / n = c \in (0, + \infty)$ and Assumption \ref{Assumpt_2} holds. Assume that $\{ \BOmega_j\}_{j \leq n}$ are simultaneously diagonalizable with the eigen-decomposition $\BOmega_j = \BU \mathrm{diag}(\lambda_{jl}(n); l \leq p) \BU^H$, and $\BA$ has the same left singular vectors as $\BOmega_j$ with $\BA = \BU \mathrm{rect}\text{-}\mathrm{diag}(\sigma_j(n); j\leq \min(p, n)) \BI_n$. Further, assume the eigenvalues $\{ \lambda_{jl}(n)\}_{l \leq p, j \leq n}$ and the singular values $\{ \sigma_j(n)\}_{j \leq \min(p, n)}$ can be parameterized by 
 \begin{equation}
    \lambda_{jl}(n) = f_{\Omega}(l / p, j / n), ~~ \sigma_j(n) = f_A(j / p)  = \widetilde{f}_A(j / p) \Ind_{\{ j / p \leq \min(1, c^{-1}) \}},
\end{equation}
where $f_{\Omega} > 0$ and $\widetilde{f}_A$ are  bounded continuous functions. Then, the functional system of equations
\begin{align}
    \tau(a, z) &= \left[ {-z\left(1 +  \int_{0}^1 f_{\Omega}(a, u) \widetilde{\tau}(u, z) \dd u \right) + \frac{f^2_A(a)}{1 + c\int_{0}^1 f_{\Omega}(x, a)\tau(x, z) \dd x}} \right]^{-1} , \\
    \widetilde{\tau}(b, z) &= \left[ -z\left(1 + c\int_{0}^1f_{\Omega}(y, b) \tau(y, z) \dd y\right) + \frac{f^2_A(b / c)}{1 +   \int_{0}^1 f_{\Omega}(b, v) \widetilde{\tau}(v, z) \dd v} \right]^{-1},
\end{align}
admits unique function solutions $\Lambda: [0, 1] \times \mathbb{C}^+ \to \C^+$  such that $u \mapsto \Lambda(u, z)$ is continuous and $z \mapsto \Lambda(u, z) \in \mathcal{S}$, where $\Lambda$ is any of $\tau$ and $\widetilde{\tau}$. The function $m_{\infty}(z) = \int_{0}^1 \tau(x, z) \dd x$ is the Stieltjes transform of a probability distribution $F^{\infty}$ and there holds
\begin{equation}
    \PP \left( F^{\BS, n} \Rightarrow F^{\infty} \right) = 1.
\end{equation}
\end{corollary}
\textit{Proof:} The proof is a direct application of Theorem \ref{Thm_conver_ESD} and is omitted here for brevity. \qed
\begin{remark}
    When $\{ \BOmega_j \}_{j \leq n}$ are diagonal and $\BA$ is rectangular diagonal, i.e., $\BU = \BI_p$, Corollary \ref{Coro_variance profile} degenerates to \cite[Theorem 2.3]{hachem2006varianceprofile} and \cite[Theorem 2.5]{hachem2008varianceProfile}. The degeneration indicates that the ESD and the LSD are rotation-invariant.
\end{remark}
% Now we give some remarks in the following.
% \begin{remark}
%     The convergence result also holds for $z \in \mathbb{C}^- \cup \mathbb{R}^+$ if we extend the definition of Stieltjes transform $m(z)$ to these domains.
% \end{remark}
% \begin{remark}
% \label{Rm_ESD_weak_convergece}
%     In the above discussion, $F^n$ is actually also a measure sequence. The convergence in distribution here is in the following sense 
%     \begin{equation}
%         \int_{0}^{\infty} f(\lambda) (F^{\BS, n}- F^n)(\dd \lambda) \xrightarrow{a.s.} 0,
%     \end{equation}
%     for any bounded continuous function $f: \mathbb{R}^+ \to \mathbb{R}$.
% \end{remark}
\subsection{No Eigenvalues outside the Support}
The weak convergence in \eqref{Eq_ESD_LSD_convergence} shows that for any interval $[a, b]$, the empirical measure $F^{\BS, n}([a, b])$ converges almost surely to $ F^n([a, b])$. By the definition in \eqref{EQ_ESD_F_S_n}, we have $F^{\BS, n}([a, b]) = (\text{the number of eigenvalues of $\BS$ in } [a, b]) / p$. If we take $[a, b]$ such that $F^n([a, b]) = 0$, Theorem \ref{Thm_conver_ESD} indicates that $F^{\BS, n}([a, b]) = o_{a.s.}(1)$. Intuitively, the weak convergence of $F^{\BS, n}$ will not be affected if finite empirical eigenvalues appear in $[a, b]$. In this section, we will prove that such appearance will not happen, i.e., the order of $F^{\BS, n}([a, b])$ is $o_{a.s.}(1/n)$.  To this end, we need more analytical properties of the solutions $\{ \delta_j \}, \{ \widetilde{\delta}_j \}$, which are given in the following assumption.
\begin{assumption}
\label{Assumpt_5}
    For all $x \in \mathbb{R}_*$, the limits $\lim_{z \in \mathbb{C}^+ \to x} m_n(z)$, $\lim_{z \in \mathbb{C}^+ \to x} \delta_j(z)$, and $\lim_{z \in \mathbb{C}^+ \to x} \widetilde{\delta}_j(z)$ exist for each $j \in [n]$.
\end{assumption}
% \begin{assumption}
% \label{Assumpt_5}
%     For all $x \in \mathbb{R} - \{0\}$, the limits $\lim_{z \in \mathbb{C}^+ \to x} m(z)$, $\lim_{z \in \mathbb{C}^+ \to x} \delta_j(z)$ and $\lim_{z \in \mathbb{C}^+ \to x} \widetilde{\delta}_j(z)$ exist for each $j \in [n]$.
% \end{assumption}
\begin{remark}
    This assumption holds in many cases. For instance, in \cite{silverstein1995analysis}, Silverstein \textit{et al.} proved the existence of these limits when $\BA = \Bzero_{p \times n}$ and $\BB_j = \BT^{1/2}$, $\forall j \in [n]$. In \cite{couillet2014analysis}, Couillet et al. performed the same study for the separable correlated random matrix model with $\BA = \Bzero_{p \times n}$, and $\BB_j = \sqrt{{\tau}_j} \BT^{1/2} $ in \eqref{Eq_DE}.   In \cite{Dozier2007analysis}, Dozier et al. established the case for the information-plus-noise model with $\BOmega_j = \sigma^2 \BI_p$. In \cite{zhou2024analysis}, Zhou \textit{et al.} proved the existence of the limits for the information-plus-noise model with $\BB_j = \BT^{1/2}$, $\forall j \in [n]$, under the assumption
    that matrices $\BT$ and $\BA\BA^H$ are commutative. We make this assumption here for technical convenience. We note that this assumption is possibly provable and leave the proof as a future study.
\end{remark}
For notational simplicity, we define the limits in Assumption \ref{Assumpt_5} as $m_n(x)$, $\delta_j(x)$, and $\widetilde{\delta}_j(x)$, respectively.
Denote the underline measure of $\delta_j(z)$ and $\widetilde{\delta}_j(z)$ as $\mu_j$ and $\widetilde{\mu}_j$, $j \in [n]$, respectively. By \cite[Lemma B.10]{bai2010spectral}, the distribution function $\mu_j([0, x])$, $\widetilde{\mu}_j([0, x])$, and $F^n(x)$ are differentiable at $x$ when Assumption \ref{Assumpt_5} holds, and their derivatives are given by $\frac{1}{\pi} \Im (m_n(x))$, $\frac{1}{\pi} \Im (\delta_j(x))$, and $\frac{1}{\pi}\Im(\widetilde{\delta}_j(x))$, respectively.  The following proposition reveals the relationship between the supports of $\mu_j$, $\widetilde{\mu}_j$, and $F^n$, which are necessary for the proof of Theorem \ref{Thm_No_Eigenvalues}.
\begin{proposition}
\label{Prop_same_support}
    Assume Assumptions \ref{Assumpt_1}, \ref{Assumpt_3}-\ref{Assumpt_5} hold. Then, the supports of the measures $\{ \mu_j \}_{1 \leq j \leq n}$ and $\{ \widetilde{\mu}_j \}_{1 \leq j\leq n}$ are subsets of the support of $F^n$ on $\R_*$, that is
    \begin{equation}
        \left( \mathsf{Supp}(\mu_j ) \cap \R_* \right) \subset \left( \mathsf{Supp}(F^n)   \cap \R_*\right), ~~ \left( \mathsf{Supp}(\widetilde{\mu}_j) \cap \R_* \right) \subset \left( \mathsf{Supp}(F^n)   \cap \R_* \right), ~~ \forall j \in [n].
    \end{equation}
\end{proposition}
\textit{Proof:} The proof is given in Appendix \ref{App_Proof_Prop_same_support}. \qed

\begin{remark} 
    An interesting case is when the supports of these measures are identical. In fact, we can release the left hand side (LHS) of the inequality \eqref{Eq_Assumpt_3} in Assumption \ref{Assumpt_3} to
      \begin{equation}
    0 < \liminf_{n \to \infty} \min_{1 \leq j \leq n} \lambda_{\min}(\BOmega_j), \label{Eq_Enhanced_Eq_3}
    \end{equation}  
    where $\lambda_{\min}(\BOmega)$ denotes the minimal eigenvalue of $\BOmega$. With the above condition, it holds that
    \begin{equation}
         \left( \mathsf{Supp}(F^n)   \cap \R_*\right) = \left( \mathsf{Supp}(\mu_j ) \cap \R_* \right) = \left( \mathsf{Supp}(\widetilde{\mu}_j)   \cap \R_* \right), ~~ \forall j \in [n].
    \end{equation}
\end{remark}
The following theorem shows the no-eigenvalue property for the generally correlated and noncentral covariance matrices. 
\begin{theorem}
\label{Thm_No_Eigenvalues}
    Assume Assumptions \ref{Assumpt_1}-\ref{Assumpt_5} hold. Let $[a, b]$ with $a > 0$ lie in an open interval outside $ \mathsf{Supp}(F^n)$ for all large $n$ and satisfy the condition $\liminf_{n \to \infty} \inf_{x \in [a, b]} \min_{1 \leq j \leq n} \big\{ |\widetilde{\delta}_j(x)|, |1 + \delta_j(x)| \big\} > 0$.  Then, we have 
    \begin{equation}
    \label{Eq_P_no_eigens}
        \PP \left[ \text{no eigenvalues of $\mathbf{S}$ appears in $[a, b]$ for all large $n$}  \right] = 1.
    \end{equation}
\end{theorem}
Theorem \ref{Thm_No_Eigenvalues} demonstrates that there exists a sufficiently large positive integer $N$ such that for any $n > N$, the eigenvalues of  $\BS$ cannot be located in the interval $[a, b]$.
\par
\textit{Proof:} Since the proof of Theorem \ref{Thm_No_Eigenvalues} is quite lengthy and involves many steps, we provide the main flow of the proof as follows:
\begin{itemize}
    \item []\textit{Step 1:}  Simplifying Assumption \ref{Assumpt_2} to get the equivalent model. In particular, without changing the asymptotic locations of the eigenvalues of $\BS$, it suffices to truncate the random variables $\{X_{ij}\}$. The details of this step are given in Appendix \ref{App_Trunc}.
    \item []\textit{Step 2:} With the simplified model, we evaluate the convergence of the resolvent. Specifically, we will show that for $z = x + j v$ with $v = n^{-\alpha / 912}$, and $\alpha = \frac{\varepsilon}{8 + 2\varepsilon}$ (The coefficient 912 comes from the constraints on the speed at which $z$ approaches the real axis in the proof, which may not be minimal), the following holds 
    % \begin{equation}
    % \label{Eq_sup_diff_resolvent}
    %     \sup_{x \in [a, b]} v \abs{\Tr \BQ(z) - \Tr \BTheta(z)} \xrightarrow[n \to \infty]{a.s.} 0, 
    % \end{equation}
    \begin{equation}
    \label{Eq_sup_diff_resolvent}
        \sup_{x \in [a, b]} vp \abs{m_{\BS, n}(z) - m_n(z)} \xrightarrow[n \to \infty]{a.s.} 0. 
    \end{equation}
     The details of the proof are given in Appendix \ref{App_Proof_Thm_noeigen}.
    \item []\textit{Step 3:} In the last step, we will prove \eqref{Eq_P_no_eigens}  by \eqref{Eq_sup_diff_resolvent}.
\end{itemize}
Next, we show \textit{Step 3}, which is similar to the procedure in \cite[Section 6]{bai1998no}. Denote $\beta = \alpha / 912$. Since $p = \frac{p}{n} n = \frac{p}{n} v^{-\frac{1}{\beta}}$ and $\liminf \frac{p}{n} > 0$, equation \eqref{Eq_sup_diff_resolvent} implies 
\begin{equation}
    \max_{k \in [N_0]}v^{1 - \frac{1}{\beta}} \sup_{x \in [a, b]} \abs{m_{\BS, n}(x + \jmath \sqrt{k}v) - m_{n}(x + \jmath \sqrt{k}v)} \xrightarrow[n \to \infty]{a.s.} 0,
\end{equation}
where $N_0$ is an integer such that $N_0 \geq 1/{(2\beta)}$. Taking the imaginary part of $m_{\BS, n}(x + \jmath \sqrt{k}v) - m_{n}(x + \jmath \sqrt{k}v)$, we have 
\begin{equation}
    \max_{k \in [N_0]}v^{2 - \frac{1}{\beta}} \sup_{x \in [a, b]} \abs{\int\frac{ (F^{\BS, n} - F^n)(\dd \lambda)}{(x - \lambda)^2 + kv^2}} \xrightarrow[n \to \infty]{a.s.} 0.
\end{equation}
Denote $L(k) = \int\frac{ (F^{\BS, n} - F^n)(\dd \lambda)}{(x - \lambda)^2 + kv^2}$. By taking the difference $L(k_1) - L(k_2)$, we can get
\begin{equation}
    \max_{k_1 \neq k_2 }v^{4 - \frac{1}{\beta}} \sup_{x \in [a, b]} \abs{\int\frac{ (F^{\BS, n}- F^n)(\dd \lambda)}{\left[(x - \lambda)^2 + k_1v^2\right]\left[(x - \lambda)^2 + k_2v^2\right]}} \xrightarrow[n \to \infty]{a.s.} 0.
\end{equation}
 By repeating the differencing process, we can obtain
 \begin{equation}
    v^{{2N_0} - \frac{1}{\beta}} \sup_{x \in [a, b]} \abs{\int\frac{ (F^{\BS, n}- F^n)(\dd \lambda)}{\prod_{k=1}^{N_0}\left[(x - \lambda)^2 + kv^2\right]}} \xrightarrow[n \to \infty]{a.s.} 0.
\end{equation}
Since $[a, b]$ lies in an open interval outside $\mathsf{Supp}(F^n)$, there exist $a' < a$ and $b' > b$ such that $[a', b']$ is outside $\mathsf{Supp}(F^n)$ for all large $n$. Splitting the integration interval and using $p = v^{-1/\beta}$, we have 
 \begin{equation}
     \sup_{x \in [a, b]} \abs{\int\frac{ v^{{2N_0} - \frac{1}{\beta}}\Ind_{[a', b']^c} (F^{\BS, n}- F^n)(\dd \lambda)}{\prod_{k=1}^{N_0}\left[(x - \lambda)^2 + kv^2\right]} + \sum_{\lambda_j \in [a', b']} \frac{v^{2N_0}}{\prod_{k=1}^{N_0}\left[(x - \lambda_j)^2 + kv^2\right]}} \xrightarrow[n \to \infty]{a.s.} 0. \label{Eq_sup_ab_Diff_FS_Fn}
\end{equation}
We write the LHS of \eqref{Eq_sup_ab_Diff_FS_Fn} as $\sup_{x \in [a, b]}\abs{A(x) + B(x)}$. By the weak convergence in \eqref{Eq_ESD_LSD_convergence}, we have $\sup_{x \in [a, b]} |A(x)| \to 0$ almost surely. If there is at least one eigenvalue in the interval $[a, b]$, then $B(x)$ is uniformly bounded away from 0, which leads to a contradiction. Therefore, \eqref{Eq_P_no_eigens} is proved. \qed
\begin{figure}[t]
    \centering
    \includegraphics[width=4.4in]{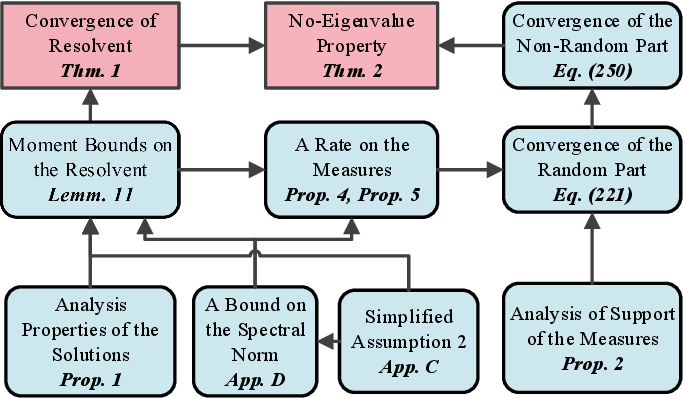}
    \caption{Logic Flow for the Proofs of Theorem \ref{Thm_conver_ESD} and Theorem \ref{Thm_No_Eigenvalues}.}
    \label{Fig_Logic_flow}
\end{figure}
\begin{remark}
The condition $\liminf_{n \to \infty} \inf_{x \in [a, b]} \min_{j \leq n} \big\{ |\widetilde{\delta}_j(x)|, |1 + \delta_j(x)| \big\} > 0$ is not necessary in many cases. In \cite{Bai2025noeigen}, Bai \textit{et al.} proved that this condition holds uniformly for $x > 0$ for the noncentral model with $\BOmega_j = \BT$, where matrices $\BT$ and $\BA\BA^H$ are commutative. Additionally, if $\BA = \mathbf{0}_{p \times n}$, it can also be shown that this condition is satisfied uniformly.  Intuitively, $\{ \delta_j(z) \}_{j \leq n}$ and $\{ \widetilde{\delta}_j(z) \}_{j \leq n}$ are continuous functions. When $x$ is larger than the right endpoint of $\mathsf{Supp}(F^n)$, we have $\Re(\delta_j(z)), \Re(\widetilde{\delta}_j(z)) < 0$, $\lim_{x \to \infty} \Re(\delta_j(z)) = 0$, and $\lim_{x \to \infty} \Re(\widetilde{\delta}_j(z)) = 0$. Similarly, when $x$ is smaller than the left endpoint of $\mathsf{Supp}(F^n)$, we have $|\Re (\delta_j(z))|, |\Re(\widetilde{\delta}_j(z))| > 0$. Therefore, the number of zeros for $|1 + \delta_j(z)| = 0$ and $\widetilde{\delta}_j(z) = 0$ is limited. Furthermore, this condition can be verified numerically by Remark \ref{Rem_iteration}. 
\end{remark}
\begin{remark}
    In the fields of signal detection and statistical inference, many subspace-based algorithms involve integrals of the resolvent in the following form  \cite{couillet_liao_2022, hachem2013bilinear}
    \begin{equation}
        \int_{\mathcal{C}^+} f(z) \digamma[\BQ(z)] \dd z, \label{Eq_integral_Resolvent}
    \end{equation}
    where $\digamma(\cdot)$ is a linear operator that can have the form of $ \frac{1}{p} \Tr \BC (\cdot) $ or $ \Bu^H(\cdot) \Bv $, and $\mathcal{C}^{+}$ is typically the pre-selected positive contour that includes the eigenvalues corresponding to the signal subspace. Theorem \ref{Thm_No_Eigenvalues} provides theoretical guarantees for these algorithms that the contour is almost surely valid. Moreover, a surprising phenomenon in random matrix theory is that the central limit theorem (CLT) for LSS holds \cite{Bai2004CLTaop}. In the key step of the proof in \cite{Bai2004CLTaop}, \eqref{Eq_integral_Resolvent} was used, where $\mathcal{C}^+$ was chosen as the contour that includes the entire support. Theorem \ref{Thm_No_Eigenvalues} also provides a theoretical foundation for the analysis of the LSS for generally correlated and noncentral random matrices.
\end{remark}
A direct application of Theorem \ref{Thm_No_Eigenvalues} is bounding the spectral norm of the sample covariance $\BS$, which is stated in the following corollary.
\begin{corollary}
\label{Coro_largest_eigenvalue}
    Denote the right endpoint of $\mathsf{Supp}(F^n) \cap \R_* $ as $e_n^+$. If 
    $
     \limsup_{n \to \infty} e_n^+ < \infty, 
    $
    there exists $E^+ <\infty$ such that
    \begin{equation}
        \PP \left( \limsup_{n \to \infty} \lambda_{\max}(\BS) \leq E^+\right) = 1. \label{Eq_Bound_on_Sp_Norm}
    \end{equation}
    % If $\liminf_{n \to \infty} e_n^- > \infty$, the following holds
    % \begin{equation}
    %     \PP \left( \liminf_{n \to \infty} \lambda_{\min}(\BS) \geq E^- \right) = 1. \label{Eq_Bound_on_Min_eigenvalue}
    % \end{equation}
\end{corollary}
\textit{Proof:} Choose $[a, b]$ such that $a$ is large enough with  $a - e_n^+ \geq  2 \sup_{n, j} \frac{p}{n} \norm{\BOmega_j}$. Then, for $x \in [a, b]$, we have 
\begin{equation}
    \delta_j(x) \overset{(a)}{=} \int_{\R^+}\frac{\mu_j(\dd\lambda)}{\lambda - x} \geq - \frac{\mu_j(\R^+)}{2 \sup_{n, j}\frac{p}{n} \norm{\BOmega_j}}, ~~~ \widetilde{\delta_j}(x) = \int_{\R^+}\frac{\widetilde{\mu}_j(\dd \lambda)}{\lambda - x} \leq  -\frac{\widetilde{\mu}_j(\R^+)}{b}.
\end{equation}
 Here, step $(a)$ follows from the analytic continuation argument and the fact $ \delta_j(z) $ is analytic over $ \mathbb{C}^+$. By the properties of the Stieltjes transform \cite[Proposition 2.2]{hachem2007deterministic}, we have $\mu_j(\R^+) = \frac{1}{n} \Tr \BOmega_j$ and $\widetilde{\mu}_j(\R^+) = 1$, which validate the condition  $\liminf_{n \to \infty} \inf_{x \in [a, b]} \min_{j \leq n} \big\{ |\widetilde{\delta}_j(x)|, |1 + \delta_j(x)| \big\} > 0$.  Due to the arbitrariness of $b$, \eqref{Eq_Bound_on_Sp_Norm} is proved.  \qed
 % The a.s. lower boundedness of $\lambda_{\min}(\BS)$ can be shown in a similar manner. \qed
\begin{remark}
     In \cite{yin1988limit}, it was demonstrated that when $\BA = \Bzero_{p \times n}$, $ \BOmega_j = \BI_p$, and $ \E |X_{11}|^4 < \infty $, the largest eigenvalue of $\BS$ converges to  $(1 + \sqrt{y})^2$ almost surely as $p/n \to y $. The proof in \cite{yin1988limit} first identified the statistical equivalent model, followed by employing the combinatorial method to evaluate the largest eigenvalue by the inequality $ \E \lambda_{\max}^k(\BS)\leq \E\Tr(\BS^k)$. In \cite{BAI19884th}, it was shown that the $4$-th moment condition is the weakest that ensures the largest eigenvalue to be bounded almost surely. However, to prove the boundedness of the largest eigenvalue of the model in \eqref{Eq_random_matrix_Model}, the combinatorial method is not applicable since $\BOmega_j$s can be non-diagonal, leading to a significantly more complex expansion of $ \Tr(\BS^k)$.
\end{remark}
\begin{remark}
In \cite{WagnerMISO2012}, the authors assumed the almost sure boundedness for the spectral norm of $\BS$ to prove the convergence of the resolvent for the covariance of the generally correlated and centered random matrix (refer to Assumption 3 in \cite{WagnerMISO2012}). Corollary \ref{Coro_largest_eigenvalue} specifies the conditions under which this boundedness is valid.
\end{remark}
To clearly illustrate the relationships between the key results in this paper and enhance the readability, we summarize the logical flow for the proofs of Theorems \ref{Thm_conver_ESD} and \ref{Thm_No_Eigenvalues} in Fig. \ref{Fig_Logic_flow}. The red blocks represent the main results of this paper, while the blue blocks indicate intermediate results.

% \subsection{Relation to Previous Works}
\section{Applications}
\label{Sec_applications}
In this section, we apply the main theoretical results to two communication scenarios.
\subsection{SINR Analysis of LMMSE Detection}
\label{Sec_LMMSE_dec}
% \subsection{Multi-User MIMO Detection with Rician Channels}
 Consider a multi-user uplink MIMO system with a base station (BS) equipped with $p$ antennas and $n + 1$ single-antenna users. Denote $x_j \sim \mathcal{CN}(0, 1)$ as the symbol transmitted by the $j$-th user. Then, the received signal $\By \in \mathbb{C}^p$ at the BS is given by
 \begin{equation}
     \By = \sum_{j=0}^n \Bh_j x_j + \Bn,
 \end{equation}
where $\mathbf{n} \sim \mathcal{CN}(0, \sigma^2 \mathbf{I}_p)$ represents the additive white Gaussian noise (AWGN) at the BS and $\Bh_j \in \mathbb{C}^p $ denotes the channel between the $j$-th user and the BS. Consider the Rician channel model \cite{Kammoun2019RZF}
\begin{equation}
    \Bh_j = \sqrt{\frac{1}{1 + \tau}} \BC_j^{1/2} \Bz_j + \sqrt{\frac{\tau}{1 + \tau}} \overline{\Bz}_j, ~~ j \in [n] \cup \{0 \},
    \label{Eq_Rician_model}
\end{equation}
where $\Bz_j \sim \mathcal{CN}(\mathbf{0}_p, \BI_p / n)$ and $\BC_j \in \mathbb{C}^{p \times p}$ is a Hermitian non-negative matrix that models the spatial correlation of the BS antennas with respective to the $j$-th user. $\overline{\Bz}_j \in \mathbb{C}^p$ is a deterministic vector that represents the line-of-sight (LoS) component. The scalar $\tau \geq 0$ is the Rician factor. With perfect channel state information at the receiver, the classical LMMSE estimator for symbol $x_0$ is given by
\begin{equation}
    \widehat{x}_0 = \Bh_0^H\left(\BH\BH^H+ \sigma^2 \BI_p \right)^{-1} \By,
\end{equation}
where $\BH = [\Bh_0, \ldots, \Bh_n]$. Rewrite the received vector as $\By =\Bh_0 x_0 + \By_{0}$ where $\By_{0} = \sum_{i \neq 0} \Bh_i x_i + \Bn$ represents the interference-plus-noise term. Thus, the SINR of the $0$-th user with the LMMSE receiver is given by
\begin{equation}
    \gamma_0 = \frac{ \E_{\Bx, \Bn}|\Bg^H\Bh_0 x_0|^2}{\E_{\Bx, \Bn}|\Bg^H \By_{0} |} = \Bh_0^H\left(\BH_{[0]} \BH_{[0]}^H + \sigma^2 \BI_p \right)^{-1} \Bh_0, \label{Eq_SINR_LMMSE}
\end{equation}
where $\Bg = \left(\BH\BH^H+ \sigma^2 \BI_p \right)^{-1}\Bh_0$  and $\BH_{[0]} = [\Bh_1, \ldots,  \Bh_n]$ is the rank-one perturbation of $\BH$ by removing the first column. As the number of antennas $p$ and the number of users $n$ approach infinity with the same speed, the asymptotic behavior of $\gamma_0$ is interesting, and many results over different channel conditions have been obtained in the literature \cite{Abla2009Ber, Kammoun2009clt, Zhuang2025Decentralized}. However, the asymptotic behavior of $\gamma_0$ over Rician fading channels with heterogeneous spatial correlation is not yet available. The following corollary elucidates this asymptotic behavior.
\begin{corollary} 
\label{Coro_SINR}
Denote $\overline{\BZ} = [\overline{\Bz}_1, \ldots. \overline{\Bz}_n]$. Assume the dimensions $p$ and $n$ grow to infinity with
\begin{itemize}
    \item [1)]  $0 < \liminf_n p / n \leq  \limsup_n p / n < \infty$. 
    \item [2)] $ 0 < \liminf_n \min_{0 \leq j \leq n} \frac{1}{p} \Tr \BC_j \leq \limsup_n \max_{0 \leq j \leq n} \norm{\BC_j} < \infty$.
    \item [3)] $\limsup_n \norm{\left[\overline{\Bz}_0, \overline{\BZ}\right]} < \infty$.
\end{itemize}
Then, we have 
\begin{equation}
    \gamma_0 - \left( \frac{\Tr \BC_0 \BTheta }{n(1 + \tau)} +  \frac{ \tau \overline{\Bz}_0^H \BTheta \overline{\Bz}_0}{1 + \tau}\right) \xrightarrow[n \to \infty]{a.s.} 0,
\end{equation}
where $\BTheta$ is obtained by setting $\BOmega_j = p_j\BC_j / (1 + \tau)$ and $\BA = \overline{\BZ}  \sqrt{\tau}/ \sqrt{1 + \tau}$ in \eqref{Eq_DE}.
\end{corollary}
\textit{Proof:} The proof can be derived from Theorem \ref{Thm_conver_ESD} and follows a similar approach as in \cite{Zhuang2025Decentralized}. Therefore, we omit it here. \qed 
\subsection{The Invertibility of the ZF Precoding Matrix}
Consider a multi-user downlink MIMO system where a BS with $p$ antennas serves $n$ single-antenna users. Denote $\Bg_i \in \mathbb{C}^p$ as the precoding vector associated with the $i$-th user. Then, the received signal $y_j$ of the $j$-th user is given by
\begin{equation}
    y_j = \Bh_j^H \sum_{i=1}^n \Bg_i x_i + n_j, ~~ j \in [n],
\end{equation}
where $x_j \sim \mathcal{CN}(0, 1)$ and $n_j \sim \mathcal{CN}(0 , \sigma^2)$ denote the data symbol intended to the $j$-th user and the AWGN, respectively. $\Bh_j \in\mathbb{C}^p$ denotes the random channel vector from the BS to the $j$-th user and is modeled as \cite{WagnerMISO2012}
\begin{equation}
    \Bh_j = \BC^{1/2}_j \Bz_j ,~~ j \in [n], 
    \label{Eq_Channel_model_rayleigh}
\end{equation}
where $\BC_j \in \mathbb{C}^{p \times p}$ is Hermitian non-negative and $\Bz_j$ is a random vector (not necessarily Gaussian distributed), whose elements are i.i.d. with $\mathbb{E}[\Bz_j]_i = 0$ and $\mathbb{E}|\sqrt{n}[\Bz_j]_i|^2 = 1 $. Furthermore, assume there exists a constant $e > 0 $ such that $\mathbb{E}|\sqrt{n}[\Bz_j]_{i}|^{4+e} < \infty$. To mitigate the inter-user interference, a classic design for the precoding matrix $ \BG = [\Bg_1, \ldots, \Bg_n] $ is ZF.  To this end, we utilize the same precoding matrix as  \cite{AblaNoEigen}, which is given by
\begin{equation}
    \BG = \left( \BH \BH^H\right)^{-1} \BH,
\end{equation}
with $\BH = [\Bh_1, \ldots, \Bh_n]$.
The performance of the ZF precoding has been extensively studied as the numbers of antennas and users increase simultaneously \cite{Hoydis2013MIMO}. However, to analyze the performance of the ZF precoder, it is essential to determine whether the Gram matrix $\BH\BH^H$ is singular. The following corollary provides a sufficient condition for the well-definedness of the precoding matrix $\BG$.
\begin{corollary}
\label{Coro_min_eig_zero_forcing}
Assume the dimensions $p$ and $n$ grow to infinity with 
\begin{itemize}
\item [1)]  $0 < \liminf_n p / n \leq  \limsup_n p / n < 1$. 
\item [2)] $ 0 < \liminf_n \min_{1 \leq j \leq n} \lambda_{\min}(\BC_j) \leq \limsup_n \max_{1 \leq j \leq n} \norm{\BC_j} < \infty$.
\end{itemize}
Then, there exists a  constant $E_-> 0$ such that
\begin{equation}
    \PP \left( \liminf_{n \to \infty} \lambda_{\min}(\BH\BH^H) \geq E_- \right) = 1.
\end{equation}
\end{corollary}
\textit{Proof:} The proof can be obtained by applying Theorem \ref{Thm_No_Eigenvalues}. For that purpose, we first extend the result to the case with $a = 0$, and the details are given in Appendix \ref{App_prof_of_coro_min_eig_zero_forcing}. \qed
\begin{remark}
In \cite{AblaNoEigen}, it was proved that when $\{ \Bz_j \}_{j \leq n}$ defined in \eqref{Eq_Channel_model_rayleigh} are complex Gaussian random vectors, the minimal eigenvalues of $\BH\BH^H$ are bounded away from 0 with probability 1. In contrast,  the small scale fading component $\Bz_j$ in this work can follow any distribution, as long as the $4 + \varepsilon$ moment exists. Similarly, in \cite{JWsilverstein1985wishartsmall}, the authors showed that when $\BOmega_j = \BI_p$, $\forall j \in [n]$ and $\BA = \Bzero_{p \times n}$, the minimal eigenvalue of $\BS$ is bounded for the Gaussian ensemble. This result was extended to the non-Gaussian matrices in \cite{Bai1993smalleig}. It is important to note that the proofs in \cite{AblaNoEigen} and \cite{JWsilverstein1985wishartsmall} rely on the Gaussianity. It is very challenging to generalize the boundedness of the smallest eigenvalue to non-Gaussian matrices.
\end{remark}

\section{Simulations}
\label{Sec_Simulations}
In this section, we validate the accuracy of the theoretical results by simulations.
\par
\begin{figure}[t]
    \centering
    \includegraphics[width=4.1in]{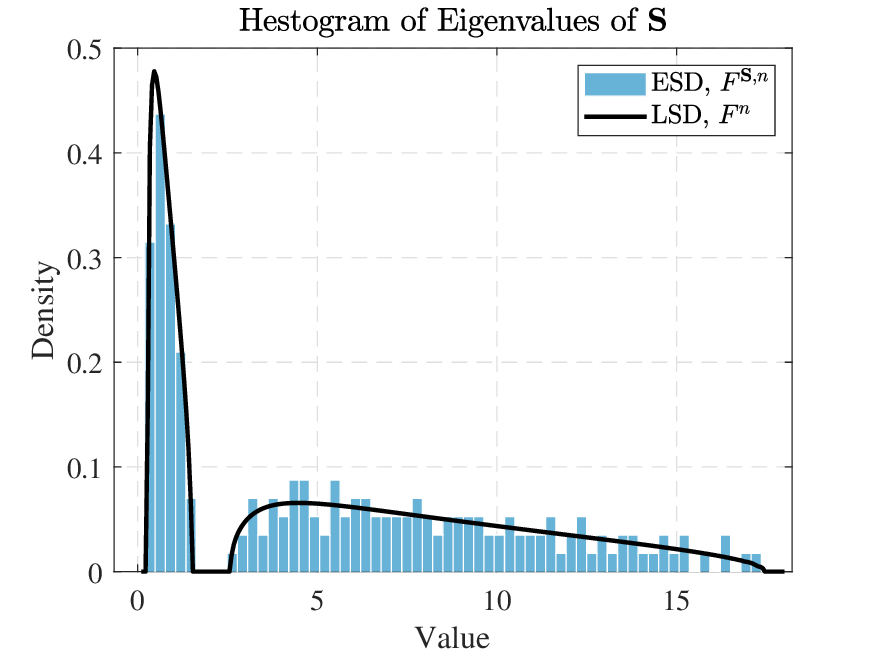}
    \caption{ESD vs LSD.}
    \label{Fig_ESD_LSD}
\end{figure}
\textbf{Convergence of the ESD:}
 Fig. \ref{Fig_ESD_LSD} compares the ESD $F^{\BS, n}$ with the LSD $F^n$. Here, the dimensions are set as $p = 200$ and $n = 400$. The correlation matrices $\{ \BOmega_j \}$ are set as $\BOmega_j = \diag\big\{\BI_{100}, 8\BI_{100} + (n + j) / (2n) \BK_j \big\}$, where $\BK_j = \diag(z_{j1}^2, \ldots, z_{j100}^2)$ with $z_{ji}$s being i.i.d. standard real Gaussian random variables. The mean of $\BSigma$ is given by $\BA = 2 \BE_{11}  -2 \exp(-0.6\pi \jmath) \BE_{22}$, where $\BE_{ij}$ denotes the all-zero matrix with only the $(i, j)$-th element equal to 1. The elements of the double array $\{X_{ij}\}$ are i.i.d. and follow the standard complex Gaussian distributions.
 From Fig. \ref{Fig_ESD_LSD}, it can be observed that the ESD closely aligns with the LSD, which validates Theorem \ref{Thm_conver_ESD}.
 \par
\begin{figure}[t]
    \centering
    \includegraphics[width=4.1in]{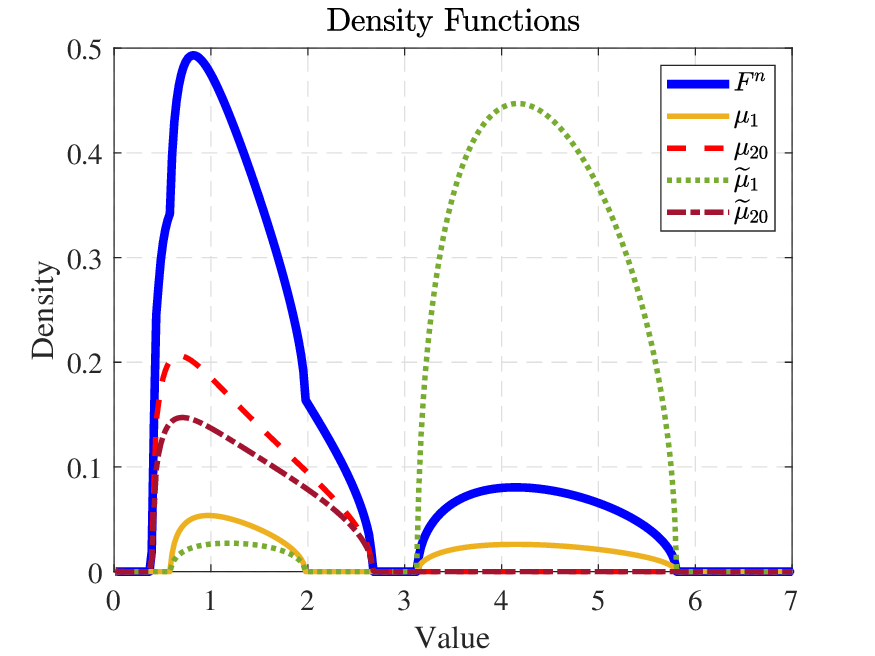}
    \caption{The Density Functions of Different Measures.}
    \label{Fig_Density_Merasures}
\end{figure}
\textbf{The Support of the Measures:} Fig. \ref{Fig_Density_Merasures} depicts the measures $\{ \mu_j \}$, $\{ \widetilde{\mu}_j\}$, and $F^n$. The dimensions are set as $p = 6$ and $n = 20$. The correlation matrices are generated as follows. First, let $\BOmega_j = \diag(1 + (i+j) / (p + n); 1 \leq i \leq p)$. Next, for $j_1 = 1, 2, \ldots 5$, we set $\BOmega_{j_1} \leftarrow \BOmega_{j_1} \BD_1$ and  for $j_2 = 16, 17, \ldots 20$, we set $\BOmega_{j_2} \leftarrow \BOmega_{j_2} \BD_2$, where $\BD_1 = \diag(\BI_3, \Bzero_3)$ and $\BD_2 = \diag(\Bzero_3, \BI_3)$. The mean of $\BSigma$ is given by $\BA = 2 \BE_{11} + \exp(0.4 \pi \jmath) \BE_{22}$. For clarity, only five densities $\mu_1$, $\mu_{20}$,  $\widetilde{\mu}_{1}$, $\widetilde{\mu}_{20}$, and $F^n$ are plotted. It can be observed from Fig. \ref{Fig_Density_Merasures} that the non-zero regions of $\{ \mu_j \}$ and $\{ \widetilde{\mu}_j\}$ are subsets of  $F^n$, which validates Proposition \ref{Prop_same_support}.
\par
\begin{figure}[t!]
\centering
 \subfloat{
\includegraphics[width=4.2in]{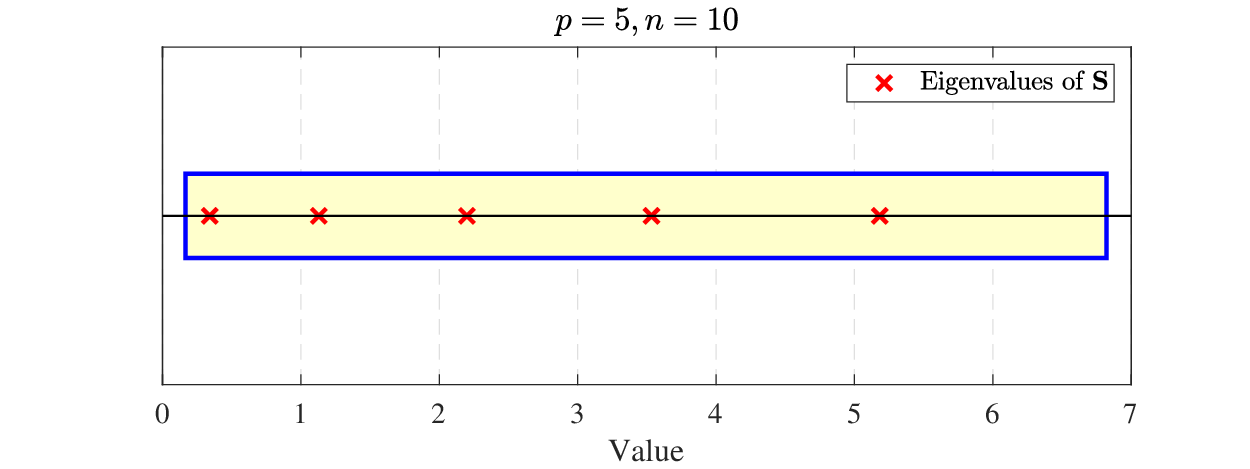}
}
\hfill
 \subfloat{
\includegraphics[width=4.2in]{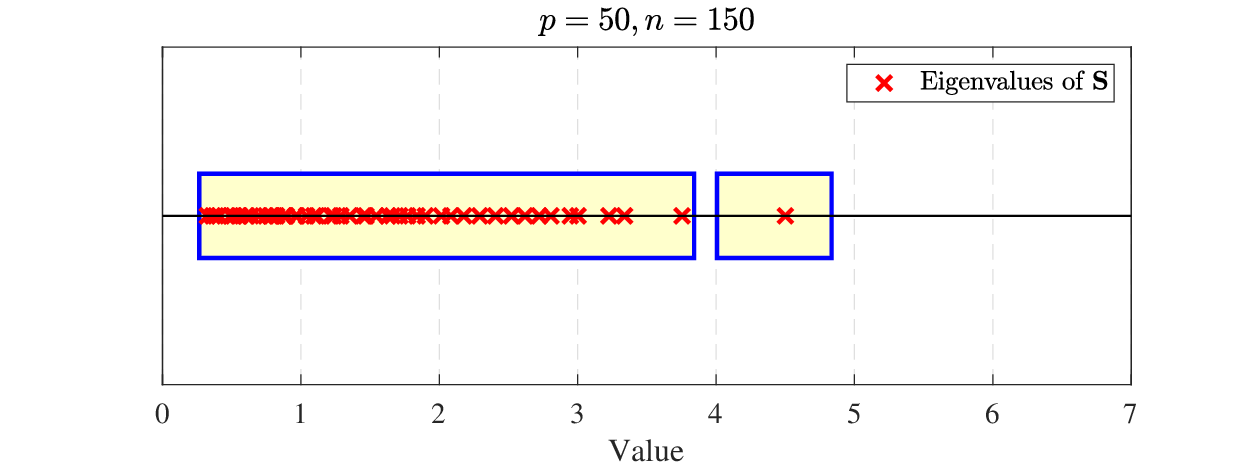}
}
\caption{Locations of the Empirical Eigenvalues. \label{Fig_No_eigenvalues}} 
\end{figure}

\textbf{Locations of the Empirical Eigenvalues:} Fig. \ref{Fig_No_eigenvalues} plots the theoretical support and the positions of the empirical eigenvalues of $\BS$. Here, the deterministic matrices are set as $\BOmega_j = \diag(1 + (i + j)/(p + n); 1 \leq i \leq p)$, and $\BA = \BE_{11} + 1.5 \exp(-0.2 \pi \jmath)\BE_{22}$. The entries of $\{X_{ij}\}$ are i.i.d. and follow the uniform distribution in $[-\sqrt{3}, \sqrt{3}]$. The blue box represents the theoretical support, while the red crosses indicate the eigenvalues of $\BS$. From Fig. \ref{Fig_No_eigenvalues}, it can be observed that there are no eigenvalues outside the support, which validates Theorem \ref{Thm_No_Eigenvalues}. Moreover, the smallest eigenvalue lies within the limiting support and is bounded away from zero, which verifies Corollary \ref{Coro_min_eig_zero_forcing}.

\begin{figure}[t]
    \centering
    \includegraphics[width=4.1in]{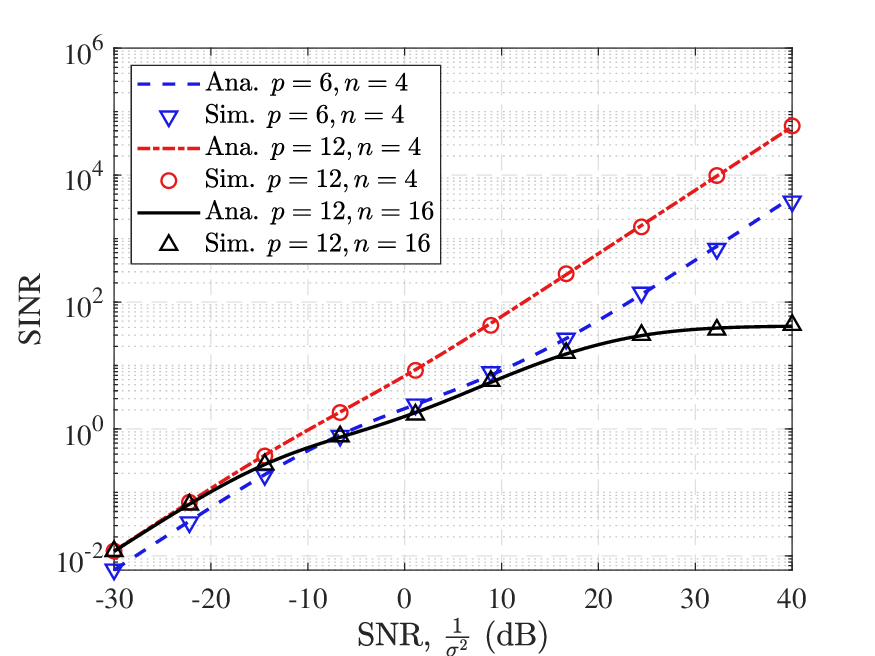}
    \caption{SINR versus SNR.}
    \label{Fig_SINR}
\end{figure}
\textbf{Analysis of the SINR with LMMSE Receivers:} Fig. \ref{Fig_SINR} shows the SINR of uplink MIMO systems with LMMSE receivers where the SNR is defined as $1 / \sigma^2$. The spatial correlation matrices $\{ \BC_j \}$ in \eqref{Eq_Rician_model} are modeled as $[\BC_j]_{k, l} = q(j)^{|k-l|}$, and the LoS components $\{\overline{\Bz}_j\}$ is generated according to the model $ \overline{\Bz}_j = [1, \ldots, \exp\{\jmath (p-1) \pi \sin\theta(j)\}]^T$. Here, the parameters $\{q(j)\}_{1 \leq j \leq n}$ and $\{\theta(j)\}_{1\leq j \leq n}$ are set as $q(j) = 0.7 + 0.2 \cdot \frac{j}{n}$ and $ \theta(j) = \frac{\pi}{2} + \frac{j \pi}{n}$, and the Rician factor is set as $\tau = 1$.  In Fig. \ref{Fig_SINR}, the markers represent the results of 1000 Monte Carlo simulations for \eqref{Eq_SINR_LMMSE}, and the lines denote the analysis results in Corollary \ref{Coro_SINR}.  It can be observed that the analysis and simulation results overlap very well, which validates the accuracy of Corollary \ref{Coro_SINR}. By comparing the red and blue curves, it can be seen that for the same SNR, when the number of users is fixed, the SINR increases as the number of BS antennas increases. Additionally, by considering the red and black curves, we can observe that for a fixed number of antennas, the SINR decreases as the number of users increases, and this reduction is not significant at low SNR levels.

\section{Conclusion}
\label{Sec_Conclusion}
In this paper, we investigated the spectrum of generally correlated and noncentral sample covariance matrices. Specifically, we demonstrated the almost sure convergence of the trace and bilinear forms of the resolvent, which indicates that the ESD converges weakly to the LSD almost surely. Furthermore, we proved that outside the support of the LSD, no eigenvalues will appear with probability one, thereby strengthening the notion of weak convergence. From a practical perspective, we applied the theoretical results to two communication systems and obtained new insights. In particular, based on the convergence of the resolvent, we derived the asymptotic SINR for multi-user uplink MIMO systems with LMMSE receivers over generally correlated Rician channels. We also showed that the minimum eigenvalue of the Gram matrix for the generally correlated Rayleigh MIMO channel is almost surely bounded away from zero, indicating that the ZF precoding matrix is valid in multi-user MIMO downlink systems. Future research directions include exploring the exact separation of eigenvalues and the fluctuations of extreme eigenvalues for the generally correlated and noncentral random matrices.
%%%%%%%%%%%%%%%%%%%%%%%%%%%%%%%%%%%%%%%%%%%%%%%%%%%%%%%%%%%%
\appendices
\section{Mathematical Tools}
\label{App_mathmatical_tools}
In this section, we present the mathematical tools that will be frequently used in the proofs throughout this paper.
\begin{lemma}(\cite[Lemma 2.8]{bai1998no})\label{Lemm_norm_mat} For $p$-by-$n$ matrices $\BA$ and $\BB$ with respective singular values $\sigma_1 \geq \sigma_2 \geq \ldots \geq \sigma_q$ and $\tau_1 \geq \tau_2 \geq \ldots \geq \tau_q$, where $q = \min(p, n)$, we have
\begin{equation}
  \max_{1 \leq k \leq q}|\sigma_k - \tau_k| \leq \norm{\BA - \BB}.
\end{equation}  
\end{lemma}
\begin{lemma}(\cite[Lemma 2.7]{bai1998no})\label{Lemm_trace_pre}
    Let $\Bx = [X_1, X_2, \ldots, X_p]^T$ be a vector whose entries $X_i$s are i.i.d. (complex) random variables with $\E X_1 = 0$ and $\E|X_1|^2 = 1$. Let $\BM$ be a $p \times p$ deterministic (complex) matrix. Then, for any $t \geq 2$, we have
    \begin{equation}
        \E \abs{\Bx^H \BM \Bx - \Tr \BM}^t \leq K_t \left\{ \left[\E|X_1|^4 \Tr (\BM\BM^H)\right]^{\frac{t}{2}} + \E |X_1|^{2t} \Tr (\BM\BM^H)^{\frac{t}{2}} \right\},
    \end{equation}
    where $K_t> 0$ is a constant only dependent on $t$. 
\end{lemma}
\begin{lemma}
  \label{Lemm_trace}
  Let $\Bx = [X_1, X_2, \ldots, X_p]^T$ be a vector with i.i.d. (complex) entries where $\E X_1 = 0$, 
  $\E |X_1|^2 = 1$, $\E |X_1|^{4 + \varepsilon} < \infty$, and $|X_1| \leq n^{\frac{1}{2} - \alpha}$, with 
  $\varepsilon > 0$ and $\alpha = \frac{\varepsilon}{8 + 2 \varepsilon}$. Let $\BM$ be a $p$-by-$p$ deterministic (complex) matrix. Then, we have, for any $t \geq 2$,
  \begin{equation}
    \E \abs{\Bx^H \BM \Bx - \Tr \BM}^t \leq  K_t \norm{\BM}^t R_t(p, n),
  \end{equation}
where $K_t > 0$ is a constant only dependent on $t$. Here, we have $R_{t}(p, n) = p^{\frac{t}{2}} + pn^{t - 2\alpha t - 2}$ when $t \geq 2 + \frac{\varepsilon}{2}$, and $R_t(p, n) = p^{\frac{t}{2}}$, otherwise. Moreover, for a deterministic vector $\Bv$, we have
\begin{align}
  &\E \abs{\Bx^H\Bv}^{2t} \leq K_t \norm{\Bv}^{2t} T_t(p, n), \\
  & \E \abs{\Bx^H\Bv\Bv^H \Bx -\Bv^H\Bv}^{t} \leq K_t \norm{\Bv}^{2t} T_t(p, n),
\end{align}
where $T_t(p, n) = 1$ for $1 \leq t \leq 2 + \frac{\varepsilon}{2}$ and $T_t(p, n) = n^{t - 2\alpha t - 2}$ when $t \geq 2 + \frac{\varepsilon}{2}$.
\end{lemma}
\textit{Proof:} The proof of Lemma \ref{Lemm_trace} is a direct application of Lemma \ref{Lemm_trace_pre} and omitted here. \qed
% \begin{lemma}
%   \label{Lemm_trace}
%   Let $\Bx = [x_1, x_2, \ldots, x_p]^T$ has i.i.d. (complex) entries with $\E x_1 = 0$, 
%   $\E |x_1|^2 = 1$, $\E |x_1|^{4 + \varepsilon} < \infty$, and $|x_1| < n^{\frac{1}{2} - \alpha}$ where 
%   $\varepsilon > 0$ and $\alpha = \frac{\varepsilon}{8 + 2 \varepsilon}$. We have, for any $k \geq 2$,
%   \begin{equation}
%     \E \abs{\Bx^H \BM \Bx - \Tr \BM}^k \leq  C_k \norm{\BM}^k R_k(p, n),
%   \end{equation}
% where $C_k > 0$ is a constant only depended on $k$ and $R_{k}(p, n) = p^{\frac{k}{2}} + pn^{k - 2\alpha k - 2}$ when $k \geq 2 + \frac{\varepsilon}{2}$ and $R_k(p, n) = p^{\frac{k}{2}}$ otherwise. Moreover, for deterministic vector $\Bv$,  
% \begin{equation}
%   \E \abs{\Bx^H\Bv}^{2k} \leq C_k \norm{\Bv}^{2k} T_k(p, n),
% \end{equation}
% where $T_k(p, n) = 1$ for $1 \leq k \leq 2 + \frac{\varepsilon}{2}$ and $T_k(p, n) = n^{k - 2\alpha k - 2}$ when $k \geq 2 + \frac{\varepsilon}{2}$.
% \end{lemma}
\begin{lemma} \label{Lemm_Mat_Ber}
(Matrix Bernstein, \cite[Theorem 1.4]{tropp2012user}, \cite[Theorem 5.29]{vershynin2010introduction}) Consider a finite sequence $\{ \BX_i\}_{1 \leq i \leq n}$ of independent Hermitian random $p$-by-$p$ matrices. Assume that $\E \BX_i = 0$, $\norm{\BX_i} \leq G$ almost surely, and $\norm{\sum_{i=1}^n \E \BX_i^2} \leq \sigma^2$. Then, for all $t > 0$, we have 
\begin{equation}
\PP \left\{ \norm{\sum_{i=1}^n\BX_i} \geq t  \right\} \leq  2p \exp \left[-c \min\left( \frac{t^2}{\sigma^2}, \frac{t}{G}\right) \right],
\end{equation}
where $c > 0$ is a absolute constant that does not depend on $p, \sigma$, and $G$.
\end{lemma}
\begin{lemma}
  \label{Lemm_Mat_norm_Pbound}
  Let $\BX = [\Bx_1, \ldots, \Bx_n] \in \mathbb{C}^{p \times n}$ be a random matrix with independent columns. Assume that $\E [\Bx_j] = \mathbf{0}_p$, $\PP(\norm{\Bx_j} \leq \sqrt{m}) = 1$, and $\E [\Bx_j\Bx_j^H] = \BOmega_j$  for $j \in [n]$. Define $u = \max_{1 \leq j \leq n} \norm{\BOmega_j}$. Then, for any $\delta \geq  0$ and deterministic matrix $\BT$, we have
  \begin{equation}
    \PP \left( \norm{n^{-\frac{1}{2}}\BX + \BT} \geq \norm{\BT} + \sqrt{u} + \delta \right) \leq 2p \exp \left(-\frac{cn \delta^2}{2m} \right),
  \end{equation}
  where $c > 0$ is a absolute constant.  
\end{lemma}
\textit{Proof:} The proof is an application of Lemma \ref{Lemm_Mat_Ber} and given in Appendix \ref{App_Prof_Mat_norm_Pbound}. \qed
\begin{lemma} 
\label{Lemm_Burk_Ineq}
(Burkholder Inequality, \cite[Lemma 2.1]{bai1998no}) Let $\{ X_k \}$ be a complex martingale difference sequence with respect to the filtration $\{\mathscr{F}_k\}$. Then, for $t \geq 1$, there exists constant $K_t > 0$ such that 
\begin{equation}
    \E \abs{\sum_{k=1}^n X_k}^{2t} \leq K_t \left[ \E\left[\sum_{k=1}^n \E(|X_k|^2| \mathscr{F}_{k-1}) \right]^t + \sum_{k=1}^n\E|X_k|^{2t} \right].
\end{equation} 
\end{lemma}
% \begin{lemma} Let $z \in c^H$
    
% \end{lemma}
\begin{lemma}(\cite[Lemma 5.2]{hachem2008varianceProfile})
  \label{Lemm_Spr_pMat}
  Let $\BA \in \mathbb{R}^{n \times n}$, $\Bu \in \mathbb{R}^n$, and $\Bv \in \mathbb{R}^{n}$. Assume that $[\BA]_{i, j} \geq 0$, $[\Bu]_j > 0$, and $[\Bv]_j > 0$, for any $i, j \in [n]$, and the following holds
  \begin{equation}
    \Bu = \BA \Bu + \Bv.
  \end{equation} 
  Then, the spectral radius of $\BA$ can be given by $\rho(\BA) \leq 1 - \frac{\min_{1 \leq i \leq n}([\Bv]_i)}{\max_{1 \leq i \leq n}([\Bu]_i)}$. Moreover, the matrix $\BI_n - \BA$ is invertible with $\min_{i, j} [(\BI_n - \BA)^{-1}]_{i, j} \geq 0$ and $\| (\BI_n - \BA)^{-1} \|_{\infty} \leq \frac{\max_{1 \leq i \leq n}([\Bu]_i)}{\min_{1 \leq i \leq n}([\Bv]_i)}$.
\end{lemma}
\begin{lemma} \label{Lemm_matrix_Spr_maxrow_Norm}
Let $\{\BA_l\}_{1 \leq l\leq N}$ be a collection of $n$-by-$n$ complex matrices, and $\{\BB_{l}\}_{1 \leq l\leq N}$ and $\{\BC_{l}\}_{1 \leq l\leq N}$ be collections of $n$-by-$n$ non-negative matrices. Assume that for each $l$, $
      \abs{[\BA_l]_{i,j}} \leq \sqrt{[\BB_l]_{i, j} [\BC_{l}]_{i, j}}
$ and denote $\BA = \prod_{l=1}^N \BA_{l}$, $\BB = \prod_{l=1}^N \BB_{l}$, and $\BC = \prod_{l=1}^N \BC_{l}$. Then, we have
\begin{equation}
    \abs{[\BA]_{i, j}} \leq \sqrt{[\BB]_{i, j}[\BC]_{i, j}}, \text{ and } \rho\left( \BA \right) \leq \sqrt{\rho\left(\BB\right)\rho( \BC )}.
\end{equation}
Moreover, if $\max(\rho(\BB), \rho(\BC)) < 1$, we have
\begin{align}
&\norm{(\BI_n - \BA)^{-1}}_{\infty} \leq \sqrt{\norm{(\BI_n - \BB)^{-1}}_{\infty}\norm{(\BI_n - \BC)^{-1}}_{\infty} }, \\
&\norm{(\BI_n - \BA)^{-1}}_{1} \leq \sqrt{\norm{(\BI_n - \BB)^{-1}}_{1}\norm{(\BI_n - \BC)^{-1}}_{1} }.
\end{align}
\textit{Proof: } This lemma is a generalization of \cite[Lemma 13]{AblaNoEigen} whose proof follows a similar logic and thus is omitted. \qed
\end{lemma}
\begin{lemma}
\label{Lemm_Spnorm_Bound_of_sum_semidef_matrices}
    Let  $\{ \BA_l\}_{1 \leq l\leq N}$ be a collection of $n$-by-$n$ Hermitian non-negative matrices and $\{a_l\}_{1 \leq l\leq N}$ be a sequence of complex numbers with $|a_l| \leq b_l$, for any $ l \in [N]$. Then, we have
    \begin{equation}
        \norm{\sum_{l=1}^N a_l\BA_l} \leq \norm{\sum_{l=1}^N b_l\BA_l}.
    \end{equation}
\end{lemma}
\textit{Proof:} By definition of the operator norm, we have
\begin{align*}
    &\norm{\sum_{l=1}^N a_l\BA_l} = \sup_{\norm{\Bx} = 1, \norm{\By} = 1} \abs{ \sum_{l=1}^N a_l \Bx^H \BA_l \By} \leq \sup_{\norm{\Bx} = 1, \norm{\By} = 1}\sum_{l=1}^N b_l \abs{\Bx^H \BA_l^{1/2} \BA_l^{1/2}  \By}  \\
    & \leq \sup_{\norm{\Bx} = 1, \norm{\By} = 1}  \sum_{l=1}^N b_l \sqrt{\Bx^H \BA_l \Bx \By^H \BA_l \By} \leq \sup_{\norm{\Bx} = 1, \norm{\By} = 1} \left(\sum_{l=1}^N b_l \Bx^H \BA_l \Bx \right)^{1/2} \left(\sum_{l=1}^N b_l\By^H \BA_l \By\right)^{1/2} \leq \norm{\sum_{l=1}^N b_l \BA_l}, \numberthis
\end{align*}
which completes the proof. \qed
\section{Proof of Lemma \ref{Lemm_Mat_norm_Pbound}}
\label{App_Prof_Mat_norm_Pbound}
Denote $\BX_j = \frac{1}{n} (\Bx_j\Bx_j^H - \BOmega_j)$ for $1 \leq j \leq n$. Next, we verify the conditions of Lemma \ref{Lemm_Mat_Ber}. The norm of ${\BX_j}$ can be bounded as
\begin{equation}
  \norm {\BX_j} = \frac{\norm{\Bx_j\Bx_j^H - \BOmega_j}}{n}  \leq \frac{\norm{\Bx_j}^2 + \norm{\BOmega_j}}{n}  \leq \frac{m + u}{n} \leq \frac{2m}{n} := G.
\end{equation}
Furthermore, the following holds
\begin{equation}
  \norm{\sum_{j=1}^n \E \BX_j^2} \leq \sum_{j=1}^n\norm{ \E \BX_j^2} \leq \sum_{j=1}^n\frac{\norm{ \E\norm{\Bx_j}^2 \Bx_j\Bx_j^H} + \norm{\BOmega_j}^2}{n^2} \leq \frac{mu + u^2}{n} \leq \frac{2mu}{n} := \sigma^2.
\end{equation}
By applying Lemma \ref{Lemm_Mat_Ber} and taking $t = \max{(\sqrt{u}\delta, \delta^2)}$ with $\delta \geq  0$, we have 
\begin{equation}
  \PP\left(n^{-1}\norm{\BX\BX^H - \overline{\BOmega}} \geq \max{(\sqrt{u}\delta, \delta^2)}\right) \overset{(a)}{\leq} 2p \exp \left(-\frac{cn \delta^2}{2m}  \right),
\end{equation}
with $c > 0$ being a absolute constant and $\overline{\BOmega} = \sum_{j=1}^n \BOmega_j$. Here step $(a)$ is due to the following identity 
\begin{equation}
  \min \left(u^{-1} \max(u\delta^2, \delta^4), \max(\sqrt{u}\delta, \delta^2) \right) = \delta^2, ~~ \delta \geq 0.
\end{equation}
For the event $\{ n^{-1}\norm{\BX\BX^H - \overline{\BOmega}} \leq \max{(\sqrt{u}\delta, \delta^2)}\}$, the following holds
\begin{equation}
  \norm{n^{-\frac{1}{2}}{\BX}}^2  = \norm{n^{-1}\BX\BX^H} \leq n^{-1}\norm{{\BX\BX^H - \overline{\BOmega}}} + n^{-1}\norm{\overline{\BOmega}} \leq \max{(\sqrt{u}\delta, \delta^2)} + u \leq (\sqrt{u} + \delta)^2,
\end{equation}
which implies $\Vert n^{-\frac{1}{2}}\BX + \BT \Vert \leq \sqrt{u} + \delta + \norm{\BT}$. Therefore, we complete the proof. \qed
\section{Truncation, Centralization, and Normalization}
\label{App_Trunc}
To study the asymptotic behavior for the eigenvalues of $\BS$, we first simplify Assumption \ref{Assumpt_2}. 
In particular, in the truncation step, we restrict the range of the random variables $X_{ij}$ such that $|X_{ij}|$ is bounded by a suitable polynomial in $n$. In the centralization and normalization step, we subtract the mean from the truncated random variables and divide the results by the standard deviation to ensure that the transformed random variables are standardized.
% In particular, we will truncate the double array of random variables $X_{ij}$ to make them bounded by suitable polynomials in $n$, such that the moments are also polynomially bounded. 
This process will not affect the asymptotic locations of the random eigenvalues of $\BS$.
\par
Define the parameter $\alpha = \frac{\varepsilon}{8 + 2\varepsilon}$. Since $\E |X_{11}|^{4 + \varepsilon} < \infty$, we have $(4 + \varepsilon)\int_{0}^{\infty} t^{3 + \varepsilon} \PP(|X_{11}| \geq t) \dd t = \E|X_{11}|^{4 + \varepsilon} < \infty$, which implies
\begin{equation}
    \sum_{m=1}^{\infty} 2^{2m} \PP\left( |X_{11}| \geq 2^{(\frac{1}{2}-\alpha)m} \right) < \infty.
\end{equation}
Define the truncated random variables $X_{ij, n} = X_{ij} \Ind\{|X_{ij}| \leq n^{\frac{1}{2} - \alpha}\}$, and construct the random matrix $\widehat{\BSigma} = \BA + \widehat{\BY}$ such that $\widehat{\BY}$ has the same structure as $\BY$ by replacing $X_{ij}$ with $X_{ij, n}$ (note that we omit the index $n$ in $\BSigma_n$ for simplicity). Let $\ell^+ = \max(\limsup_n \max_j {d_j} / n, \limsup_n p / n)$.  Then, as $n \to \infty$, the probability of the event $\{\BSigma \neq \widehat{\BSigma} ~~ \mathrm{i.o.}\}$ can be bounded by
\begin{align*}
    &\PP\left(\BSigma \neq \widehat{\BSigma} ~~ \mathrm{i.o.} \right) =  \PP\left( \bigcap_{k\geq 1} \bigcup_{n \geq 2^k} \left\{\BSigma \neq \widehat{\BSigma}\right\} \right) = \lim_{k \to \infty} \PP\left( \bigcup_{n \geq 2^k} \left\{\BSigma \neq \widehat{\BSigma}\right\} \right) \\
    & \leq \lim_{k \to \infty} \sum_{m=k}^{\infty}\PP\left( \bigcup_{2^m \leq n \leq 2^{m+1}}\bigcup_{j \leq n} \bigcup_{i \leq d_j} \left\{ |X_{ij}|\geq n^{0.5 - \alpha}  \right\} \right) \\
    &\overset{(a)}{\leq} \lim_{k \to \infty} \sum_{m=k}^{\infty}\PP\left( \bigcup_{2^m \leq n \leq 2^{m+1}}\bigcup_{j \leq n, i\leq 2\ell^+n}  \left\{ |X_{ij}| \geq 2^{(0.5 - \alpha)m} \right\} \right) \\
    & = \lim_{k \to \infty} \sum_{m=k}^{\infty}\PP\left( \bigcup_{j \leq 2^{m+1}, i \leq \ell^+2^{m+2}} \left\{ |X_{ij}| \geq 2^{(0.5 - \alpha)m} \right\} \right) \\
    & \leq \lim_{k \to \infty} \sum_{m=k}^{\infty} \ell^+ 2^{2m+3}\PP\left(  |X_{11}| \geq 2^{(0.5 - \alpha)m}  \right) = 0, \numberthis
\end{align*}
where step $(a)$ follows from $d_j \leq 2 \ell^+ n$ for sufficiently large $n$. The above derivation shows that truncation does not affect the limits for the extreme eigenvalues of $\BS = \BSigma\BSigma^H$. 
Next, we perform the centralization. Denote the centralized random variables $\widetilde{X}_{ij, n} = {X}_{ij, n} - \E {X}_{ij, n}$ and construct the matrix $\widetilde{\BSigma} = \BA + \widetilde{\BY}$ with the same structure as $\BSigma$ by replacing $X_{ij}$ with $\widetilde{X}_{ij, n}$. By Lemma \ref{Lemm_norm_mat}, we have 
\begin{align*}
    &\max_{1 \leq k \leq p} \abs{\lambda_k^{\frac{1}{2}}(\widehat{\BS}) - \lambda_k^{\frac{1}{2}}(\widetilde{\BS}) } \leq \norm{\widehat{\BSigma} - \widetilde{\BSigma}} = \frac{1}{\sqrt{n}} \sup_{\norm{\Ba} = 1} \norm{\sum_{j=1}^n a_j \BB_j(\widehat{\Bx}_j - \widetilde{\Bx}_j)} \\
    &\leq \frac{1}{\sqrt{n}} \sup_{\norm{\Ba} = 1} \sum_{j=1}^n |a_j| \norm{\BB_j}d_j^{\frac{1}{2}}\sup_{i \leq d_j} \abs{ {X}_{ij, n} - \widetilde{X}_{ij, n}} \leq K n^{\frac{1}{2}} \sup_{i, j \geq 1} \E \abs{ {X}_{ij} - {X}_{ij, n}} \overset{(a)}{\leq}  \frac{K_1}{n^{1 + \alpha}}, \numberthis
\end{align*}
where $K$ and $K_1$ are constants independent of $n$ and step $(a)$ is obtained by setting $\beta = 1$ in the following inequality
\begin{align*}
\label{Eq_central_moment_X}
    &\E |X_{ij} - X_{ij, n}|^{\beta} \leq \beta \int_{0}^{\infty} \PP\left(|X_{ij}| \geq \max(n^{(\frac{1}{2} - \alpha)}, t) \right) t^{\beta - 1} \dd t  \\
    &\leq n^{(\frac{1}{2} - \alpha)\beta}   \PP(|X_{ij}| \geq n^{\frac{1}{2} - \alpha}) + \beta \int_{n^{\frac{1}{2} - \alpha}}^{\infty} t^{\beta - 1}\PP(|X_{ij}| \geq t) \dd t \\
    &\leq  \frac{\E|X_{ij}|^{4 + \varepsilon} + \beta\int_{n^{\frac{1}{2} - \alpha}}^{\infty} t^{3 + \varepsilon}\PP(|X_{ij}| \geq t) \dd t }{n^{(\frac{1}{2} - \alpha)(4 + \varepsilon - \beta)}} \leq \frac{(4 + \varepsilon + \beta)\E |X_{ij}|^{4 + \varepsilon}}{(4 + \varepsilon)n^{2 - (\frac{1}{2} - \alpha)\beta}}, ~~ \beta \in [1, 4 + \varepsilon]. \numberthis
\end{align*}
The above analysis shows that $\limsup_{n \to \infty}\max_{1 \leq k \leq p} |\lambda_k(\widehat{\BS}) - \lambda_k(\widetilde{\BS}) | = 0$.
\par
In the following, we normalize the random variables. Denote $\sigma_n^2 = \E |\widetilde{X}_{11, n}|^2$ and $\overline{X}_{ij, n} = \widetilde{X}_{ij, n} / \sigma_n$. Similar to the above derivation, we define $\overline{\BSigma} = \BA + \overline{\BY}$. Then, we can obtain
\begin{equation}
    1 - \sigma_n^2 = 1 - \E \abs{{X}_{11, n}}^2 + \abs{\E \widehat{X}_{11, n}}^2 = \E|X_{11}|^2 \Ind{\{|X_{11}| \geq n^{\frac{1}{2} - \alpha} \}} + \abs{\E|X_{11}| \Ind{\{|X_{11}| \geq n^{\frac{1}{2} - \alpha} \}} }^2 \overset{(a)}{ \leq } \frac{K}{n^{1 + 2\alpha}},
\end{equation}
where step $(a)$ follows by setting $\beta = 2$ in \eqref{Eq_central_moment_X}. Hence, we have 
\begin{align*}
 \label{Eq_Max_lambda_k_SS}
    \max_{1 \leq k \leq p} \abs{\lambda_k^{\frac{1}{2}}(\overline{\BS}) - \lambda_k^{\frac{1}{2}}(\widetilde{\BS}) } \leq \norm{\overline{\BSigma} - \widehat{\BSigma}} \leq \frac{|1 - \sigma_n|}{\sqrt{n}|\sigma_n|} \sup_{\norm{\Ba} = 1}  \sum_{j=1}^n |a_j| \norm{\BB_j} \norm{\widetilde{\Bx}_j} \leq \frac{\sqrt{n}|1 - \sigma_n^2| K}{|\sigma_n(1 + \sigma_n)|} \sup_{i, j \geq 1} |\widetilde{X}_{ij, n}| \overset{(a)}{\leq} \frac{K_1}{n^{3 \alpha}}, \numberthis
\end{align*}
where $K$ and $K_1$ are non-negative constants and step $(a)$ in \eqref{Eq_Max_lambda_k_SS} follows from $|\widetilde{X}_{ij, n}| \leq |X_{ij, n}| + |\E (X_{ij} - X_{ij, n})| \leq 2n^{1-\alpha}$ for large $n$. By \eqref{Eq_central_moment_X}, it is straightforward to prove $\E |\overline{X}_{ij, n}|^{4 + \varepsilon} < \infty$.
\par
Based on the above discussion, we know that the asymptotic locations of the eigenvalues for $\BS$ and $\overline{\BS}$ are identical. Therefore, we can make the following simplified and stronger assumptions on the distribution of $X_{ij}$.
\begin{itemize}
    \item[1)] $\{ X_{ij} \}_{i, j \geq 1}$ is a double array whose entries are i.i.d. complex random variables.
    \item [2)] $\E X_{ij} = 0$, $\E|X_{ij}|^2 = 1$, and $\E |X_{ij}|^{4 + \varepsilon} < \infty$. 
    \item [3)] $|X_{ij}| \leq n^{\frac{1}{2} - \alpha}$.
\end{itemize}
\section{A Bound on the Spectral Norm}
In this section, we provide a loose bound on the largest eigenvalue of S. This bound will be necessary for investigating tighter results. Next, we first prove the following lemma.
\begin{lemma} 
  \label{Lemm_P_bound_Sn}
Assume Assumptions \ref{Assumpt_1}-\ref{Assumpt_4} hold and $X_{ij}$ follows the distribution as described in Appendix \ref{App_Trunc}. 
For sufficiently large $n$, the following inequality holds
\begin{equation}
  \label{Eq_P_norm_Sn}
  \PP\left[ \norm{\BSigma} \geq \log(n)\right] \leq  K_l n^{-l}, ~~ \forall l > 0,
\end{equation}
where $K_l$ is a constant that only depends on $l$. Furthermore, given $r \geq 1$, there exists a positive integer $N = N(r)$ such that the following holds
\begin{equation}
  \E \norm{\BSigma}^r \leq  2\log^r(n)
\end{equation}
for all $n \geq N(r)$.
\end{lemma}
\textit{Proof:} We first give a bound on $\PP(\norm{\BB_j\Bx_j} \geq \sqrt{Kn})$, where $K$ is a sufficiently large constant. According to Lemma \ref{Lemm_trace}, we have 
\begin{align*}
    \PP \left(\norm{\BB_j\Bx_j} \geq \sqrt{Kn} \right) \leq \PP \left(\abs{\Bx_j^H\BB_j^H\BB_j\Bx_j - \Tr \BOmega_j} \geq Kn/2 \right) \leq \frac{2^t\E\abs{\Bx_j^H\BB_j^H\BB_j\Bx_j - \Tr \BOmega_j}^t}{(Kn)^t} \leq \frac{K'_t}{n^{2\alpha t + 1}}, \numberthis
\end{align*}
where $t \geq 4$ is arbitrarily large and $K'_t$ is a constant dependent on $t$. Denote the set $A_n = \cap_{j \leq n} \left\{ \norm{\BB_j\Bx_j} \leq (Kn)^{1/2} \right\}$. Then, for any $a> 0$, we have 
\begin{equation}
\label{Eq_Norm_Bd_Sigma}
    \PP\left(\norm{\BSigma} \geq a \right) \leq \PP \left( \{\norm{\BSigma} \geq a\} \cap A_n \right) + \sum_{j=1}^n \PP\left(\norm{\BB_j\Bx_j} \geq \sqrt{Kn} \right) \leq \PP(\norm{\underline{\BSigma}} \geq a) + \frac{K_t'}{n^{2\alpha t}}.
\end{equation}
Define the truncated matrix $\underline{\BSigma} = \underline{\BY} + \BA$, where the $j$-th column of $\underline{\BY}$ is given by $n^{-\frac{1}{2}} \BB_j \Bx_j \Ind\left\{ \norm{\BB_j\Bx_j} \leq (Kn)^{1/2} \right\}$. According to Lemma \ref{Lemm_Mat_norm_Pbound},
there exist positive constants $K_1$ and $K_2$ such that for any $u \geq 0$, 
\begin{align*}
\label{Eq_Pbound_on_USigma}
    \PP\left( \norm{\underline{\BSigma}} \geq K_1 + u \right) \leq 2p \exp\left(-K_2 u^2\right). \numberthis
\end{align*}
Taking $u = \log(n) - K_1$ in \eqref{Eq_Pbound_on_USigma}, \eqref{Eq_Norm_Bd_Sigma} implies
\begin{equation}
    \PP( \norm{\BSigma} \geq \log(n)) \leq \frac{K'_t}{n^{2 \alpha t + 1}} + \frac{2p}{n^{K_2' \log(n)}},
\end{equation}
for sufficiently large $n$ and constant $K_2' > 0$. Since $t$ is arbitrary, \eqref{Eq_P_norm_Sn} is proved.
\par
Next, we prove the bound on the moment of $\norm{\BSigma}$. By $|X_{ij}| \leq n^{1/2 - \alpha}$, we can obtain
\begin{equation}
    \norm{\BSigma} \leq \norm{\BA} + \frac{1}{\sqrt{n}} \sup_{\norm{\Ba} = 1}\sum_{j=1}^n |a_j| \norm{\BB_j} \norm{\Bx_j} \leq \widetilde{K} n^{1 - \alpha}
\end{equation}
for some constant $\widetilde{K}$ and large $n$. As a result, we have 
\begin{align*}
    \E \norm{\BSigma}^r  &= \E \norm{\BSigma}^r \{ \Ind_{A_n} + \Ind_{A_n^c}\} \leq r \int_{0}^{\widetilde{K}n^{1-\alpha}} u^{r-1} \PP(\norm{\underline{\BSigma}} \geq u) \dd u +\E \norm{\BSigma}^r \Ind_{A_n^c} \\
    & \leq r \int_{0}^{\log(n)} u^{r-1} \dd u + 2pr \int_{\log(n)}^{\widetilde{K}n^{1-\alpha}} u^{r-1} e^{-K_2u^2/4} \dd u + \PP (A_n^c) (\widetilde{K}n^{1-\alpha})^r := W_1 + W_2 + W_3, \numberthis
\end{align*}
by applying \eqref{Eq_Pbound_on_USigma}. By straight calculation, we have $W_1 = \log^r(n)$. $W_2$ can be bounded by
\begin{equation}
    W_2 \leq 2p r (\widetilde{K}n^{1-\alpha})^{(r - 1)} \int_{\log(n)}^{\infty} e^{-K_2u^2 / 4} \dd u \leq \frac{K_r n^r}{n^{ \widetilde{K}_2 \log(n)}} \leq \frac{\log^r(n)}{2},
\end{equation}
for any $n$ larger than a positive $r$-related integer $N = N(r)$. Here, $\widetilde{K}_2$ is a constant and $K_r$ is a constant that depends on $r$. It follows from \eqref{Eq_Norm_Bd_Sigma} that $W_3 \leq \log^r(n) / 2$. Therefore, we complete the proof of Lemma \ref{Lemm_P_bound_Sn}. \qed
\section{Proof of Theorem \ref{Thm_conver_ESD}}
\label{App_Converge_resolvent}
As the proofs for the trace form and the bilinear form of the resolvent are similar, and the latter is more technically challenging, we will only provide the convergence proof for the bilinear form $\Bu^H (\BQ(z) - \BTheta(z)) \Bv$. Using the polarization identity
\begin{equation}
    \Bu^H \BT \Bv = \frac{1}{4} \left[ (\Bu + \Bv)^H\BT(\Bu + \Bv) - (\Bu - \Bv)^H \BT (\Bu - \Bv) +  \jmath(\Bu - \jmath\Bv)^H \BT (\Bu - \jmath \Bv) - \jmath (\Bu + \jmath \Bv)^H \BT (\Bu - \jmath \Bv) \right],
\end{equation}
it is sufficient to establish the convergence of the quadratic form $\Bu^H (\BQ(z) - \BTheta(z)) \Bu$ with $\norm{\Bu} = 1$. Here, we will prove a stronger result that not only gives the order of convergence but is also useful for demonstrating the no-eigenvalue property. To this end, we will use the simplified model presented in Appendix \ref{App_Trunc}. As discussed in Appendix \ref{App_Trunc}, the simplified process will not change the asymptotic positions of any eigenvalues of $\BS$.
\begin{lemma}
\label{Lemm_bilinear}
    Assume Assumptions \ref{Assumpt_1}-\ref{Assumpt_4} hold and $X_{ij}$ follows the distribution as described in Appendix \ref{App_Trunc}. Then, for $z = x + \jmath v \in \mathbb{C}^+$ where $v$ is less or equal to a small positive constant $K_{\mathrm{im}}$, the following holds for $t \geq 2 + \varepsilon / 2$
    \begin{equation}
        \sup_{\norm{\Bu} = 1} \E \abs{\Bu^H \left[\BQ(z) - \BTheta(z)\right] \Bu }^{2t} \lesssim_t \frac{[\mathcal{U}(z)]^{2t}}{n^{2\alpha t}}, \label{Eq_bilinear_U_t}
    \end{equation}
    where $\mathcal{U}(z) = (K_2 + x^2 + v^2)^{9} / v^{29}$ with $K_2$ being a positive constant independent of $n$.
\end{lemma}
\textit{Proof:} The main techniques follow from \cite{hachem2013bilinear}, which rely on the martingale difference argument and Burkholder inequality (Lemma \ref{Lemm_Burk_Ineq}). For simplicity of notation, we will omit the variable $z$ when there is no ambiguity. 
\par
Next, we introduce the related intermediate variables. Define $\widetilde{\BQ}(z) = (\BSigma^H\BSigma - z\BI_n)^{-1}$. Let $\BSigma_j$ be the $p$-by-$(n-1)$ matrix by deleting the $j$-th column of $\BSigma$ and denote $\BQ_j(z) = (\BSigma_j \BSigma_j^H - z \BI_p)^{-1}$. Write $\Bxi_j = \Ba_j + \By_j$, where vectors $\Bxi_j$ and $\Ba_j$ are the $j$-th column of $\BSigma$ and $\BA$, respectively, and $\By_j = n^{-1/2} \BB_j \Bx_j$. Further define
% \begin{align*}
% &\kappa_j = \frac{1}{n} \Tr \BOmega_j \E \BQ, ~~~ \eta_j =  \frac{1}{1 + \frac{1}{n}\Tr \BOmega_j\BQ_j + \Ba_j^H\BQ_j\Ba_j}, ~~~ \widetilde{q}_j = \frac{1}{-z(1 + \Bxi_j^H\BQ_j\Bxi_j)}, ~~ \widetilde{\kappa}_j = \E \widetilde{q}_j  \\
% &\alpha_{j, j} = \frac{1}{n} \Tr \BOmega_j \BQ_j, ~~ \Delta_j = \Bxi_j^H \BQ_j \Bxi_j - \frac{1}{n} \Tr \BOmega_j \BQ_j - \Ba_j^H \BQ_j\Ba_j, ~~ \BD_{\kappa} = \diag(\kappa_1, \ldots, \kappa_n)\\
% &\BR(z) = \left[-z\left(\BI_p + \sum_{j=1}^n\frac{\BOmega_j}{n} \widetilde{\kappa}_j(z)\right) + \BA \left(\BI_n + \BD_{\kappa}(z)\right)^{-1}\BA^H \right]^{-1}, \\
% &\widetilde{\BR}(z) = \left[ -z(\BI_n + \BD_{\kappa}(z)) + \BA^H \left( \BI_p +  \sum_{j=1}^n\frac{\BOmega_j}{n} \widetilde{\kappa}_j(z) \right) \BA \right]^{-1}.
% \end{align*}
\begin{align*}
&\kappa_j = \frac{1}{n} \Tr \BOmega_j \E \BQ, ~~~ \eta_j =  \frac{1}{1 + \frac{1}{n}\Tr \BOmega_j\BQ_j + \Ba_j^H\BQ_j\Ba_j}, ~~~ \widetilde{q}_j = \frac{1}{-z(1 + \Bxi_j^H\BQ_j\Bxi_j)},  \\
&\widetilde{\kappa}_j = \E \widetilde{q}_j, ~~ \alpha_{j, j} = \frac{1}{n} \Tr \BOmega_j \BQ_j, ~~ \Delta_j = \Bxi_j^H \BQ_j \Bxi_j - \frac{1}{n} \Tr \BOmega_j \BQ_j - \Ba_j^H \BQ_j\Ba_j, \\
&\BG = \left[-z\left( \BI_p + \sum_{j=1}^n \frac{{\BOmega_j \widetilde{\kappa}_j}}{n}\right) \right]^{-1}, ~~ \widetilde{\BG} = \diag\left(\frac{-1}{z(1 + \kappa_j)}; 1 \leq j \leq n \right), \\
& \BR = \left( \BG^{-1} -z\BA\widetilde{\BG}\BA^H \right)^{-1}, ~~ \widetilde{\BR} = \left( \widetilde{\BG}^{-1} -z\BA^H \BG \BA \right)^{-1}. \label{Label_Def_of_terms_bilinear}\numberthis
\end{align*}
Similar to \cite{hachem2013bilinear}, the random quantity $\Bu^H[\BQ - \BTheta]\Bu$ is divided into three parts as
\begin{equation}
    \Bu^H[\BQ- \BTheta]\Bu = \Bu^H[\BQ - \E \BQ]\Bu + \Bu^H[\E \BQ - \BR]\Bu + \Bu^H[ \BR -  \BTheta]\Bu,
\end{equation}
and we will handle them separately.  We use $\E_j[\cdot]$ to denote the conditional expectation $\E[\cdot | \mathscr{F}_j]$, where $\mathscr{F}_j$ is the $\sigma$-algebra generated by $\{ X_{ij'}, 1 \leq j' \leq j, 1 \leq i \leq m_{j'} \}$, and $\E_0 \equiv \E $. Some useful transforms and estimations are given in Appendix \ref{App_Useful_Estimations}. In the following, we will frequently use the asymptotic notation ``$\lesssim$'', and the asymptotic regime is stated in Assumption \ref{Assumpt_1} with $n \to \infty$.
\paragraph*{Part I} In this part, we will prove
\begin{equation}
\label{Eq_Part_I}
    \E \abs{\Bu^H[\BQ - \E \BQ]\Bu}^{2t}\lesssim_t \frac{\log^{6t}(n)}{v^{12t}n^{2 \alpha t + 0.5}},
\end{equation}
for $t \geq 2 + \varepsilon / 2$. First, write $\Bu^H (\BQ - \E \BQ) \Bu$ as the sum of martingale difference sequence, 
\begin{align*}
    &\Bu^H (\E_n \BQ - \E_0 \BQ) \Bu = \sum_{j=1}^n \Bu^H (\E_{j}\BQ - \E_{j-1} \BQ) \Bu \overset{(a)}{=}  \sum_{j=1}^n  [\E_{j}- \E_{j-1}]\Bu^H(\BQ - \BQ_j) \Bu \overset{(b)}{=}  \sum_{j=1}^n  [\E_{j}- \E_{j-1}]\frac{-\Bu^H\BQ_j\Bxi_j \Bxi_j^H\BQ_j \Bu}{1 + \Bxi_j^H \BQ_j \Bxi_j}  \\
    & \overset{(c)}{=} \sum_{j=1}^n  [\E_{j}- \E_{j-1}]\left( -\eta_j\Bu^H\BQ_j\Bxi_j \Bxi_j^H\BQ_j \Bu + \frac{\eta_j \Delta_j \Bu^H\BQ_j\Bxi_j \Bxi_j^H\BQ_j \Bu}{1 + \Bxi_j^H \BQ_j \Bxi_j} \right) = \sum_{j=1}^n \Gamma_{1j} + \Gamma_{2j}, \numberthis
\end{align*}
where step $(a)$ follows from $\E_{j}\BQ_j = \E_{j-1} \BQ_{j}$, step $(b)$ is due to the identity in \eqref{Eq_Wood_resolvent_Id}, and step $(c)$ follows by writing $z\widetilde{q}_j = -\eta_j - z\widetilde{q}_j \eta_j \Delta_j$.
One can check that both $\{\Gamma_{1j}\}_{j \leq n}$ and $\{ \Gamma_{2j} \}_{j \leq n}$ are martingale difference sequences with respective to the filtration $\{\mathscr{F}_j\}_{j \leq n}$. Next, we handle $\E |{\sum_j} \Gamma_{1j}|^{2t}$.  Express $\Gamma_{1j}$ as
\begin{align*}
    -\Gamma_{1j} &= [\E_{j} - \E_{j-1}]\eta_j\Bu^H\BQ_j\Ba_j \By_j^H\BQ_j \Bu + \eta_j\Bu^H\BQ_j\By_j \Ba_j^H\BQ_j \Bu + \eta_j\Bu^H\BQ_j\Ba_j \Ba_j^H\BQ_j \Bu + \eta_j\Bu^H\BQ_j\By_j \By_j^H\BQ_j \Bu \\
    & = Z_{1j} + Z_{2j} + Z_{3j} + Z_{4j}. \numberthis
\end{align*}
By  Lemma \ref{Lemm_Burk_Ineq}, we have
\begin{align*}
\label{Eq_Bound_Gamma_1}
    \E \abs{\sum_{j=1}^n \Gamma_{1j}}^{2t} \leq \sum_{i = 1}^4 K_t  \left[ \E \left(\sum_{j=1}^n\E_{j-1}|Z_{ij}|^2 \right)^t + \sum_{j=1}^n \E|Z_{ij}|^{2t} \right]. \numberthis 
\end{align*}
In the following, we will evaluate the right hand side (RHS) of \eqref{Eq_Bound_Gamma_1} for different $i \in \{1, 2, 3, 4\}$. Since $\By_j$ is independent of $\mathscr{F}_{j-1}$, we can obtain
\begin{equation}
    Z_{4j} = \E_j \eta_j\Bu^H\BQ_j\By_j \By_j^H\BQ_j \Bu  - \frac{1}{n}\eta_j\Bu^H\BQ_j\BOmega_j\BQ_j \Bu. 
\end{equation}
It follows from Lemma \ref{Lemm_Bound_on_qterms} and Lemma \ref{Lemm_trace} that
\begin{equation}
    \sum_{j=1}^n \E_{j-1} |Z_{4j}|^2  \lesssim v^{-2}\sum_{j=1}^n\E_{j-1} \abs{ \Bu^H\BQ_j\By_j \By_j^H\BQ_j \Bu  - \frac{1}{n}\Bu^H\BQ_j\BOmega_j\BQ_j \Bu }^2 \lesssim \frac{1}{nv^6}, \label{Eq_Burk_bound_Z4j}
\end{equation}
where we have used the fact that $|\eta_j| \leq K_0 / v$ since $v$ is small and the bound $\|\BQ_j \| \leq 1/v$. By Lemma \ref{Lemm_trace}, we have 
\begin{align*}
    \sum_{j=1}^n \E |Z_{4j}|^{2t} \lesssim_t v^{-2t} \sum_{j=1}^n \E\abs{ \Bu^H\BQ_j\By_j \By_j^H\BQ_j \Bu  - \frac{1}{n}\Bu^H\BQ_j\BOmega_j\BQ_j \Bu }^{2t} \lesssim_t \frac{1}{v^{6t} n^{4 \alpha t + 1}}. \numberthis \label{Eq_trace_2t_Z4j}
\end{align*}
Obviously, $Z_{3j}$ will not contribute to the RHS of \eqref{Eq_Bound_Gamma_1} since $Z_{3j} = 0$. Using $Z_{2j} = \E_{j}\eta_j\Bu^H\BQ_j\By_j \Ba_j^H\BQ_j \Bu $ and Lemma \ref{Lemm_Usedul_est}, we have 
\begin{align*}
    \E\left( \sum_{j=1}^n \E_{j-1}|Z_{2j}|^2 \right)^t &\leq \E\left(\sum_{j=1}^n \E_{j-1}\abs{\eta_j\Bu^H\BQ_j\By_j \Ba_j^H\BQ_j \Bu}^2 \right)^t \\
    &\lesssim_t \frac{1}{n^tv^{4t}}\E\left(\sum_{j=1}^n\E_{j-1}(\Bu^H\BQ_j^H \Ba_j\Ba_j^H\BQ_j \Bu) \right)^t \overset{(a)}{\lesssim_t} \frac{\log^{2t}(n)}{n^tv^{8t}}, \numberthis
\end{align*}
where step $(a)$ follows by the estimation in \eqref{Eq_useful_est_4}.
Similar to \eqref{Eq_trace_2t_Z4j}, we can get 
\begin{align*}
    \sum_{j=1}^n \E|Z_{2j}|^{2t}  & \lesssim_t \frac{1}{v^{2t}} \sum_{j=1}^n\E\abs{\Bu^H\BQ_j\By_j \Ba_j^H\BQ_j \Bu}^{2t} \lesssim_t \frac{1}{v^{4t} n^{2\alpha t + 2}} \sum_{j=1}^n\E\abs{ \Ba_j^H\BQ_j \Bu}^{2t} \\
    & \lesssim_t \frac{1}{v^{6t - 2} n^{2\alpha t + 2}} \sum_{j=1}^n\E\abs{ \Ba_j^H\BQ_j \Bu}^{2} \lesssim \frac{\log^2(n)}{v^{6t + 2} n^{2\alpha t + 2}}. \numberthis \label{Eq_Burk_Bound_2t_Z2j}
\end{align*}
The evaluation of the term related to $Z_{1j}$ is quite similar to that related to $Z_{2j}$, and is thus omitted. The following is obtained by gathering the bounds \eqref{Eq_Burk_bound_Z4j}-\eqref{Eq_Burk_Bound_2t_Z2j}
\begin{equation}
    \E \abs{\sum_{j=1}^n \Gamma_{1j}}^{2t} \lesssim_t \frac{1}{n^tv^{6t}} + \frac{1}{v^{6t} n^{4 \alpha t + 1}} +\frac{\log^{2t}(n)}{n^tv^{8t}} + \frac{\log^2(n)}{v^{6t + 2} n^{2\alpha t + 2}}.
\end{equation}
\par
Next, we handle $\mathbb{E}|\sum_j \Gamma_{2j}|^{2t}$. To this end, we decompose $\Gamma_{2j}$ as
\begin{align*}
    -\Gamma_{2j} &= [\E_j - \E_{j-1}](z \widetilde{q}_j\eta_j \Delta_j \Bu^H\BQ_j\Ba_j \By_j^H\BQ_j \Bu + z \widetilde{q}_j\eta_j \Delta_j \Bu^H\BQ_j\By_j \Ba_j^H\BQ_j \Bu + z \widetilde{q}_j\eta_j \Delta_j \Bu^H\BQ_j\Ba_j \Ba_j^H\BQ_j \Bu \\
    &+ z \widetilde{q}_j\eta_j \Delta_j \Bu^H\BQ_j\By_j \By_j^H\BQ_j \Bu) = \chi_{1j} + \chi_{2j} + \chi_{3j} + \chi_{4j}. \numberthis
\end{align*}
Define the set $A = \{\norm{\BSigma} \leq \log(n)\}$. Using $\E_{j-1}|(\E_j - \E_{j-1})a|^2 \leq 2\E_{j-1}|a|^2$, and the inequalities $|z\widetilde{q}_j| \lesssim  \norm{\BSigma}^2 / v$ and $|\eta_j| \lesssim v^{-1}$ in Lemma \ref{Lemm_Bound_on_qterms}, we have 
\begin{equation}
    \sum_{j=1}^n \E_{j-1} |\chi_{4j}|^2 = \sum_{j=1}^n \E_{j-1} \left[ |\chi_{4j}|^2 (\Ind_{A} + \Ind_{A^c}) \right] \lesssim \frac{\log^4(n)}{v^4}\sum_{j=1}^n \E_{j-1}\abs{\Delta_j}^2\abs{ \Bu^H\BQ_j\By_j \By_j^H\BQ_j \Bu}^2  + \sum_{j=1}^n  \E_{j-1} \left[|\chi_{4j}|^2 \Ind_{A^c} \right]. 
\end{equation}
According to the truncation $|X_{ij}| \leq n^{1/2 - \alpha} < \sqrt{n}$, we know that $\|\By_j\|^2 \lesssim n$,  $\norm{\BSigma} \lesssim n$, and $|\Delta_j| \lesssim (\norm{\Bxi_j}^2 + 1)  / v \lesssim n^2 / v$, which further give $|\chi_{4j}| \leq n^6/v^5$. As a result, the following is obtained by Hölder's inequality and Lemma \ref{Lemm_bound_Delta}
\begin{align*}
\label{Eq_bound_X4j}
    &\E \left[ \sum_{j=1}^n \E_{j-1} |\chi_{4j}|^2 \right]^t \lesssim_t \E \left[ \frac{\log^4(n)}{v^4}\sum_{j=1}^n \E_{j-1}\abs{\Delta_j}^2\abs{ \Bu^H\BQ_j\By_j \By_j^H\BQ_j \Bu}^2  + \frac{n^{12}}{v^{10}}\sum_{j=1}^n  \E_{j-1}  \Ind_{A^c}\right]^t \\
    &\lesssim_t \frac{\log^{4t}(n)}{v^{4t}} \E \left[ \sum_{j=1}^n \E_{j-1}^{4 / (4 + \varepsilon)}\abs{\Delta_j}^{2 + \varepsilon/2} \E_{j-1}^{\varepsilon / (4 + \varepsilon)}\abs{ \Bu^H\BQ_j\By_j \By_j^H\BQ_j \Bu}^{(8 + 2 \varepsilon) / \varepsilon} \right]^t  + \frac{n^{13t}\PP(A^c)}{v^{10t}} \\
    &\lesssim_t \frac{\log^{4t}(n)}{v^{10t}n^{8 \alpha t}} + \frac{n^{13t}\PP(A^c)}{v^{10t}} \lesssim_t \frac{\log^{4t}(n)}{v^{10t}n^{8 \alpha t}}. \numberthis
\end{align*}
By Lemma \ref{Lemm_P_bound_Sn}, the probability of $A^c$ is smaller than any polynomial of $n^{-1}$. From the above derivation, we can observe that if $|a| \leq \mathcal{P}$, where $\mathcal{P}$ is a polynomial in $n$, then we have $\mathbb{E} |a| \leq \mathbb{E} |a|\Ind_A + \mathcal{P} \PP(A^c) = \mathbb{E} |a|\Ind_A + o(K_l n^{-l})$ by splitting $1 = \Ind_A + \Ind_{A^c}$. Therefore, in the following discussion, if the evaluated term is upper bounded by a polynomial in $n$,  we always assume $\norm{\BSigma} \leq \log(n)$ for simplicity. The evaluation for $\sum_j \E |\chi_{4j}|^{2t}$ is given by
\begin{align*}
    \sum_{j=1}^n \E |\chi_{4j}|^{2t} & \lesssim_t \sum_{j=1}^n \E |z\widetilde{q}_j\eta_j \Delta_j \Bu^H\BQ_j\By_j \By_j^H\BQ_j \Bu|^{2t} \lesssim_t \frac{\log^{4t}(n)}{v^{4t}} \sum_{j=1}^n\E^{1/2} |\Delta_j|^{4t} \E^{1/2}|\Bu^H\BQ_j\By_j \By_j^H\BQ_j \Bu|^{4t} \\
    & \lesssim_t\frac{\log^{4t}(n) n}{v^{4t}} \sqrt{\frac{1}{v^{4t}}\left(\frac{1}{n^{8\alpha t + 1}} + \frac{1}{n^{4\alpha t + 2}} \right) \cdot \frac{1}{v^{8t}n^{8\alpha t + 2}}} \lesssim \frac{\log^{4t}(n)}{v^{10 t}} \left(\frac{1}{n^{8 \alpha t + 0.5}}+\frac{1}{n^{6\alpha t + 1}} \right). \numberthis \label{Eq_chi_4j_2t}
\end{align*}
% By using the same method, we can get
% \begin{align*}
%     &\E \left[ \sum_{j=1}^n \E_{j-1} |\chi_{3j}|^2 \right]^t \lesssim_t \frac{\log^{6t}(n)}{ v^{12 t} n^t}, ~~\sum_{j=1}^n \E |\chi_{3j}|^{2t} \lesssim_t \frac{\log^{4t + 2}(n)}{v^{10 t + 2} }\left(\frac{1}{n^{4\alpha t + 0.5}}+ \frac{1}{n^{2\alpha t + 1}} \right), \\
%     &\E \left[ \sum_{j=1}^n \E_{j-1} |\chi_{ij}|^2 \right]^t \lesssim_t \frac{\log^{4t}(n)}{v^{10 t} n^{6\alpha t} }, ~~ \sum_{j=1}^n \E |\chi_{ij}|^{2t} \lesssim_t \frac{\log^{4t}(n)}{v^{10 t}  }\left(\frac{1}{n^{6\alpha t + 0.5}} + \frac{1}{n^{4\alpha t + 1}}\right), i = 1, 2. \numberthis
% \end{align*}
Next, we evaluate the $\chi_{3j}$-related terms. By Lemma \ref{Lemm_Usedul_est}, we have 
\begin{align*}
    \E \left[\sum_{j=1}^n \E_{j-1} |\chi_{3j}|^2 \right]^t \lesssim_t \frac{\log^{4t}(n)}{v^{4t}} \E \left[\sum_{j=1}^n \E_{j-1}|\Delta_j|^2 \abs{\Bu^H\BQ_j\Ba_j\Ba_j^H\BQ_j\Bu}^2 \right]^t  \lesssim_t \frac{\log^{6t}(n)}{ v^{12 t} n^t}. \numberthis
\end{align*}
Define $\E_{-\Bx_j}(\cdot) = \E(\cdot | \mathscr{F}_{-\Bx_j})$ with $\mathscr{F}_{-\Bx_j}$ being the $\sigma$-algebra generated by $\{\Bx_{j'} : j' \leq n, j' \neq j\}$. The evaluation for $\sum_{j} \E |\chi_{3j}|^{2t}$ can be given by
\begin{equation}
    \sum_{j=1}^n \E |\chi_{3j}|^{2t} \lesssim_t \frac{\log^{4t}(n)}{v^{4t}} \sum_{j=1}^n \E\left[ \left(\E_{-\Bx_j}|\Delta_j|^{2t} \right) \abs{\Bu^H\BQ_j\Ba_j\Ba_j^H\BQ_j\Bu}^{2t}\right] \lesssim_t \frac{\log^{4t+2}(n)}{v^{10t+2}} \left(\frac{1}{n^{4\alpha t + 1}} + \frac{1}{n^{2\alpha t + 2}}\right).
\end{equation}
We then handle the $\chi_{2j}$-related terms. By the independence between $\By_j$ and $\BQ_j$, Lemma \ref{Lemm_bound_Delta}, and Lemma \ref{Lemm_trace}, we have 
\begin{align*}
    &\E_{j-1}\abs{\chi_{2j}}^2  \lesssim \frac{\log^{4}(n)}{v^{4}} \E_{j-1}  \abs{\Delta_j}^2 |\Bu^H\BQ_j \By_j|^2 \abs{\Ba_j^H \BQ_j \Bu}^2  = \frac{\log^{4}(n)}{v^{4}} \E_{j-1} \left\{ \E_{-\Bx_j}\left[\abs{\Delta_j}^2 |\Bu^H\BQ_j \By_j|^2\right] \abs{\Ba_j^H \BQ_j \Bu}^2 \right\} \\
    & \leq \frac{\log^{4}(n)}{v^{4}} \E_{j-1} \left[ \left(\E_{-\Bx_j}^{4 / (4 + \varepsilon)}\abs{\Delta_j}^{2 + \varepsilon/2} \right) \left( \E_{-\Bx_j}^{\varepsilon / (4 + \varepsilon)}|\Bu^H\BQ_j \By_j|^{(8 + 2 \varepsilon) / \varepsilon} \right) \abs{\Ba_j^H \BQ_j \Bu}^2 \right] \lesssim \frac{\log^{4}(n)}{v^{8}n^{6 \alpha + 1}} \E_{j-1} \abs{\Ba_j^H \BQ_j \Bu}^2, \numberthis
\end{align*}
which implies
\begin{align*}
    \E \left[ \sum_{j=1}^n \E_{j-1}\abs{\chi_{2j}}^2 \right]^t \lesssim_t \frac{\log^{4t}(n)}{v^{8t} n^{(6\alpha + 1)t}} \E\left[ \sum_{j=1}^n \E_{j-1}  \Bu^H \BQ_j^H \Ba_j \Ba_j^H \BQ_j \Bu \right]^t \lesssim_t \frac{\log^{6t}(n)}{v^{12t}n^{(6\alpha + 1)t}}. \numberthis
\end{align*}
Similar to \eqref{Eq_chi_4j_2t}, the evaluation for $\sum_j \E |\chi_{2j}|^{2t}$ is given by
\begin{align*}
    &\sum_{j=1}^n \E|\chi_{2j}|^{2t} \lesssim_t \frac{\log^{4t}(n)}{v^{4t}} \sum_{j=1}^n \E \left\{ \left(\E^{1/2}_{-\Bx_j}|\Delta_j|^{4t}\right) \left(\E^{1/2}_{-\Bx_j}|\Bu^H\BQ_j\By_j|^{4t} \right) \abs{\Ba_j^H\BQ_j\Bu}^{2t} \right\} \\
    &\lesssim_t \frac{\log^{4t}(n)}{v^{8t}} \left(n^{-(6\alpha t + 1.5)} + n^{-(4\alpha t + 2)} \right) \sum_{j=1}^n \E |\Ba_j^H\BQ_j \Bu|^{2t} \lesssim_t \frac{\log^{4t + 2}(n)}{v^{10t + 2}} \left(n^{-(6\alpha t + 1.5)} + n^{-(4\alpha t + 2)} \right). \numberthis
\end{align*}
The evaluation for the $\chi_{1j}$-related terms is quite similar to that for $\chi_{2j}$. With the above analysis, it can be verified that when $v$ is close to the real axis, the bound of $\E|\Bu^H[\BQ - \E \BQ]\Bu|^{2t}$ is dominated by $\log^{6t}(n)/(v^{12t}n^{2 \alpha t + 0.5})$, which completes the proof of \eqref{Eq_Part_I}.
\paragraph*{Part II} \label{Para_Part_II}
In this part, we prove that 
\begin{equation}
\label{Eq_Part_II}
    \abs{\Bu^H (\E \BQ - \BR) \Bu} \lesssim \frac{1}{n^{\alpha}v^7}.
\end{equation}
We note that the terms $\Bu$, $\E \BQ$, and $\BR$ are deterministic. By using the resolvent identity $\BA - \BB = \BA(\BB^{-1} - \BA^{-1})\BB$, we have 
\begin{align*}
\Bu^H(\BQ - \BR)\Bu &= \Bu^H\BQ\left(\sum_{j=1}^n \frac{\BOmega_j}{n}(-z\widetilde{\kappa}_j) -z \BA\widetilde{\BG}\BA^H - \sum_{j=1}^n \Bxi_j\Bxi_j^H \right) \BR \Bu \\
&= \sum_{j=1}^n  -z\widetilde{\kappa}_j \frac{\Bu^H\BQ\BOmega_j\BR\Bu}{n} + \frac{\Bu^H\BQ\Ba_j\Ba_j^H\BR\Bu}{1 + \kappa_j} - \frac{\Bu^H\BQ_j\Bxi_j\Bxi_j^H\BR\Bu}{1 + \Bxi_j^H\BQ_j\Bxi_j} \\
&= \sum_{j=1}^n  -z\widetilde{\kappa}_j \frac{\Bu^H\BQ\BOmega_j\BR\Bu}{n} + \frac{\Bu^H\BQ_j\Ba_j\Ba_j^H\BR\Bu}{1 + \kappa_j} - \frac{\Bu^H\BQ_j\Bxi_j\Bxi_j^H\BQ_j\Ba_j\Ba_j^H\BR\Bu}{(1 + \kappa_j)(1 + \Bxi_j^H\BQ_j\Bxi_j)} - \frac{\Bu^H\BQ_j\Bxi_j\Bxi_j^H\BR\Bu}{1 + \Bxi_j^H\BQ_j\Bxi_j} \\
&= \sum_{j=1}^n W_{1j} + W_{2j} + W_{3j} + W_{4j}, \numberthis
\end{align*}
where
\begin{align*}
  W_{1j} &= -z\widetilde{\kappa}_j \frac{\Bu^H\BQ\BOmega_j\BR\Bu}{n} - \frac{\Bu^H\BQ_j\By_j\By_j^H\BR\Bu}{1 + \Bxi_j^H\BQ_j\Bxi_j}, ~~ W_{2j} = \frac{(\By_j^H\BQ_j\By_j - \kappa_j)\Bu^H\BQ_j\Ba_j\Ba_j^H\BR\Bu}{(1 + \kappa_j)(1 + \Bxi_j^H\BQ_j\Bxi_j)}, \\
  W_{3j} &= - \frac{\Bu^H\BQ_j\By_j\By_j^H\BQ_j\Ba_j\Ba_j^H\BR\Bu}{(1 + \kappa_j)(1 + \Bxi_j^H\BQ_j\Bxi_j)}, \\
  W_{4j} &= - \frac{\Bu^H\BQ_j\Ba_j\By_j^H\BQ_j\Ba_j\Ba_j^H\BR\Bu + \Bu^H\BQ_j\By_j\Ba_j^H\BQ_j\Ba_j\Ba_j^H\BR\Bu}{(1 + \kappa_j)(1 + \Bxi_j^H\BQ_j\Bxi_j)} \\
  &+ \frac{\Ba_j^H\BQ_j\By_j \Bu^H\BQ_j\Ba_j\Ba_j^H\BR\Bu + \By_j^H\BQ_j\Ba_j \Bu^H\BQ_j\Ba_j\Ba_j^H\BR\Bu}{(1 + \kappa_j)(1 + \Bxi_j^H\BQ_j\Bxi_j)} \\
  &- \frac{\Bu^H\BQ_j\By_j\Ba_j^H\BR\Bu + \Bu^H\BQ_j\Ba_j\By_j^H\BR\Bu}{1 + \Bxi_j^H\BQ_j\Bxi_j}. \numberthis
\end{align*}
First, we evaluate $\sum_j \E W_{1j}$. Express $\E W_{1j}$ as
\begin{align*}
\E W_{1j} &=\E \left[ (-z \widetilde{q}_j - \eta_j) \left(\frac{\Bu^H\BQ_j\BOmega_j\BR\Bu}{n} - \Bu^H\BQ_j\By_j\By_j^H\BR\Bu\right) \right] + \E\left[ (z\widetilde{q}_j -z\widetilde{\kappa}_j) \frac{\Bu^H\BQ_j\BOmega_j\BR\Bu}{n} \right] \\
&-z\widetilde{\kappa}_j \E \frac{\Bu^H(\BQ - \BQ_j)\BOmega_j\BR\Bu}{n} = \zeta_{1j} + \zeta_{2j} + \zeta_{3j}. \numberthis
\end{align*}
By the Cauchy-Schwarz inequality and Lemma \ref{Lemm_trace}, we have 
\begin{align*}
   \sum_{j=1}^n|\zeta_{1j}| &\leq \sum_{j=1}^n\E^{\frac{1}{2}}\abs{z\widetilde{q}_j \eta_j{\Delta_j}}^2 \E^{\frac{1}{2}}\abs{\frac{1}{n}\Bu^H\BQ_j\BOmega_j\BR\Bu - \By_j^H\BR\Bu\Bu^H\BQ_j\By_j }^2  \lesssim n \cdot \frac{\log^2(n)}{v^3n^{1/2}} \frac{1}{v^2 n} = \frac{\log^2(n)}{v^5 n^{1/2}}. \numberthis
\end{align*}
Here, we have utilized the result $\norm{\BR} \leq 1 / v$, which can be obtained by the properties of the matrix-valued Stieltjes transform \cite[Proposition 2.2]{hachem2007deterministic}. By identity $\E(a - \E a)b = \E a(b - \E b)$ and Hölder’s inequality, the term $|\zeta_{2j}|$ can be bounded as 
\begin{align*}
    \sum_{j=1}^n |\zeta_{2j}| &\leq \sum_{j=1}^n\frac{1}{n}\E |z\widetilde{q}_j|\abs{\Bu^H (\BQ_j - \E \BQ_j) \BOmega_j \BR \Bu}  \lesssim \sum_{j=1}^n\frac{\log^2(n)}{nv} \E^{1/2q} \abs{\Bu^H (\BQ_j - \E \BQ_j) \BOmega_j \BR \Bu}^{2q} \overset{(a)}{\lesssim} \frac{\log^5(n)}{v^{7}n^{\alpha  + 1/4q}} \lesssim \frac{1}{v^{7}n^{ \alpha}}, \numberthis
\end{align*}
where step $(a)$ follows from \eqref{Eq_Part_I} and $q$ is a large enough given number. By using the resolvent decomposition in \eqref{Eq_Wood_resolvent_Id}, $|\zeta_{3j}|$ can be bounded as
\begin{align*}
|\zeta_{3j}| \leq \frac{|z\widetilde{\kappa}_j|}{n} \E \abs{z\widetilde{q}_j \Bu^H\BQ_j\Bxi_j\Bxi_j^H\BQ_j\BOmega_j\BR\Bu} \lesssim \frac{\log^4(n)}{nv^2}\E^{\frac{1}{2}}|\Bu^H\BQ_j\Bxi_j|^2 \E^{\frac{1}{2}}|\Bxi_j^H\BQ_j\BOmega_j\BR\Bu|^2. \numberthis
\end{align*}
By summing over the subscript $j$ and utilizing Lemma \ref{Lemm_Usedul_est}, we can get
\begin{align*}
\sum_{j=1}^n|\zeta_{3j}| \lesssim \frac{\log^4(n)}{nv^2}\left( \sum_{j=1}^n\E|\Bu^H\BQ_j\Bxi_j|^2 \right)^{1/2} \left( \sum_{j=1}^n\E|\Bxi_j^H\BQ_j\BOmega_j\BR\Bu|^2 \right)^{1/2} \lesssim \frac{\log^6(n)}{v^7 n}. \numberthis
\end{align*}
Hence, we have 
\begin{equation}
\label{Eq_Part_II_W1j}
    \abs{\sum_{j=1}^n \E W_{1j}} \lesssim \frac{\log^2(n)}{v^5n^{1/2}} + \frac{1}{v^7n^{\alpha}} + \frac{\log^6(n)}{v^7n}.
\end{equation}
\par
Next, we handle $|\sum_{j=1}^n \E W_{2j}|$. By inequality $|1 + \kappa_j| \lesssim 1 / v$ in Lemma \ref{Lemm_Bound_on_qterms} and the Cauchy-Schwarz inequality, we get
\begin{align*}
\sum_{j=1}^n |\E W_{2j}| &\leq \sum_{j=1}^n  \E \abs{ \frac{(\By_j^H\BQ_j\By_j - \kappa_j)\Bu^H\BQ_j\Ba_j\Ba_j^H\BR\Bu}{(1 + \kappa_j)(1 + \Bxi_j^H\BQ_j\Bxi_j)}} \lesssim \sum_{j=1}^n \frac{\log^2(n)}{v^2} |\Ba_j^H\BR\Bu|\E^{1/2}|\Bu^H\BQ_j\Ba_j|^2 \E^{1/2}\abs{\By_j^H\BQ_j\By_j - \kappa_j}^2 \\
&\lesssim \sum_{j=1}^n \frac{\log^2(n)}{v^2} \left( |\Bu^H\BR^H\Ba_j\Ba_j^H\BR\Bu +  \E \Bu^H\BQ_j\Ba_j\Ba_j^H\BQ_j^H\Bu \right) \E^{1/2}\abs{\By_j^H\BQ_j\By_j - \kappa_j}^2. \numberthis
\end{align*}
To evaluate $\E\abs{\By_j^H\BQ_j\By_j - \kappa_j}^2$, we write 
\begin{equation}
  \By_j^H\BQ_j\By_j - \kappa_j = \By_j^H\BQ_j\By_j - \frac{\Tr \BOmega_j \BQ_j}{n} + \frac{\Tr \BOmega_j (\BQ_j - \BQ) }{n}+ \frac{\Tr \BOmega_j (\BQ - \E \BQ)}{n}. \label{Eq_Diff_yQy_to_kappa}
\end{equation}
The third term on the RHS of \eqref{Eq_Diff_yQy_to_kappa} can be written as the sum of a martingale difference sequence 
\begin{equation}
  \frac{1}{n} \Tr \BOmega_j (\BQ - \E \BQ) = \sum_{j=1}^n \frac{1}{n} \Tr \BOmega_j [\E_j - \E_{j-1}] (\BQ - \BQ_j) = \sum_{j=1}^n [\E_j - \E_{j-1}] \frac{1}{n} z\widetilde{q}_j{\Bxi_j^H\BQ_j\BOmega_j\BQ_j\Bxi_j}.
\end{equation}
Using the fact $|z\widetilde{q}_j{\Bxi_j^H\BQ_j\BOmega_j\BQ_j\Bxi_j}| \leq \norm{\BOmega_j}/v$ and Lemma \ref{Lemm_Burk_Ineq}, we can obtain $
  \E \abs{\frac{1}{n} \Tr \BOmega_j (\BQ - \E \BQ)}^2 \lesssim n^{-1}v^{-2}
$. Thus, we get the bound $\E\abs{\By_j^H\BQ_j\By_j - \kappa_j}^2 \lesssim n^{-1}v^{-2}$, which implies 
\begin{align*}
\sum_{j=1}^n |\E W_{2j}| &\lesssim  \frac{\log^2(n)}{v^3n^{1/2}} \sum_{j=1}^n\left( |\Bu^H\BR^H\Ba_j\Ba_j^H\BR\Bu +  \E \Bu^H\BQ_j\Ba_j\Ba_j^H\BQ_j^H\Bu \right)  \lesssim \frac{\log^2(n)}{v^3 n^{1/2}} \left( \frac{1}{v^2} + \frac{\log^{2}(n)}{v^4} \right) \lesssim \frac{\log^4(n)}{n^{1/2}v^7}. \numberthis \label{Eq_Part_II_W2j}
\end{align*}
The evaluation for $\sum_j |\E W_{3j}|$ is given by
\begin{align*}
\sum_{j=1}^n |\E W_{3j}| &\leq \sum_{j=1}^n\E \abs{ \frac{\Bu^H\BQ_j\By_j\By_j^H\BQ_j\Ba_j\Ba_j^H\BR\Bu}{(1 + \kappa_j)(1 + \Bxi_j^H\BQ_j\Bxi_j)}} \lesssim \frac{\log^2(n)}{v^2} \sum_{j=1}^n |\Ba_j^H\BR\Bu| \E^{1/2}|\Bu^H\BQ_j\By_j|^2 \E^{1/2} |\By_j^H\BQ_j\Ba_j|^2 \\
& \lesssim    \frac{\log^2(n)}{nv^4} \sum_{j=1}^n|\Ba_j^H\BR\Bu| \leq \frac{\log^2(n)}{nv^4} \cdot \sqrt{n} \cdot \left(\sum_{j=1}^n \Bu^H\BR^H \Ba_j\Ba_j^H\BR\Bu \right)^{1/2}\lesssim  
 \frac{\log^2(n)}{n^{1/2}v^5}. \numberthis \label{Eq_Part_II_W3j}
\end{align*}
\par
To estimate $|\sum_j\E W_{4j}|$, we write $W_{4j} = P_{j}/[(1 + \kappa_j)(1 + \Bxi_j^H\BQ_j\Bxi_j)] + z\widetilde{q}_jQ_j$, where 
\begin{align*}
    P_j &= - \Bu^H\BQ_j\Ba_j\By_j^H\BQ_j\Ba_j\Ba_j^H\BR\Bu - \Bu^H\BQ_j\By_j\Ba_j^H\BQ_j\Ba_j\Ba_j^H\BR\Bu  + \Ba_j^H\BQ_j\By_j \Bu^H\BQ_j\Ba_j\Ba_j^H\BR\Bu \\
    &+ \By_j^H\BQ_j\Ba_j \Bu^H\BQ_j\Ba_j\Ba_j^H\BR\Bu,  \\
    Q_j &= \Bu^H\BQ_j\By_j\Ba_j^H\BR\Bu +\Bu^H\BQ_j\Ba_j\By_j^H\BR\Bu. \numberthis
\end{align*}
It can be verified that $\E P_j =\E Q_j = 0$ and 
% \begin{align*}
%    \sum_{j=1}^n \E|P_j|^2 &\lesssim \sum_{j=1}^n  \frac{1}{n} \Big[v^{-2}\E|\Bu^H\BQ_j\Ba_j\Ba_j^H\BR\Bu|^2 + v^{-2} \E|\Ba_j^H\BQ_j\Ba_j\Ba_j^H\BR\Bu|^2 + v^{-2}\E|\Bu^H\BQ_j\Ba_j\Ba_j^H\BR\Bu|^2 \\
%     &+ v^{-2} \E |\Bu^H\BQ_j\Ba_j\Ba_j^H\BR\Bu|^2  + \left(v^{-2} + 1\right)^2 \left( v^{-2}\E |\Ba_j^H\BQ_j\BR\Bu|^2 + v^{-4} \E |\Bu^H\BQ_j\Ba_j|^2  \right) \Big] \\
%     & \lesssim \frac{1}{nv^6} + \frac{\log^2(n)}{nv^{10}}. \numberthis \\
%     \sum_{j=1}^n \E|Q_j|^2 &\lesssim \left( v^{-2}\E |\Ba_j^H\BQ_j\BR\Bu|^2 + v^{-4} \E |\Bu^H\BQ_j\Ba_j|^2  \right) 
% \end{align*}
\begin{align*}
   \sum_{j=1}^n \E|P_j|^2 &\lesssim \sum_{j=1}^n  \frac{1}{n} \Big[v^{-2}\E|\Bu^H\BQ_j\Ba_j\Ba_j^H\BR\Bu|^2 + v^{-2} \E|\Ba_j^H\BQ_j\Ba_j\Ba_j^H\BR\Bu|^2 + v^{-2}\E|\Bu^H\BQ_j\Ba_j\Ba_j^H\BR\Bu|^2 \\
    &+ v^{-2} \E |\Bu^H\BQ_j\Ba_j\Ba_j^H\BR\Bu|^2  \Big]  \lesssim \frac{1}{nv^6}. \numberthis \\
    \sum_{j=1}^n \E|Q_j|^2 &\lesssim \frac{1}{n}\sum_{j=1}^n \left( v^{-2}\E |\Ba_j^H\BR\Bu|^2 + v^{-2} \E |\Bu^H\BQ_j\Ba_j|^2  \right)  \lesssim \frac{\log^2(n)}{n v^6}.
\end{align*}
By writing $-z \widetilde{q}_j = \eta_j + z\widetilde{q}_j\eta_j \Delta_j$, we can obtain
\begin{align*}
\abs{ \sum_{j=1}^n \E W_{4j}} &\leq \sum_{j=1}^n \E \abs{\frac{P_j\Delta_j\eta_j z\widetilde{q}_j}{1 + \kappa_j}} + \E \abs{Q_j\Delta_j\eta_j z\widetilde{q}_j} \lesssim \frac{\log^2(n)}{v^3} \sum_{j=1}^n \E^{\frac{1}{2}}|\Delta_j|^2 \E^{\frac{1}{2}} |P_j|^2 + \frac{\log^2(n)}{v^2} \sum_{j=1}^n \E^{\frac{1}{2}}|\Delta_j|^2 \E^{\frac{1}{2}} |Q_j|^2 \\
&\lesssim \frac{\log^2(n)}{v^4}\left( \sum_{j=1}^n \E|P_j|^2 \right)^{\frac{1}{2}} + \frac{\log^2(n)}{v^3}\left( \sum_{j=1}^n \E|Q_j|^2 \right)^{\frac{1}{2}} \lesssim\frac{\log^3(n)}{n^{\frac{1}{2}}v^7}. \numberthis \label{Eq_Part_II_W4j}
\end{align*}
Combining the results from \eqref{Eq_Part_II_W1j}, \eqref{Eq_Part_II_W2j}, \eqref{Eq_Part_II_W3j}, and \eqref{Eq_Part_II_W4j}, we can prove \eqref{Eq_Part_II}.
\paragraph*{Part III} In this part, we will prove
\begin{equation}
\label{Eq_Part_III}
  \abs{\Bu^H [\BR - \BTheta] \Bu} \lesssim \frac{\mathcal{U}(z)}{n^{\alpha}},
\end{equation}
where $\mathcal{U}(z) = (K_2 + x^2 + v^2)^{9} / v^{29}$.
Before devoting to the detail of the proof, we first evaluate two differences $-z \widetilde{\kappa}_j - z [\widetilde{\BR}]_{j, j}$, and $\kappa_j - \Tr \BOmega_j \BR / n$, which are given in the following proposition
\begin{proposition}
\label{Prop_Co_Res_trace}
    For $z = x + \jmath v$ with $v \leq K_{\mathrm{im}}$, we have
    \begin{align*}
        \abs{-z\widetilde{\kappa}_j - (-z[\widetilde{\BR}]_{j, j})} \lesssim \frac{1}{n^{\alpha}v^9}, ~~ \abs{\kappa_j - \frac{1}{n}\Tr \BOmega_j \BR} \lesssim \frac{1}{n^{\alpha}v^7}. \numberthis
    \end{align*}
\end{proposition}
\textit{Proof: } We first handle the first inequality. By the Woodbury matrix identity, we have 
\begin{equation}
[\widetilde{\BR}]_{j, j} =-\frac{1}{z(1 + \kappa_j)} + \frac{\Ba_j^H\BR\Ba_j}{z(1 + \kappa_j)^2}. \numberthis
\end{equation}
Then, the difference $-z \widetilde{\kappa}_j - (-z [\widetilde{\BR}]_{j, j})$ can be evaluated as
\begin{align*}
&-z \widetilde{\kappa}_j - (-z [\widetilde{\BR}]_{j, j}) = \E \left[\frac{1}{1 + \Bxi_j^H\BQ_j\Bxi_j} - \eta_j \right] + \E \left[\eta_j - \frac{1}{1 + \alpha_{j,j}} + \frac{\Ba_j^H\BQ\Ba_j}{(1 + \alpha_{j, j})^2} \right] + \E \left[ \frac{1}{1 + \alpha_{j,j}} - \frac{1}{1 + \kappa_j} \right] \\
&+ \frac{\Ba_j^H(\BR - \E \BQ)\Ba_j }{(1 + \kappa_j)^2} + \E \left[ \Ba_j^H \BQ \Ba_j\left( \frac{1 }{(1 + \kappa_j)^2} - \frac{1}{(1 + \alpha_{j,j})^2} \right)\right] = X_{1j} + X_{2j} + X_{3j} + X_{4j} + X_{5j}. \numberthis \label{Eq_kappa_j_wR_jj}
\end{align*}
Next, we evaluate the term $ X_{ij}$ for $i \in \{1, \ldots, 5\}$. For $X_{1j}$,  it is straightforward to show that $|X_{1j}| \leq \mathbb{E} |z\widetilde{q}_j \eta_j \Delta_j| \lesssim \log^2(n)v^{-3}n^{-1/2}$. Using \eqref{Eq_Wood_resolvent_Id} and $-z\widetilde{q}_j = \eta_j + z\widetilde{q}_j \Delta_j \eta_j $, we have 
\begin{align*}
    X_{2j} &= \E \left[\eta_j - \frac{1}{1 + \alpha_{j,j}} + \frac{\Ba_j^H\BQ_j\Ba_j - (\eta_j + z\widetilde{q}_j \Delta_j \eta_j) \Ba_j^H \BQ_j \Bxi_j \Bxi_j^H \BQ_j\Ba_j }{(1 + \alpha_{j, j})^2} \right] \\
    & = \E \left[ \eta_j - \frac{1}{1 + \alpha_{j,j}} + \frac{\Ba_j^H\BQ_j\Ba_j - \eta_j \Ba_j^H \BQ_j \Ba_j \Ba_j^H \BQ_j\Ba_j }{(1 + \alpha_{j, j})^2} \right] \\
    &- \E \left[\frac{\eta_j(\Ba_j^H\BQ_j\Ba_j\By_j^H\BQ_j\Ba_j + \Ba_j^H\BQ_j\By_j\Ba_j^H\BQ_j\Ba_j + \Ba_j^H\BQ_j\By_j\By_j^H\BQ_j\Ba_j)}{(1 + \alpha_{j,j})^2} \right] \\
    & -\E \left[\frac{ z\widetilde{q}_j \Delta_j \eta_j \Ba_j^H \BQ_j \Bxi_j \Bxi_j^H \BQ_j\Ba_j }{(1 + \alpha_{j, j})^2} \right] = X_{2j, 1} + X_{2j, 2} + X_{2j, 3}. \numberthis \label{Eq_kappa_j_wR_jj_X_2j}
\end{align*}
It can be shown that $X_{2j, 1} = 0$. As $\By_j$ is independent of $\eta_j$, $\alpha_j$, and $\BQ_j$, we have $|X_{2j, 2}| \lesssim n^{-1}v^{-5}$. By the same method as in \cite[Lemma 2.6]{SILVERSTEIN1995Empirical}, we have $|z\widetilde{q}_j\Ba_j^H \BQ_j \Bxi_j \Bxi_j^H \BQ_j\Ba_j| \leq \norm{\Ba_j}^2 / v \leq \norm{\BA}^2 / v$ and 
\begin{align*}
    |X_{2j, 3}| \lesssim \frac{1}{v^4} \E^{1/2} |\Delta_j|^2 \lesssim \frac{1}{v^5 n^{1/2}}. \numberthis
\end{align*}
Hence, we get $|X_{2j}| \lesssim v^{-5}n^{-1/2}$. The bounds $|X_{3j}| \lesssim v^{-3}n^{-1}$ and $|X_{5j}| \lesssim v^{-7}n^{-1}$ can be obtained by Lemma \ref{Lemm_Bound_on_qterms}. According to \eqref{Eq_Part_II} and  Lemma \ref{Lemm_Bound_on_qterms}, we have $|X_{4j}| \lesssim n^{-\alpha}v^{-9}$. Therefore, it can be concluded that
$
    |-z \widetilde{\kappa}_j - (-z [\widetilde{\BR}]_{j, j})|\lesssim {v^{-9} n^{-\alpha}},
$
which proves the first inequality.
\par
For the second inequality, we write the eigen-decomposition of $\BOmega_j$ as $\BOmega_j = \sum_{i = 1}^p \lambda^j_i\Bu^j_i(\Bu^j_i)^H  $. Then, the following holds by \eqref{Eq_Part_II}
\begin{align*}
    \abs{\kappa_j - \frac{1}{n} \Tr \BOmega_j \BR} \leq \frac{1}{n}\sum_{i=1}^p \lambda^j_i\abs{ (\Bu^j_i)^H[\E \BQ - \BR] \Bu^j_i} \lesssim \frac{1}{v^{7}n^{\alpha}}. \numberthis
\end{align*}
Therefore, we have completed the proof. \qed
\par
Denote $\Bkappa = [\kappa_1, \ldots, \kappa_n]^T$, $\widetilde{\Bkappa} = [\widetilde{\kappa}_1, \ldots, \widetilde{\kappa}_n]^T$, $\Bdelta = [\delta_1, \ldots, \delta_n]^T$, and $\widetilde{\delta}_j = [\widetilde{\delta}_1, \ldots, \widetilde{\delta}_n]^T$. Given $\BR$ is a function of $(\Bkappa, \widetilde{\Bkappa})$, and $\BTheta $ is a function of $(\Bdelta, \widetilde{\Bdelta})$, it is essential to evaluate $\|\Bkappa - \Bdelta\|_{\infty}$ and $\| \widetilde{\Bkappa} - \widetilde{\Bdelta} \|_{\infty}$ in order to analyze $\Bu^H[\BR - \BTheta]\Bu $. By writing $\kappa_j - \delta_j = \kappa_j - \frac{1}{n} \Tr \BOmega_j \BR + \frac{1}{n} \Tr \BOmega_j \BR - \delta_j$, and then applying the resolvent identity $ \BR - \BTheta = \BR(\BTheta^{-1}- \BR^{-1})\BTheta $, we can obtain the following system of equations
\begin{equation}
  \begin{bmatrix}
   \Bkappa - \Bdelta \\
   z\widetilde{\Bkappa} - z\widetilde{\Bdelta}  
  \end{bmatrix} = \begin{bmatrix}
    \Bd \\
    z\widetilde{\Bd}  
   \end{bmatrix} + \begin{bmatrix}
    z^2 \BPsi_{\kappa, \delta} & \BGamma_{\kappa, \delta} \\
    z^2\widetilde{\BGamma}_{\kappa, \delta} & z^2 \BPhi_{\kappa, \delta}
   \end{bmatrix} \begin{bmatrix}
    \Bkappa - \Bdelta \\
    z\widetilde{\Bkappa} - z\widetilde{\Bdelta}  
   \end{bmatrix}, \label{Eq_Maitrx_kappa_delta}
\end{equation} 
where $[\Bd]_j = \frac{1}{n} \Tr \BOmega_j (\E \BQ - \BR), ~~ [\widetilde{\Bd}]_j = \widetilde{\kappa}_{j} - [\widetilde{\BR}]_{j,j}$, and the matrix blocks are defined as
\begin{align*}
&[\BPsi_{\kappa, \delta}]_{j, l} = \frac{1}{n}[\widetilde{\BF} \BA^H \BTheta \BOmega_{j} \BR \BA \widetilde{\BG}]_{l, l},  ~~ [\BGamma_{\kappa, \delta}]_{j, l} = \frac{1}{n^2}\Tr \BOmega_{j}\BR \BOmega_l \BTheta, \\
&[\widetilde{\BGamma}_{\kappa, \delta}]_{j, l} = [\widetilde{\BR}]_{j, l} [\widetilde{\BTheta}]_{l, j}, ~~[\BPhi_{\kappa, \delta}]_{j, l} = \frac{1}{n} [\widetilde{\BR}\BA^H \BG\BOmega_l \BF \BA \widetilde{\BTheta}]_{j, j}. \numberthis \label{Eq_Def_of_Gamma_Phi_Psi_delta_kappa}
% & \BF = \left[-z \left( \BI_p + \textstyle\sum_{j=1}^n {\BOmega_j \widetilde{\delta}_j}\big/{n} \right) \right]^{-1}, ~~ \widetilde{\BF} = \left[ -z(\BI_n + \BD_{\delta}) \right]^{-1}, \\
% &\BG = \left[-z\left( \BI_p + \textstyle\sum_{j=1}^n {\BOmega_j \widetilde{\kappa}_j}\big/{n} \right) \right]^{-1}, ~~ \widetilde{\BG} = \left[ -z(\BI_n + \BD_{\kappa}) \right]^{-1}. 
\end{align*}
According to Proposition \ref{Prop_Co_Res_trace}, we know $\|\Bd\|_{\infty} \leq K v^{-9}n^{-\alpha}$ and $\|z\widetilde{\Bd}\|_{\infty} \leq K v^{-7}n^{-\alpha}$ for some constant $K$. To control the difference vector $(\Bkappa - \Bdelta, \widetilde{\Bkappa} - \widetilde{\Bdelta})$, we need to solve the above system of equations. By the Cauchy-Schwarz inequality, we have 
\begin{align*}
    \abs{[\BPi_{\kappa, \delta}]_{j, l}} \leq \sqrt{[\BPi_{ \kappa}]_{j, l} [\BPi_{ \delta}]_{j, l}}, \numberthis
\end{align*}
where $\BPi$ can be any of $\{ \BPsi, \BGamma, \widetilde{\BGamma}, \BPhi\}$, which are given by
\begin{align*}
&[\BPsi_{\kappa}]_{j, l} = \frac{1}{n}[\widetilde{\BG}^H \BA^H \BR^H \BOmega_{j} \BR \BA \widetilde{\BG}]_{l, l},  ~~ [\BGamma_{\kappa}]_{j, l} = \frac{1}{n^2}\Tr \BOmega_{j}\BR \BOmega_l \BR^H, \\
&[\widetilde{\BGamma}_{\kappa}]_{j, l} = [\widetilde{\BR}]_{j, l} [\widetilde{\BR}^H]_{l, j}, ~~[\BPhi_{\kappa}]_{j, l} = \frac{1}{n} [\widetilde{\BR}\BA^H \BG\BOmega_l \BG^H \BA \widetilde{\BR}^H]_{j, j},\\
&[\BPsi_{\delta}]_{j, l} = \frac{1}{n}[\widetilde{\BF} \BA^H \BTheta \BOmega_{j} \BTheta^H \BA \widetilde{\BF}^H]_{l, l},  ~~ [\BGamma_{\delta}]_{j, l} = \frac{1}{n^2}\Tr \BOmega_{j}\BTheta^H \BOmega_l \BTheta, \\
&[\widetilde{\BGamma}_{\delta}]_{j, l} = [\widetilde{\BTheta}^H]_{j, l} [\widetilde{\BTheta}]_{l, j}, ~~[\BPhi_{\delta}]_{j, l} = \frac{1}{n} [\widetilde{\BTheta}^H\BA^H \BF^H \BOmega_l \BF \BA \widetilde{\BTheta}]_{j, j}. \numberthis \label{Eq_Def_of_Gamma_Phi_Psi_delta_kappa_split}
\end{align*}
We note that $\BT^H = \BT(z^*)$ for $\BT \in \{\BF, \BTheta, \widetilde{\BF}, \widetilde{\BTheta},\BG, \BR, \widetilde{\BG}, \widetilde{\BR} \}$. Intuitively, by Lemma \ref{Lemm_matrix_Spr_maxrow_Norm}, to control the error, it is sufficient to study the properties of $\BPi_\kappa$ and $\BPi_\delta$, which are more tractable since each term of $\BPi_\kappa$ and $\BPi_{\delta}$ only involves the terms with respect to $(\Bkappa, \widetilde{\Bkappa})$ and $(\Bdelta, \widetilde{\Bdelta})$.
% Similar to \eqref{Eq_Maitrx_kappa_delta}, we define $\BC_{\Bdelta}$ and $\BC_{\kappa}$ as 
% \begin{equation}
%     \BC_{\kappa} = \begin{bmatrix}
%     |z|^2 \BPsi_{\kappa} & \BGamma_{\kappa} \\
%     |z|^2\widetilde{\BGamma}_{\kappa} & |z|^2 \BPhi_{\kappa}
%    \end{bmatrix}, ~~~ \BC_{\delta} = \begin{bmatrix}
%     |z|^2 \BPsi_{\delta} & \BGamma_{\delta} \\
%     |z|^2\widetilde{\BGamma}_{\delta} & |z|^2 \BPhi_{\delta}
%    \end{bmatrix}.
% \end{equation}
% By \cite[Lemma 13]{AblaNoEigen}, we know that the spectral radius $\rho(\BC_{\kappa, \delta}) \leq \sqrt{\rho(\BC_{\kappa}) \rho(\BC_{\delta})}$ and the max-row norm $ \| (\BI_{2n} - \BC_{\kappa, \delta})^{-1}\|_{\infty} \leq \sqrt{\| (\BI_{2n} - \BC_{\kappa})^{-1}\|_{\infty}\| (\BI_{2n} - \BC_{\delta})^{-1}\|_{\infty}} $. Therefore, it suffices to study the properties of $\BC_\kappa$ and $\BC_\delta$. 
In fact, by the identity $\Im a = (a - a^*) / 2 \jmath$ and using the same method to get \eqref{Eq_Maitrx_kappa_delta}, we can obtain the following system of equations
\begin{align*}
&\begin{bmatrix}
      \Im \Bkappa \\
      \Im z \widetilde{\Bkappa} 
  \end{bmatrix} = \begin{bmatrix}
    \Im \Bd \\
   \Im z\widetilde{\Bd}
\end{bmatrix}  + v \begin{bmatrix}
       \Bv_{\kappa}\\
      |z|^2 \Bu_{\kappa}
  \end{bmatrix} + \begin{bmatrix}
      |z|^2\BPsi_{\kappa} &  \BGamma_{\kappa} \\
      |z|^2\widetilde{\BGamma}_{\kappa} & |z|^2 \BPhi_{\kappa}
  \end{bmatrix}\begin{bmatrix}
      \Im \Bkappa \\
      \Im z \widetilde{\Bkappa} 
  \end{bmatrix}, \numberthis \label{Eq_Im_kappa_Im_wkappa}\\
  &\begin{bmatrix}
      \Im \Bdelta \\
      \Im z \widetilde{\Bdelta} 
  \end{bmatrix} = v \begin{bmatrix}
       \Bv_{\delta}\\
      |z|^2 \Bu_\delta
  \end{bmatrix} + \begin{bmatrix}
      |z|^2\BPsi_{\delta}&  \BGamma_{\delta} \\
      |z|^2\widetilde{\BGamma}_{\delta} & |z|^2 \BPhi_{\delta}
  \end{bmatrix}\begin{bmatrix}
      \Im \Bdelta \\
      \Im z \widetilde{\Bdelta} 
  \end{bmatrix}, \numberthis \label{Eq_Im_delta_Im_wdelta}
\end{align*}
with $[\Bv_\delta]_j = \frac{1}{n} \Tr \BOmega_j \BTheta^H \BTheta$, $[\Bu_\delta]_j = [\widetilde{\BTheta}^H \BA^H\BF^H\BF\BA\widetilde{\BTheta}]_{j, j}$, $[\Bv_\kappa]_j = \frac{1}{n} \Tr \BOmega_j \BR \BR^H$, and $[\Bu_\kappa]_j = [\widetilde{\BR} \BA^H \BG\BG^H \BA\widetilde{\BR}^H]_{j, j}$. By Lemma \ref{Lemm_Spr_pMat}, we can obtain $\rho(|z|^2 \BPhi_{\delta}) < 1$ for $z \in \mathbb{C}^+$ and 
\begin{align*}
    \norm{(\BI_{n} - |z|^2\BPhi_{\delta})^{-1}}_{\infty} &\leq \frac{\max_j(\Im (z \widetilde{\delta}_j))}{\min_j(v|z|^2 [\Bu_{\delta}]_j + \sum_{l} |z|^2 [\widetilde{\BGamma}_{\delta}]_{j, l} \Im(\delta_l))} \leq \frac{\max_j(\Im (z \widetilde{\delta}_j))}{\min_j (|z|^2 |\widetilde{\delta}_{j}|^2 \Im(\delta_j) )} \\
    & \overset{(a)}{\leq} \frac{\frac{K_0}{v}}{v^2 \cdot \frac{v^3K_1^3}{(K_2 + x^2+ v^2)^3}} = \frac{K_0(K_2 + x^2+ v^2)^3}{K_1^3 v^6}, \numberthis \label{Eq_I_minus_z_Phi_delta}
\end{align*}
where step $(a)$ follows by $|z| \geq v$, $|\widetilde{\delta}_j| \geq |\Im \widetilde{\delta}_j|$, and Lemma \ref{Lemm_Bound_on_qterms}. Here, the constants $K_j$, $j =0, 1, 2$ are the same as those in Lemma \ref{Lemm_Bound_on_qterms}. Without loss of generality, we assume $K_2 > 1$ in the following. Solving $\Im z\widetilde{\Bdelta}$ in \eqref{Eq_Im_delta_Im_wdelta}, we can get
\begin{align*}
\Im \Bdelta &= \BJ_{\delta} \Im\Bdelta  + v \Bv_{\delta} + v|z|^2\BGamma_{\delta}\left(\BI_n - |z|^2 \BPhi_{\delta}\right)^{-1} \Bu_{\delta}, \numberthis \label{Eq_J_delta}
\end{align*}
where $\BJ_\delta = |z|^2 \BPsi_{\delta} + |z|^2 \BGamma_{\delta} \left(\BI_n - |z|^2 \BPhi_{\delta}\right)^{-1}\widetilde{\BGamma}_\delta $. Using Lemma \ref{Lemm_Spr_pMat}, we have 
\begin{equation}
    \norm{(\BI_n - \BJ_\delta)^{-1}}_{\infty} \leq \frac{\max_j(\Im\delta_j)}{\min_j(v[\Bv_{\delta}]_j)} \leq \frac{K_0 (K_2 + x^2 + v^2)^2}{v^4 {K}_1^2}. \label{Eq_I_minus_J_delta_inv}
\end{equation}
\par
Next, we derive a similar norm bound for the matrices in \eqref{Eq_Im_kappa_Im_wkappa}. It can be observed that the RHS of \eqref{Eq_Im_kappa_Im_wkappa} involves $\Im \Bd$ and $\Im z \Bd$, whose elements are not always positive. To apply Theorem \ref{Lemm_Spr_pMat}, we need to determine the range of $z$. To that end, we define the following sets
\begin{equation}
 \mathfrak{A} = \left\{z: \frac{K(K_2 + x^2 + v^2)^3}{n^{\alpha}K_1^3v^{12}} \leq \frac{1}{2} \right\}, ~~ \mathfrak{E} = \{z: v \leq K_{\mathrm{im}}\}. \label{Eq_Def_of_range_z_A_E}
\end{equation}
Here, the constant $K$ is from the bounds $\|\Bd \|_{\infty} \leq K v^{-9}n^{-\alpha}$ and $\| z \widetilde{\Bd} \|_{\infty} \leq K v^{-7}n^{-\alpha}$. Then, for $z \in \mathfrak{A} \cap \mathfrak{E}$, we have $\rho(|z|^2 \BPhi_{\kappa}) < 1$ and 
\begin{align*}
  \norm{(\BI_n - |z|^2 \BPhi_{\kappa})^{-1}}_{\infty} &\leq \frac{\max_j (\Im(z\widetilde{\kappa_j}))}{\min_j(|z|^2|\widetilde{\kappa}_j|^2\Im(\kappa_j)) - \norm{z\Bd}_{\infty}}  \leq \frac{\frac{K_0}{v}}{\frac{v^5K_1^3}{(K_2 + x^2 + v^2)^3} - \frac{K}{n^\alpha v^7}}  \leq \frac{2K_0(K_2 + x^2 + v^2)^3}{K^3_1v^6}.  \numberthis \label{Eq_norm_I_minus_Phi}
\end{align*}
Solving \eqref{Eq_Im_kappa_Im_wkappa} with respective to $\Im (z\widetilde{\Bkappa})$, we can obtain
\begin{equation}
    (\BI_n - \BJ_{\kappa})\Im(\Bkappa) = \Im(\Bd) + v\Bv_{\kappa} + \BGamma_\kappa (\BI_n - |z|^2\BPhi_\kappa)^{-1} \Im(z\widetilde{\Bd})
     + |z|^2\BGamma_{\kappa} (\BI_n - |z|^2 \BPhi_\kappa)^{-1} \Bu_{\kappa},
\end{equation}
with $\BJ_\kappa = |z|^2 \BPsi_{\kappa} + |z|^2 \BGamma_{\kappa} \left(\BI_n - |z|^2 \BPhi_{\kappa}\right)^{-1}\widetilde{\BGamma}_\kappa$.
To evaluate $\| (\BI - \BJ_\kappa)^{-1} \|_{\infty}$,  we define another set $\mathfrak{B} \subset \mathbb{C}^+$ as
\begin{equation}
  \mathfrak{B} = \left\{ z: \frac{K(2K_0^3 + K_1^3)(K_2 + x^2 + v^2)^5}{K_1^5 v^{18}n^{\alpha}} \leq \frac{1}{2}\right\}. \label{Eq_Def_of_range_z_B}
\end{equation}
Hence, for $z \in \mathfrak{A} \cap \mathfrak{B} \cap \mathfrak{E}$, we have 
\begin{align*}
  \norm{(\BI_n - \BJ_{\kappa})^{-1}}_{\infty} & \leq \frac{\max_j(\Im(\kappa_j))}{\min_j(v [\Bv_{\kappa}]_j)- \norm{\Im(\Bd)}_{\infty} - \|\BGamma_\kappa (\BI_n-  |z|^2\BPhi_\kappa)^{-1} \Im(z\widetilde{\Bd})\|_{\infty}} \\
  &\leq \frac{\frac{K_0}{v}}{\frac{v^3{K}_1^2}{(K_2 + x^2 + v^2)^2} - \frac{K}{n^{\alpha}v^9} - \frac{2K_0^3(K_2 + x^2 + v^2)^3}{K_1^3v^8} \frac{K}{n^{\alpha }v^7}} \\
  & \leq \frac{\frac{K_0}{v}}{\frac{v^3{K}_1^2}{(K_2 + x^2 + v^2)^2} - \frac{K(2K_0^3 + K_1^3)(K_2 + x^2 + v^2)^3}{K_1^3v^{15}n^{\alpha}} } \leq \frac{2K_0 (K_2 + x^2 + v^2)^2}{v^4 {K}_1^2}, \numberthis
\end{align*}
given $\norm{\BGamma_{\kappa}}_{\infty} \norm{(\BI_n - |z|^2 \BPhi_{\kappa})^{-1}}_{\infty} \leq {2K_0^3(K_2 + x^2 + v^2)^3}/{(K_1^3v^8)}$ by \eqref{Eq_norm_I_minus_Phi} and Lemma \ref{Lemm_Bound_on_qterms}. Thus, by \eqref{Eq_Maitrx_kappa_delta}, the error can be controlled as
\begin{align*}
    \norm{\Bkappa - \Bdelta}_{\infty} &\leq  \norm{(\BI_n - \BJ_{\kappa, \delta})^{-1}}_{\infty}\left( \norm{\Bd}_{\infty} + \norm{\BGamma_{\kappa, \delta}}_{\infty} \norm{(\BI_n - z^2\BPhi_{\kappa, \delta})^{-1}}_{\infty} \norm{z\widetilde{\Bd}}_{\infty} \right) \\
    &\overset{(a)}{\leq} \frac{\sqrt{2}K_0 (K_2 + x^2 + v^2)^2}{v^4 {K}_1^2}\left(\frac{K}{v^9n^{\alpha} } + \frac{\sqrt{2}K_0^3(K_2 + x^2 + v^2)^3}{K_1^3v^8} \cdot \frac{K}{v^7n^{\alpha}}\right),\\
    & \leq  \frac{K'_0 (K_2 + x^2 + v^2)^5}{n^{\alpha}v^{19}}, ~~~ z \in \mathfrak{A} \cap \mathfrak{B} \cap \mathfrak{E}, \numberthis
\end{align*}
where $K'_0 > 0$ is a constant independent of $n$ and step $(a)$ follows from
$ \|(\BI_{n} - \BJ_{\kappa, \delta})^{-1} \|_{\infty} \leq \sqrt{\|(\BI_{n} - \BJ_{\kappa})^{-1} \|_{\infty}\|(\BI_{n} - \BJ_{ \delta} )^{-1}\|_{\infty}}$ and $\|(\BI_n - z^2 \BPhi_{\kappa, \delta})^{-1}\|_{\infty} \leq \sqrt{\|(\BI_n - |z|^2 \BPhi_{\kappa})^{-1}\|_{\infty}\|(\BI_n - |z|^2 \BPhi_{ \delta})^{-1}\|_{\infty} } $, which can be obtained by
Lemma \ref{Lemm_matrix_Spr_maxrow_Norm}. When $z \in (\mathfrak{A} \cap \mathfrak{B})^{c} \cap \mathfrak{E}$ and by definition of the three sets defined in \eqref{Eq_Def_of_range_z_A_E} and \eqref{Eq_Def_of_range_z_B}, we have
\begin{equation}
\label{Eq_Bound_on_Diff_kappa_delta}
  \norm{\Bkappa - \Bdelta}_{\infty} \leq \frac{4K_0}{v} \cdot \frac{1}{2} \leq \frac{K''_0}{n^{\alpha}}\max\left\{ \frac{(K_2 + x^2 + v^2)^5}{v^{19}}, \frac{(K_2 + x^2 + v^2)^3}{v^{13}} \right\} \overset{(a)}{\leq} \frac{K_0''(K_2 + x^2 + v^2)^5}{v^{19}n^{\alpha}},
\end{equation}
where constant $K_0'' > 0$ and step $(a)$ follows from $K_2 > 1$ and $v \leq K_{\mathrm{im}}$. By \eqref{Eq_Maitrx_kappa_delta} and \eqref{Eq_Bound_on_Diff_kappa_delta}, we get 
\begin{align*}
    \norm{\widetilde{\Bkappa} - \widetilde{\Bdelta}}_{\infty} &= \norm{(\BI_n - z^2\BPhi_{\kappa, \delta})^{-1} \left( \widetilde{\Bd} + z \widetilde{\BGamma}_{\kappa, \delta}(\Bkappa - \Bdelta)\right)}_{\infty} \leq 
    % \frac{K_0'''(K_2 + x^2+ v^2)^3}{K_1^3 v^6 n^{\alpha}}\left(\frac{1}{v^8} + \frac{(K_2 + x^2 + v^2)^{5.5}}{v^{21}} \right) \\
    % &\leq 
    \frac{K_0'''(K_2 + x^2+ v^2)^{8.5}}{ v^{27} n^{\alpha}}, ~~ z \in \mathfrak{E}, \numberthis \label{Eq_Bound_on_Diff_wkappa_wdelta}
\end{align*}
where we have used $1/|z| \leq 1 / v$ and $|z| \leq (K_2 + x^2 + v^2)^{1/2}$. 
% Thus, combining the bounds in \eqref{Eq_Bound_on_Diff_kappa_delta} and \eqref{Eq_Bound_on_Diff_wkappa_wdelta}, we have
% \begin{equation}
%   \max\left\{ \norm{\widetilde{\Bkappa} - \widetilde{\Bdelta}}_{\infty}, \norm{\Bkappa - \Bdelta}_{\infty} \right\} \leq \frac{\widetilde{K}_0(K_2 + x^2+ v^2)^{10}}{v^{27} n^{\alpha}},
% \end{equation}
% with $\widetilde{K}_0 > 0$ an absolute constant.  
Now we can give the bound on $|\Bu^H(\BR - \BTheta)\Bu|$ for any deterministic vector $\Bu$. In fact, we can bound the spectral norm of the difference $\BR - \BTheta$.
By the resolvent identity $\BR - \BTheta = \BR (\BTheta^{-1} - \BR^{-1})\BTheta$ and $\norm{\BF}$, $\norm{\BG} \leq 1 / v$, the following holds
\begin{align*}
  &\norm{\BR - \BTheta} \leq  \sum_{j=1}^n\frac{ \norm{\BR\BOmega_j\BTheta}}{n} \abs{z\widetilde{\kappa}_j - z\widetilde{\delta}_j} + |z|^2\norm{\BR\BA\BG(\BD_{\kappa} - \BD_{\delta})\BF\BA^H\BTheta} \\
  &\leq \widetilde{K}_0 \left( v^{-2}|z|\norm{\widetilde{\Bkappa} - \widetilde{\Bdelta}}_{\infty} +  v^{-4} |z|^2\norm{\Bkappa - \Bdelta}_{\infty}\right) 
  % \leq \widetilde{K}_0 \left( \frac{(K_2 + x^2+ v^2)^{9}}{ v^{29} n^{\alpha}} + \frac{(K_2 + x^2 + v^2)^6}{v^{23}n^{\alpha}} \right) 
  \leq  \frac{\widetilde{K}_0(K_2 + x^2+ v^2)^{9}}{ v^{29} n^{\alpha}}, ~~ v \leq K_{\mathrm{im}}, \numberthis \label{Eq_Boud_norm_diff_R_Theta}
\end{align*}
where $\BD_{\delta} = \diag(\delta_j; 1 \leq j \leq n)$ and $\BD_{\kappa} = \diag(\kappa_j; 1 \leq j \leq n)$. Then, \eqref{Eq_Part_III} follows by $|\Bu^H(\BR - \BTheta) \Bu| \leq \norm{\Bu}^2 \norm{\BR - \BTheta}$. Gathering the results in \eqref{Eq_Part_I}, \eqref{Eq_Part_II}, and \eqref{Eq_Part_III}, Lemma \ref{Lemm_bilinear} is proved. \qed
\par
We remark that in Lemma \ref{Lemm_bilinear}, we require $z$ to be close to the real axis such that $ \Im z \leq K_{\mathrm{im}} $. However, the convergence in Lemma \ref{Lemm_bilinear} still holds with a slightly different polynomial $\mathcal{U}(z)$ when $\Im z > K_{\mathrm{im}}$. Furthermore, the behavior of the resolvent near the real axis is more essential for investigating the local behavior of the ESD \cite{erdHos2017dynamical}. We will also use Lemma \ref{Lemm_bilinear} to prove the no-eigenvalue property. 
\par
Given $z \in \mathbb{C}^+$ and choosing $t$ to be large enough, the RHS of \eqref{Eq_bilinear_U_t} is summable. By the Borel-Cantelli lemma, we have $\Bu^H(\BQ - \BTheta) \Bu$  converges to 0 almost surely. The convergence proof for the trace form is similar, and we omit it for simplicity. Therefore, we have completed the proof of Theorem \ref{Thm_conver_ESD}.

% Define
% $
%   \BXi = |z|^2 \BPsi + |z|^2 \BGamma \left[\BI_n - |z|^2 \BPhi\right]^{-1}\widetilde{\BGamma} 
% $
% By solving \eqref{Eq_Im_delta_wdelta}, the following is obtained 
% \begin{equation}
%     (\BI_n - \BXi)\Im \Bdelta =  v \Bv + v|z|^2\BGamma \left[\BI_n - |z|^2 \BPhi\right]^{-1} \Bu.
% \end{equation} 
\section{Useful Results for the Proof of Lemma \ref{Lemm_bilinear}}
\label{App_Useful_Estimations}
In this section, we first provide some useful transformations and equations regarding the resolvent, and then give the estimation for certain terms that will be used to prove Lemma \ref{Lemm_bilinear}.
\subsection{Useful Identities}
The following identities can be proved using the Woodbury matrix identity
\begin{align*}
    & \widetilde{q}_j(z) = \frac{1}{-z(1 + \Bxi_j^H \BQ(z) \Bxi_j)} = [\widetilde{\BQ}(z)]_{j, j}, ~~ \BQ(z) = \BQ_j(z) +z \widetilde{q}_j(z)\BQ(z)_j \Bxi_j\Bxi_j^H \BQ_j(z), \\
    &\BQ_j(z) = \BQ(z) + \frac{\BQ(z)\Bxi_j\Bxi_j^H\BQ(z)}{1 - \Bxi_j^H\BQ(z)\Bxi_j}, ~~ \BQ(z)\Bxi_j = -z\widetilde{q}_j(z)\BQ_j(z)\Bxi_j, ~~~ 1 + \Bxi_j^H\BQ_j(z)\Bxi_j = \frac{1}{1 - \Bxi_j^H\BQ(z)\Bxi_j}. \numberthis \label{Eq_Wood_resolvent_Id}
\end{align*}
\subsection{Useful Bounds}
\begin{lemma}
\label{Lemm_Bound_on_qterms}
    For $z = x + \jmath v \in \mathbb{C}^+$, there exists constant $K_0$ (independent of $n$) such that the following holds
    \begin{align}
        & \max_{j, l \leq n} \left\{ |\kappa_j|, |\widetilde{\kappa}_j|, |\delta_j|, |\widetilde{\delta}_j|, |[\widetilde{\BTheta}]_{j, l}|, |[\widetilde{\BR}]_{j, l}| \right\} \leq \frac{K_0}{v}, \label{Eq_Upper_bounds_terms_1}\\
        &\max_{j, l \leq n} \left\{ \frac{1}{n} \Tr \BOmega_j \BTheta \BOmega_l \BTheta^H, \frac{1}{n} \Tr \BOmega_j \BR \BOmega_l \BR^H,\frac{1}{n} \abs{\Tr \BOmega_j \BR \BOmega_l \BTheta^H} \right\} \leq \frac{K_0^2}{v^2}, \label{Eq_Upper_bounds_terms_1_5}\\
        &\max_{j \leq n} \left\{ |\Im(z \kappa_j)|, |\Im(z \widetilde{\kappa}_j)|, |\Im(z \delta_j)|, |\Im(z \widetilde{\delta}_j)| \right\}  \leq \frac{K_0}{v},  \label{Eq_Upper_bounds_terms_2}\\
        &\max_{j \leq n} |z\widetilde{q}_j| \leq \max\left\{K_0 \norm{\BSigma}^2 / v, 2\right\}, ~~ \max_{j \leq n} \left\{ |\eta_j|,\frac{1}{\abs{1 + \kappa_j}}, \frac{1}{\abs{1 + \delta_j}}, \frac{1}{\abs{1 + \E \alpha_{j, j}}} \right\} \leq \max\left\{K_0 / v, 2\right\}. \label{Eq_Upper_bounds_terms_3}
        % &\frac{1}{\abs{1 + \kappa_j}}, \frac{1}{\abs{1 + \delta_j}}, \frac{1}{\abs{1 + \E \alpha_{j, j}}} \leq \max\left\{K_0 / v, 2\right\}. \label{Eq_Upper_bounds_tersm}
    \end{align}
    Moreover, there exist constants $K_1$ and $K_2$ (independent of $n$) such that
    \begin{align}
  &\min_{j \leq n} \left\{ |\Im \kappa_j|, |\Im \widetilde{\kappa}_j|, |\Im\frac{1}{n}\Tr\BOmega_j \BR|, |\Im[\widetilde{\BR}]_{j, j}|, |\Im \delta_j|, |\Im \widetilde{\delta}_j| \right\} \geq \frac{v K_1}{(K_2 + x^2 + v^2)}, \label{Eq_Lower_bounds_terms_1} \\
  &\min_{j \leq n} \left\{\frac{1}{n} \Tr \BOmega_j \BTheta \BTheta^H, \frac{1}{n} \Tr \BOmega_j \BR\BR^H \right\} \geq \frac{v^2 {{K}_1^2}}{(K_2 + x^2 + v^2)^2}. \label{Eq_Lower_bounds_terms_2} 
\end{align}
\end{lemma}
\textit{Proof:}  We will prove the upper bounds and the lower bounds separately.
\subsubsection{Proof of the Upper Bounds} 
The proof of the first inequality in \eqref{Eq_Upper_bounds_terms_1} follows from $\|\BQ\|, \|\widetilde{\BQ}\|, \|\BR\|, \|\BTheta\|,  \|\widetilde{\BR}\|, \|\widetilde{\BTheta}\| \leq 1 / v$, and the property that if $m$ is a Stieltjes transform, then $|m| \leq  \int \frac{\dd \mu}{|\lambda - z|} \leq \mu(\R) / v$. The second inequality \eqref{Eq_Upper_bounds_terms_1_5} follows from $\norm{\BTheta}, \norm{\BR} \leq 1 / v$.
\par
For simplicity, we will only prove the inequality for $\widetilde{\delta}_j$ and $\kappa_j$ in \eqref{Eq_Upper_bounds_terms_2}. The proof mainly relies on the following identity \cite[Lemma C1]{hachem2007deterministic}.
\begin{equation}
    \int_{\R^+} \lambda \mu(\dd \lambda) = \lim_{y \uparrow \infty} \Re\left\{ - \jmath y \left[\jmath y m(\jmath y) + \mu(\R^+) \right] \right\}. \label{Eq_int_lambda_dmu}
\end{equation}
By the properties of the matrix-valued Stieltjes transform \cite[Proposition 2.2]{hachem2007deterministic}, we have
\begin{equation}
    \BTheta(z) = \int_{\R^+} \frac{\Bmu(\dd \lambda)}{\lambda - z}, \label{Eq_representation_Theta}
\end{equation}
where $\Bmu = ([\Bmu]_{i, j})_{1 \leq i, j \leq p}$ is a matrix-valued (complex) measure such that $\Bmu(A)$ is Hermitian non-negative for any Borel set $A$ and  $\Bmu(\R^+) = \BI_p$. Thus, we have
\begin{equation}
    \lim_{y \uparrow \infty } -\jmath y \BTheta(\jmath y) = \BI_p, ~~ \lim_{y \uparrow \infty} - \jmath y \delta_j(\jmath y) = \frac{1}{n} \Tr \BOmega_j.
\end{equation}
By the resolvent identity, we can obtain 
\begin{equation}
    \widetilde{\delta}_j(z) = -\frac{1}{z(1 + \delta_j(z))} + \frac{\Ba_j^H \BTheta(z) \Ba_j}{z(1 + \delta_j(z))^2}. \label{Eq_Res_identity_wdeltaj}
\end{equation}
Hence, it holds that $\lim_{y \uparrow \infty} -\jmath y \widetilde{\delta}_j(\jmath y) = 1$, which implies $\widetilde{\mu}_j (\R^+) = 1$. Thus, by \eqref{Eq_int_lambda_dmu}, we have  
\begin{align*}
    \int_{\R^+} \lambda \widetilde{\mu}_j(\dd \lambda) &= \lim_{y \uparrow \infty} \Re\left\{ - \jmath y \left[\jmath y \widetilde{\delta}_j(\jmath y) + \widetilde{\mu}_j(\R^+) \right] \right\} = \lim_{y \uparrow \infty} \Re\left\{ \left[\frac{- \jmath y {\delta}_j(\jmath y)}{1 + \delta_j(\jmath y) } + \frac{\Ba_j^H (-\jmath y \BTheta(\jmath y)) \Ba_j}{[1 + \delta_j(\jmath y)]^2} \right] \right\} \\
    &= \frac{1}{n} \Tr \BOmega_j + \norm{\Ba_j}^2 < \infty, \numberthis \label{Eq_Bound_of_wmu_j}
\end{align*}
which implies 
\begin{equation}
    \Im(z \widetilde{\delta}_j) = \int_{\R^+} \frac{\lambda v\widetilde{\mu}_j(\dd \lambda)}{(x - \lambda)^2 + v^2} \leq \frac{\frac{1}{n} \Tr \BOmega_j + \norm{\Ba_j}^2}{v}.
\end{equation}
\par
Denote the eigen-decomposition of $\BS$ as $\BS = \sum_{j=1}^p \lambda_j \Bu_j^H\Bu_j$ and the underline measure of ${\kappa}_j$ as $\mu_{\kappa, j}$. Then, we can obtain 
\begin{align*}
    \kappa_j = \frac{1}{n} \Tr \BOmega_j \E \BQ = \sum_{j=1}^p \E 
 \frac{\Bu_j^H \BOmega_j \Bu_j}{n(\lambda_j - z)} = \E \int_{\R^+} \frac{F^{\kappa, j}(\dd \lambda)}{\lambda - z}, \numberthis
\end{align*}
where $F^{\kappa, j} = \frac{1}{n}\sum_j \Bu_j^H \BOmega_j \Bu_j \Ind\{\lambda \leq \lambda_j\}$. Hence, we have $\mu_{\kappa, j} = \E F^{\kappa, j}$ and 
\begin{align*}
    &\int_{\R^+}\lambda \mu_{\kappa, j}(\dd \lambda) = \frac{1}{n}\E \left( \sum_{j=1}^n \lambda_j\Bu_j^H \BOmega_j \Bu_j\right) = \frac{1}{n} \Tr \BOmega_j \E \BSigma \BSigma^H \\
    &= \frac{1}{n} \Tr \BOmega_j \left(\sum_{l=1}^n \frac{1}{n} \BOmega_l + \Ba_l\Ba_l^H \right) = \frac{1}{n^2} \sum_{l=1}^n \Tr (\BOmega_j \BOmega_l) + \frac{1}{n}\Tr \BOmega_j \BA\BA^H, \numberthis
\end{align*}
which also implies $\Im(z \kappa_j) \leq K_0 / v$.
\par
The proof for \eqref{Eq_Upper_bounds_terms_3} follows from the property in \cite[Lemma 2.3]{SILVERSTEIN1995strong} which states that for positive $x$ and Stieltjes transform $m(z)$ with underline measure $\mu$, it holds that $|m(z)| \leq \mu(\mathbb{R}) / v$ and $|1 + m(z)x|^{-1} \leq \max(4\mu(\mathbb{R})x /v, 2)$.
\subsubsection{Proof of the Lower Bounds} For \eqref{Eq_Lower_bounds_terms_1}, we only show the proof for the inequality of $\widetilde{\delta}_j$ as proofs for the other terms are quite similar.  Denote $\widetilde{\mu}_j = \widetilde{\mu}_j^n$. By the bound \eqref{Eq_Bound_of_wmu_j}, we know that  $\sup_{j, n} \int \lambda \widetilde{\mu}^n_j(\dd \lambda) < \infty$. This implies that for each $j$, the positive measure sequence  $\{\widetilde{\mu}_j^n\}_{n \geq j}$ is tight and thus there exists a constant $K_\mu$ (independent of $n$) such that $\widetilde{\mu}_j([0, K_\mu]) = \widetilde{\mu}_j^n([0, K_\mu]) \geq 1 / 2$ for each $n$. Then, we have
\begin{align*}
    \Im (\widetilde{\delta}_j(z)) = \int_{0}^{\infty} \frac{v \widetilde{\mu}_j(\dd \lambda)}{(x - \lambda)^2 + v^2} \geq \int_{0}^{K_\mu}\frac{v \widetilde{\mu}_j(\dd \lambda)}{(x - \lambda)^2 + v^2} \geq \frac{\widetilde{\mu}_j([0, K_\mu])v}{2(x^2 + K_\mu^2) + v^2} \geq \frac{v}{4(x^2 + K_\mu^2 + v^2)}, \numberthis
\end{align*}
which gives the lower bound for $|\Im \widetilde{\delta}_j|$. The lower bounds for $\frac{1}{n} \Tr \BOmega_j \BTheta\BTheta^H$ in \eqref{Eq_Lower_bounds_terms_2} follows  by the lower bound for $|\Im \delta_j|$ in \eqref{Eq_Lower_bounds_terms_1} and the Cauchy-Schwarz inequality
\begin{equation}
    |\Im(\delta_j)| \leq |\delta_j| = \abs{ \frac{1}{n}\Tr \BOmega_j^{1/2} \BOmega_j^{1/2}\BTheta  } \leq \sqrt{\left( \frac{1}{n}\Tr \BOmega_j \right) \left(  \frac{1}{n}\Tr\BOmega_j\BTheta \BTheta^H \right)}.
\end{equation}
The evaluation for $\frac{1}{n} \Tr \BOmega_j \BR\BR^H$ is similar. Therefore, we have completed the proof of Lemma \ref{Lemm_Bound_on_qterms}. \qed
\subsection{Useful Estimations}
\begin{lemma}
\label{Lemm_bound_Delta}
    For $t \geq 2$, we have 
    \begin{equation}
        \E\left[\abs{\Delta_j}^t | \widetilde{\mathscr{F}}_{j} \right] \lesssim_t \left\{\begin{aligned}
        &\frac{1}{n^{t/2} v^{t}} & \text{if } 2 \leq t \leq 2 + \varepsilon/2, \\
        &\frac{1}{v^t} \left(\frac{1}{n^{2\alpha t + 1}} + \frac{1}{n^{t/2}} \right) & \text{if } 2 + \varepsilon/2 < t \leq 4 + \varepsilon, \\
        &\frac{1}{v^t}\left(\frac{1}{n^{2\alpha t + 1}} + \frac{1}{n^{\alpha t + 2}} \right) & \text{if } t > 4 + \varepsilon,
    \end{aligned}\right. 
    \end{equation}
    where $\widetilde{\mathscr{F}}_{j}$ is a  $\sigma$-algebra independent of $\sigma(X_{ij}, 1 \leq i \leq d_j)$.
\end{lemma}
\textit{Proof:} By Lemma \ref{Lemm_trace} , $\norm{\BQ_j} \leq 1/v$, and inequality $(\sum_{i=1}^L a_i)^t \leq L^{t-1} \sum_{i=1}^L|a_i|^t$, we have 
\begin{align*}
\E\left[|\Delta_j|^t|\widetilde{\mathscr{F}}_{j}\right]  \leq 3^{t-1} \E\left[\abs{\By_j^H \BQ_j \By_j - \frac{1}{n} \Tr \BOmega_j \BQ_j}^t + \abs{\By_j^H \BQ_j \Ba_j}^t + \abs{\Ba_j^H \BQ_j \By_j}^t \Big|\widetilde{\mathscr{F}}_{j}\right]\lesssim_t\frac{R_t}{v^tn^t} + \frac{T_{t/2}}{n^{t/2}v^t}, \numberthis
\end{align*}
which yields the desired result. \qed
\begin{lemma}
\label{Lemm_Usedul_est}
For deterministic vector $\Bu \in \mathbb{C}^H$ with $\norm{\Bu} = 1$ and $z = x + \jmath v \in \mathbb{C}^+$ with $v$ small ($v \leq K_{\mathrm{im}}$), we have 
\begin{align}
&\E \left[ \sum_{j=1}^n \E_{j-1} \left( \Bu^H \BQ\Ba_j\Ba_j^H\BQ^H\Bu \right) \right]^t \lesssim_t \frac{1}{v^{2t}}, \label{Eq_useful_est_1}\\
&\E \left[ \sum_{j=1}^n \E_{j-1} \left( \Bu^H \BQ\Bxi_j\Bxi_j^H\BQ^H\Bu \right) \right]^t \lesssim_t  \frac{\log^{2t}(n)}{v^{2t}}, \label{Eq_useful_est_2}\\
&\sum_{j=1}^n \E (\Bu^H \BQ_j \Ba_j\Ba_j^H\BQ_j^H\Bu)^2 \lesssim \frac{\log^4(n)}{v^8},\label{Eq_useful_est_3}  \\
&\E \left[ \sum_{j=1}^n \E_{j-1}(\Bu^H\BQ_j\Ba_j\Ba_j^H\BQ^H_j\Bu) \right]^t \lesssim_t \frac{\log^{2t}(n)}{v^{4t}}. \label{Eq_useful_est_4}
\end{align}
\end{lemma}
\textit{Proof:} The proof of \eqref{Eq_useful_est_1} is similar to that of \cite[Lemma 3.5]{hachem2013bilinear} and omitted. Next we prove \eqref{Eq_useful_est_2}.
When $t = 1$, we have
\begin{equation}
     \E \sum_{j=1}^n \E_{j-1} \Bu^H\BQ\Bxi_{j}\Bxi_j^H\BQ^H\Bu  = \E \Bu^H\BQ\BSigma\BSigma^H\BQ^H\Bu \overset{(a)}{\leq} \frac{\E \norm{\BSigma}^2}{v^2} \overset{(b)}{\lesssim} \frac{\log^2(n)}{v^2},
\end{equation}
where step $(a)$ follows from $\norm{\BQ} \leq 1 / \Im z$ and step $(b)$ is due to Lemma \ref{Lemm_P_bound_Sn}. Denote the set $A = \{ \norm{\BSigma} \leq \log(n) \}$. If $t \geq 2$, we have 
\begin{align*}
    & \E \left[ \sum_{j=1}^n \E_{j-1} \Bu^H\BQ\Bxi_j\Bxi_j^H\BQ^H\Bu \right]^t (\Ind_A + \Ind_{A^c}) = \sum_{j_{1}, \ldots j_t} \E \left[ \prod_{i=1}^t\E_{j_i-1} \Bu^H\BQ\Bxi_{j_i}\Bxi_{j_i}^H\BQ^H\Bu \right](\Ind_A + \Ind_{A^c}) \\
    &\overset{(a)}{\leq} t! \sum_{j_{1}\leq \ldots \leq j_t} \E \left[ \left(\prod_{i=1}^{t-1}\E_{j_i-1}  \Bu^H\BQ\Bxi_{j_i}\Bxi_{j_i}^H\BQ^H\Bu \right) \left( \Bu^H\BQ\Bxi_{j_t}\Bxi_{j_t}^H\BQ^H\Bu \right)\right]\Ind_A + \PP(A^c)n^{3t}v^{-2t} \\
    &\leq t! \sum_{j_{1}\leq \ldots \leq j_{t-1}} \E \left[ \left(\prod_{i=1}^{t-1}\E_{j_i-1} \Bu^H\BQ\Bxi_{j_i}\Bxi_{j_i}^H\BQ^H\Bu \right) \sum_{j_t = j_{t-1}}^n \left( \Bu^H\BQ\Bxi_{j_t}\Bxi_{j_t}^H\BQ^H\Bu \right)\right] + \PP(A^c)n^{3t}v^{-2t} \\
    &\overset{(b)}{\leq} t! \E \left[ \sum_{j=1}^n \Bu^H\BQ\Bxi_{j}\Bxi_{j}^H\BQ^H\Bu \right]^{t-1}\frac{\log^2(n)}{v^2} + o(K_ln^{-l}), ~~ \forall l. \numberthis
\end{align*}
where step $(a)$ follows from $\norm{\Bxi_j} \leq \sqrt{d_j} + \norm{\BA} \leq n$ (for large $n$) and $\Bu^H\BQ \Bxi_j \Bxi_j^H\BQ^H \Bu \leq v^{-2}n^2$, and step $(b)$ is due to the inequality
\begin{equation}
    \sum_{j_t = j_{t-1}}^n \left( \Bu^H\BQ\Bxi_{j_t}\Bxi_{j_t}^H\BQ^H\Bu \right) \leq \Bu^H\BQ\BSigma\BSigma^H\BQ^H\Bu \leq v^{-2} \norm{\BSigma}^2
\end{equation}
and Lemma \ref{Lemm_P_bound_Sn}. As a result, \eqref{Eq_useful_est_2} is proved.
\par
To prove \eqref{Eq_useful_est_3}, we write $\Bu^H\BQ_j\Ba_j\Ba_j^H\BQ_j^H\Bu$ as 
\begin{align*}
    \Bu^H\BQ_j\Ba_j\Ba_j^H\BQ_j^H\Bu &= \Bu^H(\BQ_j - \BQ)\Ba_j\Ba_j^H\BQ^H\Bu + \Bu^H\BQ\Ba_j\Ba_j^H(\BQ_j - \BQ)^H\Bu + \Bu^H(\BQ - \BQ_j)\Ba_j\Ba_j^H(\BQ - \BQ_j)^H\Bu \\
    &+ \Bu^H\BQ\Ba_j\Ba_j^H\BQ^H\Bu = X_{1j} + X_{2j} + X_{3j} + X_{4j}. \numberthis \label{Eq_Decomp_of_uQaaQu}
\end{align*}
By $\E |X_{1j}|^2 \leq \frac{1}{2}(\E X_{3j}^2 + \E X_{4j}^2 )$ and $\E |X_{2j}|^2 \leq \frac{1}{2}(\E X_{3j}^2 + \E X_{4j}^2 )$, we have 
\begin{equation}
    \sum_{j=1}^n \E (\Bu^H \BQ_j \Ba_j\Ba_j^H\BQ_j^H\Bu)^2 \leq 8 \sum_{j=1}^n \left(\E X_{3j}^2 + \E X_{4j}^2\right).
\end{equation}
By $\Bu^H\BQ\Ba_j\Ba_j^H\BQ^H\Bu \leq \sum_i \Bu^H\BQ\Ba_i\Ba_i^H\BQ^H\Bu = \Bu^H\BQ\BA\BA^H\BQ^H\Bu$, we can obtain $\sum_j \E X_{4j}^2 \leq \E(\Bu^H\BQ\BA\BA^H\BQ^H\Bu)^2 \leq \norm{\BA}^4 / v^4$. According to \eqref{Eq_Wood_resolvent_Id}, we have 
\begin{align*}
    \sum_{j=1}^n \E X_{3j}^2 &= \sum_{j=1}^n \E \left( \Bu^H\BQ\Bxi_j\Bxi_j^H\BQ\Ba_j\Ba_j^H \BQ^H \Bxi_j \Bxi_j^H \BQ^H\Bu\abs{\frac{1 + \Bxi_j^H\BQ_j\Bxi_j}{1 - \Bxi_j^H \BQ \Bxi_j}} \right)^2 \\
    &= \sum_{j=1}^n\E \abs{\frac{\Bxi_j^H\BQ\Ba_j\Ba_j^H \BQ^H \Bxi_j}{1 - \Bxi_j^H \BQ \Bxi_j} }^2  \abs{\Bu^H\BQ\Bxi_j \Bxi_j^H \BQ^H\Bu}^2\abs{1 + \Delta_j + \frac{1}{n} \Tr \BOmega_j\BQ_j + \Ba_j \BQ_j \Ba_j}^2. \numberthis
\end{align*}
Write the eigen-decomposition of $\BQ$ as $\BQ = \BU \diag{(\lambda_i - z; 1 \leq i \leq p) \BU^H}$, and denote $\overline{\Bxi}_j = \BU^H \Bxi_j = [\overline{\xi}^j_1, \ldots \overline{\xi}^j_p]^T$. Then, we have 
\begin{align*}
\label{Eq_Resol_Im_Bound}
\abs{\frac{\Bxi_j^H\BQ\Ba_j\Ba_j^H \BQ^H \Bxi_j}{1 - \Bxi_j^H \BQ \Bxi_j}} \leq \norm{\Ba_j}^2\abs{\frac{\Bxi_j^H\BQ \BQ^H \Bxi_j}{\Im (1 - \Bxi_j^H \BQ \Bxi_j)}} \leq \norm{\Ba_j}^2\abs{\frac{\sum_i \frac{|\overline{\xi}^j_i|^2}{|\lambda_i - z|^2}}{\Im\left(\sum_i \frac{|\overline{\xi}^j_i|^2}{\lambda_i - z}\right)}} \leq \frac{\norm{\Ba_j}^2}{v}. \numberthis
\end{align*}
Therefore, by $\abs{\frac{1}{n} \Tr \BOmega_j \BQ_j} \leq \frac{p}{n}\norm{\BQ_j} \norm{\BOmega_j}\lesssim v^{-1}$, $\abs{\Ba_j^H \BQ_j \Ba_j} \leq \norm{\Ba_j}^2 \norm{\BQ_j} \lesssim v^{-1}$ and Lemma \ref{Lemm_bound_Delta}, we can obtain  
\begin{align*}
    \sum_{j=1}^n \E X_{3j}^2 &\lesssim v^{-2} \sum_{j=1}^n\E  \abs{\Bu^H\BQ\Bxi_j \Bxi_j^H \BQ^H\Bu}^2(|\Delta_j|^2 + 1 + v^{-1} +  v^{-2}) \\
    &\leq v^{-2}(1 + v^{-1} +  v^{-2})\E(\Bu^H\BQ\BSigma\BSigma^H\BQ^H\Bu)^2 + \frac{1}{v^6}\sum_{j=1}^n \E \norm{\BSigma}^4 |\Delta_j|^2 \lesssim \frac{\log^4(n)}{v^8}. \numberthis
\end{align*}
Here, we note that $v$ is small enough by assumption, which implies $v^{-b} \leq v^{-a}$ for $0 < a \leq b$. 
\par
In the following, we prove the last inequality in \eqref{Eq_useful_est_4}. With the decomposition in \eqref{Eq_Decomp_of_uQaaQu}, we have
\begin{align*}
    \E \left[ \sum_{j=1}^n \E_{j-1}(\Bu^H\BQ_j\Ba_j\Ba_j^H\BQ^H_j\Bu) \right]^t \lesssim_t \sum_{i=1}^4 \E \abs{ \sum_{j=1}^n \E_{j-1}X_{ij} }^t. \numberthis
\end{align*}
By \eqref{Eq_useful_est_1}, we can obtain $\E | \sum_{j} \E_{j-1}X_{4j} |^t \lesssim_t v^{-2t}$. Similar to the derivation of \eqref{Eq_useful_est_3}, we have
\begin{align}
    \E_{j-1} X_{3j} \lesssim (v^{-1} + v^{-2}) \E_{j-1} \abs{\Bu^H\BQ\Bxi_j\Bxi_j^H\BQ^H\Bu} + v^{-1} \E_{j-1} |\Delta_j|\abs{\Bu^H\BQ\Bxi_j\Bxi_j^H\BQ^H\Bu}.
\end{align}
According to the Cauchy-Schwarz inequality $|\E[ab|\mathscr{F}]| \leq \E^{1/2}[|a|^2|\mathscr{F}]\E^{1/2}[|b|^2|\mathscr{F}]$ and \eqref{Eq_useful_est_2}, we have 
\begin{align*}
    &\E\abs{\sum_{j=1}^n  \E_{j-1} |\Delta_j|\abs{\Bu^H\BQ\Bxi_j\Bxi_j^H\BQ^H\Bu}}^t \leq \E\left( \sum_{j=1}^n  \E_{j-1} |\Delta_j|^2 \right)^{t/2} \E\left( \sum_{j=1}^n  \E_{j-1} |\Bu^H\BQ\Bxi_j\Bxi_j^H\BQ^H\Bu|^2 \right)^{t/2} \\
    &\lesssim_t v^{-t} \cdot v^{-t} \E\left( \sum_{j=1}^n  \E_{j-1} \norm{\BSigma}^2(\Bu^H\BQ\Bxi_j\Bxi_j^H\BQ^H\Bu) \right)^{t/2} \lesssim_t  \frac{\log^{2t}(n)}{v^{3t}}. \numberthis
\end{align*}
Hence, we get $\E | \sum_{j} \E_{j-1}X_{3j}|^t \lesssim_t \log^{2t}(n) v^{-4t}$. 
Using Cauchy-Schwarz again, we can obtain $|\E_{j-1}X_{ij}| \leq \sqrt{\E_{j-1}X_{3j}\E_{j-1}X_{4j}} $ for $i=1,2$, and
\begin{align}
    \E \abs{\sum_{j=1}^n \E_{j-1}X_{ij}}^t \leq  \E \left( \sum_{j=1}^n \E_{j-1}X_{3j} \right)^{t/2}\E \left(\sum_{j=1}^n \E_{j-1}X_{4j}\right)^{t/2}\lesssim_t \frac{\log^t(n)}{v^{3t}}, ~~ i = 1, 2.
\end{align}
Therefore, \eqref{Eq_useful_est_4} is proved, which completes the proof of Lemma \ref{Lemm_Usedul_est}. \qed
\section{Proof of Proposition \ref{Prop_same_support}}
\label{App_Proof_Prop_same_support}
Fix $x_0 \in \R_*$, and since $\mathsf{Supp}(F^n) \subset \mathbb{R}^+$, we can assume $x_0 > 0$ without loss of generality. We then define $\Bdelta(z) = [\delta_1(z), \ldots, \delta_n(z)]^T$ and $\widetilde{\Bdelta}(z) = [\delta_1(z), \ldots, \delta_n(z)]^T $, and denote three limits $\Bdelta_0 = \lim_{z \in \mathbb{C}^+ \to x_0} \Bdelta(z)$, $\widetilde{\Bdelta}_0 = \lim_{z \in \mathbb{C}^+ \to x_0} \widetilde{\Bdelta}(z)$, and $m_{ n, 0} = \lim_{z \in \mathbb{C}^+ \to x_0} m_n(z)$. 
According to \cite[Theorem B.10]{bai2010spectral}, it is enough to show the following holds
\begin{equation}
    \Im(m_{ n, 0}) = 0 \Rightarrow \Im[\Bdelta_0] = \Bzero_{n} \Rightarrow  \Im[\widetilde{\Bdelta}_0] = \Bzero_{n}. \label{Eq_Support_arrows}
\end{equation} 
To this end, we choose sequence $\{z_k\} \subset \mathbb{C}^+ \to x_0$. By the integration representation \eqref{Eq_representation_Theta}, we have 
\begin{equation}
    \Im(\BTheta(z_k)) = [\BTheta(z_k) - \BTheta^H(z_k)] / 2\jmath = \int_{\R^+} \frac{v\Bmu(\dd\lambda)}{|\lambda - z|^2}  > \mathbf{0}_{p \times p}. \label{Eq_Tm_Theta_geq_0}
\end{equation}
Then, by the inequality $ |\Tr \BC\BD| \leq \norm{\BC} \Tr \BD$ for Hermitian non-negative $\BD$, we can obtain 
\begin{equation}
    0 \leq \Im (\delta_j(z_k)) = \frac{1}{2\jmath n} \Tr \BOmega_j (\BTheta(z_k) - \BTheta^H(z_k)) \leq \frac{\norm{\BOmega_j}}{n } \Tr\Im (\BTheta(z_k)) \leq \frac{p \norm{\BOmega_j}\Im(m(z_k))}{n }.
\end{equation}
Letting $z_k \to x_0$ in the above inequality, we have $\Im \Bdelta_0 = \Bzero_n$. 
\par
To prove the second arrow in \eqref{Eq_Support_arrows}, we will use \eqref{Eq_Im_delta_Im_wdelta} and prove by contradiction. Assume that there exists an index $m$ such that $\Im[\widetilde{\Bdelta}_0]_m > 0$. By \eqref{Eq_Im_delta_Im_wdelta}, we can write the equation regarding $\Im\delta_m$ as
\begin{equation}
    \Im(\delta_m) = v_k[\Bv_{\delta}]_m + |z_k|^2 [\BPsi_{\delta} \Im(\Bdelta)]_m + [\BGamma_\delta \Im(z_k\widetilde{\Bdelta})]_m,
\end{equation}
which implies
$
    \Im(\delta_m(z_k)) \geq [\BGamma_{\delta}]_{m, m} \Im (z_k\widetilde{\delta}_m(z_k))
$. Letting $z_k \to x_0$, we have 
$
    [\BGamma_{\delta}]_{m, m} \to 0.
$
By the Cauchy-Schwarz inequality, the following holds 
\begin{equation}
    [\BGamma_{\delta}]_{m, m} = \frac{1}{n}\Tr \BOmega_m^{1/2} \BTheta(z_k) \BOmega_m \BTheta^H(z_k)\BOmega_m^{1/2} \geq \frac{1}{np}\abs{\Tr \BOmega_m \BTheta(z_k)} ^2 \geq \frac{n|\delta_m(z_k)|^2}{p}.
\end{equation}
Hence, we have $\delta_m(z_k) \to 0$, as $z_k \to x_0$. By resolvent identity \eqref{Eq_Res_identity_wdeltaj}, we can obtain
\begin{equation}
    \Ba_m^H \BTheta(z_k) \Ba_m = z_k(1 + \delta_m(z_k))^2 \widetilde{\delta}_m(z_k) + (1 + \delta_m(z_k)).
\end{equation}
By taking the imaginary part of the above identity, we have
\begin{equation}
    \lim_{k \to \infty} \Im(\Ba_m^H \BTheta(z_k) \Ba_m ) = \lim_{k \to \infty} \Ba_m^H \Im(\BTheta(z_k)) \Ba_m = x_0 \Im([\widetilde{\Bdelta}_0]_m) + 1 > 0.
\end{equation}
However, this cannot happen. In particular, because $\Im[\BTheta(z_k)] > \mathbf{0}_{p \times p}$ according to \eqref{Eq_Tm_Theta_geq_0}, we have  $\Im(m_n(z_k)) = \Im(\Tr\BTheta(z_k)) / p \to 0$, which yields $\Im[\BTheta(z_k)] \to \mathbf{0}_{p \times p}$. Therefore, we have completed the proof. \qed

\section{Proof of \eqref{Eq_sup_diff_resolvent}}
\label{App_Proof_Thm_noeigen}
 In this section, we use the same notations as defined in Appendix \ref{App_Converge_resolvent}. Since $[a, b]$ is a subset of a open interval outside the limiting support, we assume that there exists $\underline{\epsilon} > 0$ such that $[a - 2\underline{\epsilon}, b + 2\underline{\epsilon}]$ is outside the limiting support, and denote $a' = a - \underline{\epsilon}$ and $b' = b + \underline{\epsilon}$.
\subsection{A Rate on $F^{\BS, n}([a, b])$}
\label{Sec_A_rate_one_F}
 In this section, we provide a preliminary order for $F^{\BS, n}([a, b]) = (\text{number of eigenvalues of $\BS$ appears in $[a, b]$}) / p$. To this end, the following proposition is useful. 
\begin{proposition}
\label{Prop_bilinear_of_resolvent_uniform_x}
Assume $r \geq 1$ and $z = x + jv$ with $v = n^{-{\alpha}/(48 + 20 r)}$. Let $\BT \in \mathbb{C}^{p \times p}$ be a Hermitian non-negative matrix with bound spectral norm and $\Bu$ be a deterministic vector with bounded Euclidean norm. For any $l > 0$, the following holds.
\begin{align}
\max_{0 \leq j \leq n} \E_{j} \left\{ v^{-rl}p^{-l}\sup_{x \in \mathbb{R}}\abs{\Tr \BT [\BQ(z) -\BTheta(z)] }^l\right\} \xrightarrow{a.s.} 0, \\
\max_{0 \leq j \leq n} \E_{j} \left\{ v^{-{rl}}\sup_{x \in \mathbb{R}}{\abs{\Bu^H[\BQ(z) - \BTheta(z)]\Bu}^l}\right\} \xrightarrow{a.s.} 0. \label{Eq_Bilinear_martingale_conv}
\end{align}
\end{proposition}
\textit{Proof: } We first evaluate $\sup_{x \in S_n}{\left\{ v^{-{rl}}\sup_{x \in \mathbb{R}}{\abs{\Bu^H[\BQ(z) - \BTheta(z)]\Bu}^l}\right\}}$, where $S_n \subset \mathbb{R}$ contains at most $n$ elements. Denote  $m_{\Bu,\BS, n}(z) = \Bu^H\BQ(z)\Bu$ and $m_{\Bu, n}(z) =  \Bu^H\BTheta(z)\Bu$. We have the following decomposition
\begin{align*}
  \sup_{x\in S_n} v^{-r}\abs{m_{\Bu, \BS, n}(z) - m_{\Bu, n}(z)} &= \sup_{x\in S_n} v^{-r}\abs{m_{\Bu, \BS, n}(z) - m_{\Bu, n}(z)}\Ind\{|x| \leq {\varsigma}^rv^{-(r+1)}, \norm{\BSigma} \leq \log(n)\} \\
  &+ \sup_{x\in S_n} v^{-r}\abs{m_{\Bu, \BS, n}(z) - m_{\Bu, n}(z)}\Ind{\{|x| > {\varsigma^r}v^{-(r+1)}, \norm{\BSigma} \leq \log(n)\}}\\
  &+ \sup_{x\in S_n} v^{-r}\abs{m_{\Bu, \BS, n}(z) - m_{\Bu, n}(z)} \Ind{\{\norm{\BSigma} > \log(n)\}} := X_1 + X_2 + X_3, \numberthis
\end{align*}
where $\varsigma = n^{\alpha / (864 + 360r)}$. By inequalities $|m_{\Bu, \BS, n}(z)|, |m_{\Bu, n}(z)| \leq \norm{\Bu}^2 / v$, we have  $X_3 \leq {2 \norm{\Bu}^2}v^{-{(r+1)}}\Ind\{\norm{\BSigma} > \log(n)\}$. Write the eigen-decomposition for $\BS$ as $\BS = \BU \diag(\lambda_i; 1\leq i \leq p )\BU^H$.  When $|x| > \varsigma^r v^{-(r+1)}$ and $\norm{\BSigma} \leq \log(n)$,  we can obtain 
\begin{equation}
    \frac{1}{|\lambda_i - z|} \leq \frac{1}{|\log(n) - \varsigma^r v^{-(r+1)}|} \leq \frac{2v^{r+1}}{\varsigma^r}
\end{equation}
for sufficiently large $n$. Thus, it holds that 
\begin{equation}
|m_{\Bu, \BS, n}(z)| \leq \abs{\Bu^H \BU \diag \left( 1/|\lambda_i - z| ; 1 \leq i\leq p \right)\BU^H\Bu}  \leq \frac{2 \norm{\Bu}^2 v^{r+1}}{\varsigma^r}.
\end{equation} 
Denote the underline measure of $m_{\Bu, n}$ as $F^{\Bu, n}$. Using the same method as in the proof of Lemma \ref{Lemm_Bound_on_qterms}, it can be shown that $\int_{\mathbb{R}^+} \lambda F^{\Bu, n}(\dd \lambda) < \infty$. Then, we have 
\begin{equation}
    |m_{\Bu, n}(z)| \leq \int_{\R^+} {v^{-1}}{ \Ind_{\{\lambda \geq \varsigma^r / (2v^{r+1}) \}}}F^{\Bu, n}(\dd \lambda) + \int_{\R^+} \frac{\Ind_{\{\lambda < \varsigma^r / (2v^{r+1}) \}}F^{\Bu, n}(\dd \lambda)}{|\lambda - \varsigma^rv^{-(r+1)}|} \leq \frac{2\int_{\mathbb{R}^+} \lambda F^{\Bu, n}(\dd \lambda) v^r + 2\norm{\Bu}^2 v^{r+1}}{\varsigma^r}.
\end{equation}
As a result, we have $X_2 \leq {K_{\Bu, 2}}/ \varsigma^r$ for some constant $K_{\Bu, 2}$. Next, we evaluate $X_1$. When $|x| \leq \varsigma^r v^{-(r+1)}$, the following holds according to \eqref{Eq_Part_II} and \eqref{Eq_Part_III}
\begin{align*}
  X_1 &\leq  v^{-r}\sup_{x \in S_n} \abs{\Bu^H[\BQ - \E\BQ]\Bu + \Bu^H[\E\BQ - \BR]\Bu + \Bu^H [\BR - \BTheta]\Bu } \\
  & \leq  v^{-r}\sup_{x \in S_n}\abs{\Bu^H[\BQ - \E\BQ]\Bu} + \frac{K \norm{\Bu}^2}{v^rn^{\alpha}}\left[ \frac{1}{v^7} + \frac{(K_2 + x^2 + v^2)^{9}}{v^{29}}\right] \\
  & \leq  v^{-r}\sup_{x \in S_n}\abs{\Bu^H[\BQ - \E\BQ]\Bu} + \frac{K \norm{\Bu}^2 \varsigma^{18r}}{n^{\alpha}v^{47 + 19r}} \leq  v^{-r}\sup_{x \in S_n}\abs{\Bu^H[\BQ - \E\BQ]\Bu}+ \frac{K \norm{\Bu}^2}{n^{{\alpha/(48 + 20r)}}}, \numberthis
\end{align*}
where $K$ is a constant independent of $n$ and $K_2$ is defined in Lemma \ref{Lemm_bilinear}. We note that, unlike Lemma \ref{Lemm_bilinear}, the imaginary part of $z$ approaches 0 as  $n \to \infty$. Fortunately, the bound in Lemma \ref{Lemm_bilinear} still holds. As a result, for any $\epsilon > 0$ and sufficiently large $n$,  the following holds by Lemma \ref{Lemm_P_bound_Sn} and \eqref{Eq_Part_I}
\begin{equation}
  \begin{split}
  &\PP \left( \sup_{x \in S_n}v^{-r}\abs{m_{\Bu, \BS, n}(z) - m_{\Bu, n}(z)} \geq \epsilon \right) \leq \sum_{x \in S_n} \PP\left(v^{-r}\abs{\Bu^H[\BQ - \E\BQ]\Bu} \geq \epsilon' \right) \\
  &+ \PP\left(2{\norm{\Bu}^2v^{-(r+1)}\Ind\{\norm{\BSigma} > \log(n)\}} \geq \epsilon' \right) \leq \sum_{x \in S_n} \frac{\E\abs{\Bu^H[\BQ - \E\BQ]\Bu}^{2d}}{v^{2rd}(\epsilon')^{2d}} + o(K_{d'}n^{-d'})  \leq \frac{K_d\log^{6d}(n)n^{0.5}}{\epsilon^{2d}v^{(12 + 2r)d }n^{2\alpha d}},
\end{split}
\end{equation}
where $\epsilon' = (\epsilon - {K\norm{\Bu}^2}n^{-\alpha /(48 + 20r)} - K_{\Bu, 2} n^{-r\alpha / (864 + 360r)} ) \big/ 2 > \epsilon / 4 $, $d \geq 2 + \varepsilon / 2$, and $d' > 0$ are arbitrary.
Hence, for any positive $\epsilon$ and $l$, we have
\begin{equation}
\label{Eq_Prob_bound_on_m}
  \PP \left( \max_{x \in S_n}v^{-r}\abs{m_{\Bu, \BS, n}(z) - m_{\Bu, n}(z)} \geq \epsilon \right) \leq K_{d} \epsilon^{-d} n^{-l},
\end{equation}
for any large $d$. Next we assume that the $n$ elements of $S_n$ are equally spaced in the interval $[-\sqrt{n}, \sqrt{n}]$. Since $|x_1 - x_2| \leq 2n^{-1/2}$, it holds 
\begin{equation}
\label{Eq_max_m_1}
  \max \left\{\abs{m_{\Bu, \BS, n}(x_1 + \jmath v) - m_{\Bu, \BS, n}(x_2 + \jmath v)}, \abs{m_{\Bu, n}(x_1 + \jmath v) - m_{\Bu, n}(x_2 + \jmath v)} \right\} \leq \frac{2 \norm{\Bu}^2}{n^{1/2}v^2}, 
\end{equation}
and when $|x| > \sqrt{n}$, for $n$ large, we have 
\begin{equation}
\label{Eq_max_m_2}
  |m_{\Bu, \BS, n}(x + \jmath v)| \leq \frac{2 \norm{\Bu}^2}{n^{1/2}} + \frac{\norm{\Bu}^2 \Ind\{\norm{\BSigma} > \log(n)\}}{v} , ~~ |m_{\Bu, n}(x + \jmath v)| \leq \frac{2 \norm{\Bu}^2}{n^{1/2}} + \frac{2 \int \lambda F^{\Bu, n}(\dd \lambda )}{n^{1/2}v}.
\end{equation}
Combining \eqref{Eq_Prob_bound_on_m}, \eqref{Eq_max_m_1}, and \eqref{Eq_max_m_2}, we can obtain 
\begin{equation}
\label{Eq_Prob_bound_on_all_x}
\PP \left( \max_{x \in \mathbb{R}}v^{-r}\abs{m_{\Bu, \BS, n}(z) - m_{\Bu, n}(z)} \geq \epsilon \right) \leq K_{d} \epsilon^{-d} n^{-l}.
\end{equation}
Assume $l'$ and $l > 0$. We know $\E_j v^{-{rl'}} \sup_{x\in \mathbb{R}} |m_{\Bu, \BS, n}(z) - m_{\Bu, n}(z)|^{l'}$, $j = 0, 1, \ldots, n$, forms a discrete martingale. By letting $h = \frac{ld}{\underline{l}l'}$ with $\underline{l} > l$, we have 
\begin{equation}
  \PP\left( \max_{j \leq n} \E_j v^{-{rl'}} \sup_{x\in \mathbb{R}} |m_{\Bu, \BS, n}(z) - m_{\Bu, n}(z)|^{l'} \geq \epsilon \right) \overset{(a)}{\leq} \epsilon^{-h} \E\left[v^{-rhl'}\sup_{x\in \mathbb{R}} |m_{\Bu, \BS, n}(z) - m_{\Bu, n}(z)|^{l'h} \right] \overset{(b)}{\leq}   \frac{\epsilon^{-h}  K_d^{\frac{l}{\underline{l}}}\underline{l}}{\underline{l} - l}n^{-l},
\end{equation} 
where step $(a)$ follows from Jensen's inequality and the Maximal inequality. Step $(b)$ is due to \eqref{Eq_Prob_bound_on_all_x}. Therefore, we have completed the proof of \eqref{Eq_Bilinear_martingale_conv}. The proof for the second equation (the trace form) is quite similar and omitted here. \qed
\par
By the properties of the Stieltjes transform, the quadratic forms $\Bu^H \BQ(z)\Bu$ and $\Bu^H\BTheta(z)\Bu$ are also Stieltjes transforms and their underline positive measures are denoted by $F^{\Bu, \BS, n}$ and $F^{\Bu, n}$, respectively. In particular, denoting $\BS = \BU \diag(\lambda_i; 1 \leq i\leq p)\BU^H$ and $\BU^H\Bu = \overline{\Bu}$, the distribution $F^{\Bu, \BS, n}(x)$ is given by
\begin{equation}
  F^{\Bu, \BS, n}(x) = \sum_{j=1}^n \abs{[\overline{\Bu}]_j}^2\Ind_{\{\lambda_j \leq x\}}.
\end{equation}
Next, we use Proposition \ref{Prop_bilinear_of_resolvent_uniform_x} to give a preliminary estimate on $F^{\Bu, \BS, n}([a, b])$ and $F^{\BS, n}([a, b])$ (defined in \eqref{EQ_ESD_F_S_n}), respectively. These estimates are important for proving the no-eigenvalue property.
%%%%%%%%%%%%%%%%%%%%
\par
We evaluate $F^{\Bu, \BS, n}([a, b])$ first and keep the notations defined in the Proof of Proposition \ref{Prop_bilinear_of_resolvent_uniform_x}.
Denote $\Im m_{\Bu, \BS, n}(z) = \Im m_{\Bu, \BS, n}^{\text{out}}(z) + \Im m_{\Bu, \BS, n}^{\text{in}}(z)$, where 
\begin{align}
    \Im m_{\Bu, \BS, n}^{\text{in}}(z) = \int_{[a', b']^c} \frac{v F^{\Bu, \BS, n}(\dd \lambda)}{(\lambda - x)^2 + v^2} = \sum_{\lambda_j \in [a', b']^c} \frac{v |[\overline{\Bu}]_j|^2}{(x - \lambda_j)^2 + v^2}.
\end{align}
 Similar to $\Im m_{\Bu, \BS, n}$, denote $\Im m_{\Bu, n}(z) = \Im m_{\Bu, n}^{\text{out}}(z) + \Im m^{\text{in}}_{\Bu, n}(z) $, 
where 
\begin{equation}
    \Im m^{\text{in}}_{\Bu, n}(z) = \int_{[a', b']^c} \frac{v F^{\Bu, n}(\dd \lambda)}{(\lambda - x)^2 + v^2}.
\end{equation}
We evaluate $v^{-1}|\Im m^{\text{in}}_{\Bu, n}(z) - \Im m^{\text{in}}_{\Bu,\BS, n}(z)|$ to control the order of $F^{\Bu,\BS, n}([a, b])$. To this end, we first show that $F^{\Bu, n}([a', b']) = 0$. In fact, according to the integration representation $\BTheta(z) = \int_{\R^+} \frac{\Bmu(\dd \lambda)}{\lambda - z}$ in \eqref{Eq_representation_Theta}, we have the measure $F^n =  \Tr \Bmu / p$. Hence the interval $[a', b']$ is outside the support of the $[\Bmu]_{j, j}$ for each $j$, which implies $\Bmu([a', b']) = \mathbf{0}_{p \times p}$. As a result, it holds true that $F^{\Bu, n}([a', b']) = \Bu^H\Bmu\Bu([a', b']) = 0$ and 
$
     \Im m_{\Bu, n}^{\text{out}}(z) = 0.
$
\par
Let $l$ be a positive integer and define the sequence $\{G_q\}_{q \geq 0}$ of positive measures over $\R^l$ by
\begin{equation}
    G_{n(n+1)/2 + j}(\lambda_1, \ldots, \lambda_l) = \E_j \left(\prod_{i=1}^l F^{\Bu, \BS, n}(\lambda_i) \right), ~~ j \in [n] \cup \{0\}.
\end{equation}
Then, for $q = n(n+1)/2 + j$, the $l$-dimensional Stieltjes transform of $G_q$ is given by
$
    m_{G_q}(z_1, \ldots, z_l) = \E_j \prod_{i \leq l} m_{\Bu, \BS, n}(z_i)
$. By the normal family theorem \cite{hachem2007deterministic, bai2010spectral}, it can be verified that $G_q(\dd \lambda_1,\ldots, \dd \lambda_l) \Rightarrow \prod_{i \leq l} F^{\Bu, n}(\dd \lambda_i)$ almost surely, as $q \to \infty$. 
Notice that on $(-\infty, a'] \cup [b' + \infty)$ the set of functions in $\lambda$,
$
    S_i = \left\{ (\lambda - x)^{-i}: ~~  x \in [a, b] \right\}
$,
forms a uniformly bounded and equicontinuous family for given $i$. By \cite[Problem 8, p. 17]{billingsley2013convergence}, the following holds
\begin{align}
    \sup_{x \in [a, b]} \abs{\int_{ \left[ [a', b']^c \right]^l } \frac{G_q(\dd \lambda_1, \ldots, \dd \lambda_l) - \prod_{i=1}^l F^{\Bu, n}(\dd \lambda_i)}{\prod_{i=1}^l[(\lambda_i - x)^2 + v^2]}} \xrightarrow[q \to \infty]{a.s.} 0,  \label{Eq_G_F_Convg}
    % &\sup_{x \in [a, b]} \abs{\int_{[a', b']^c} \frac{\E_{j} F^{\Bu, \BS, n}(\dd \lambda) - F^{\Bu, n}(\dd \lambda)}{(\lambda - x)^2 + v^2} }\xrightarrow[q \to \infty]{a.s.} 0,
\end{align}
% \begin{align}
%     \sup_{x \in [a, b]} \abs{\int_{[a', b']^c \times [a', b']^c} \frac{G_q(\dd \lambda_1, \ldots, \dd \lambda_l)}{\prod_{i=1}^l[(\lambda_i - x)^2 + v^2]} - v^{-2l}\prod_{i=1}^l \Im m^{\text{in}}_{\Bu, n}(x + \jmath v) } \xrightarrow[q \to \infty]{a.s.} 0,  \label{Eq_G_F_Convg}
%     % &\sup_{x \in [a, b]} \abs{\int_{[a', b']^c} \frac{\E_{j} F^{\Bu, \BS, n}(\dd \lambda) - F^{\Bu, n}(\dd \lambda)}{(\lambda - x)^2 + v^2} }\xrightarrow[q \to \infty]{a.s.} 0,
% \end{align}
for any $v \to 0$. Since \eqref{Eq_G_F_Convg} holds for arbitrary given $l$, we can obtain
\begin{align*}
    & v^{-{2s}}\E_j \abs{\Im m^{\text{in}}_{\Bu, \BS, n}(z)- \Im m_{\Bu, n}(z)}^{2s} =  \sum_{i=0}^{2s} \binom{2s}{i}(-1)^{2s - i}\frac{\E_j[\Im m^{\text{in}}_{\Bu, \BS, n}(z)]^i}{v^i} \frac{[\Im m^{\text{in}}_{\Bu, n}(z)]^{2s-i}}{v^{2s-i}} \to 0 \numberthis
    % \to & \frac{[\Im m^{\text{in}}_{\Bu, n}(z)]^{2a}}{v^{2a}}\sum_{i=0}^{2a} \binom{2a}{i}(-1)^{2a - i} = 0,
\end{align*}
almost surely as $q \to \infty$ for any given positive integer $s$. Using the definition of the index $q$, Hölder's inequality $\E_j|xy| \leq  \E_j^{1/c}|x|^c \E_j^{1/d}|y|^d $ for positive $c, d$ with $c + d = 1$, and letting $r = 1$ in \eqref{Eq_Bilinear_martingale_conv}, we can obtain
\begin{equation}
    \sup_{x \in [a, b]}\max_{0 \leq j \leq n} v^{-l}\E_j\abs{\Im m^{\text{out}}_{\Bu, \BS, n}(z)}^l \xrightarrow[n \to \infty]{a.s.} 0,
\end{equation}
for $v = n^{-\alpha / 68}$.
Observe that for any $x \in [a, b]$, we have 
\begin{equation}
    v^{-1}\Im m^{\text{out}}_{\Bu, \BS, n}(x + \jmath v) = \int_{[a', b']} \frac{F^{\Bu, \BS, n}(\dd \lambda)}{(\lambda - x)^2 + v^2} \geq \int_{[a, b] \cap [x - v, x + v]} \frac{F^{\Bu, \BS, n}(\dd \lambda)}{(\lambda - x)^2 + v^2} \geq \frac{ F^{\Bu, \BS, n}([a, b] \cap [x - v, x + v])}{2v^2}.
\end{equation}
Select points $\{x_k\}_{k \leq L}$ such that the intervals $I_k = [x_k - v, x_k + v]$ cover $[a, b]$ and $x_{k} - x_{k-1} > v$. Hence, we have $L \leq (b - a) / v$ and
\begin{align*}
    &v^{-l} \max_{0 \leq j \leq n} \E_j \left\{ F^{\BS, n}([a, b]) \right\}^l \leq v^{-l} \max_{0 \leq j \leq n} \E_j \left\{ \sum_{k=1}^L F^{\BS, n}([a, b] \cap I_k) \right\}^l \leq 2^l \max_{0 \leq j \leq n} \E_j \left\{ \sum_{k=1}^L \Im m^{\text{out}}_{\Bu, \BS, n}(x_k + \jmath v) \right\}^l
    \\
    &\overset{(a)}{\leq} 2^l L^{l-1} \max_{0 \leq j \leq n} \sum_{k=1}^L \E_j  \left[ \Im m^{\text{out}}_{\Bu, \BS, n}(x_k + \jmath v) \right]^l \leq 2^l(b-a)^l v^{-l} \sup_{x \in [a, b]}\max_{0 \leq j \leq n}\left[ \Im m^{\text{out}}_{\Bu, \BS, n}(x + \jmath v) \right]^l \xrightarrow[n \to \infty]{a.s.} 0, \numberthis
\end{align*}
where step $(a)$ follows by Hölder's inequality. Therefore, we have proved that
\begin{align}
  \max_{0 \leq j \leq n} \E_j \left\{ F^{\Bu, \BS, n}([a, b])\right\}^l = o_{a.s.}\left(v^l\right) = o_{a.s.}\left(n^{-l\alpha / 68}\right),  ~~ \forall l \geq 1, \label{Eq_ESD_u_S_rate}
\end{align}
which give the rate for $F^{\Bu, \BS, n}([a, b])$. By using the same method, the rate for $F^{ \BS, n}([a, b])$ can be also given by
\begin{align}
  \max_{0 \leq j \leq n} \E_j \left\{ F^{\BS, n}([a, b])\right\}^l = o_{a.s.}\left(v^l\right) = o_{a.s.}\left(n^{-l\alpha / 68}\right),  ~~ \forall l \geq 1. \label{Eq_ESD_S_rate}
\end{align}
More generally, the estimation of the rate $n^{-l\alpha / 68}$ in  \eqref{Eq_ESD_u_S_rate} and \eqref{Eq_ESD_S_rate} also holds for the interval $[a', b']$.
% With the same lines as in \cite{bai1998no} and setting $r = 1$ and $\BT = \BI_p$ in Proposition \ref{Prop_bilinear_of_resolvent_uniform_x}, the following is obtained
% \begin{align}
%   &\max_{j \leq n} \E_j \left\{ F^{\BS, n}([a, b])\right\}^l = o_{a.s.}\left(v^l\right) = o_{a.s.}\left(n^{-l\alpha  / 68}\right), \label{Eq_ESD_S_rate}\\
%   &\max_{j \leq n} \E_j \left\{ F^{\Bu, \BS, n}([a, b])\right\}^l = o_{a.s.}\left(v^l\right) = o_{a.s.}\left(n^{-l\alpha / 68}\right),  ~~ \forall l \geq 1, \label{Eq_ESD_u_S_rate}
% \end{align}
% which gives the rate for $F^{\BS, n}([a, b])$ and $F^{\Bu, \BS, n}([a, b])$. More generally, the estimation of the rates in  \eqref{Eq_ESD_S_rate} and \eqref{Eq_ESD_u_S_rate} also holds for the interval $[a', b']$.
\subsection{A Rate on $F^{\Bu, \BS_j, n}([a, b])$}
Besides estimating $F^{\BS, n}([a, b])$, the evaluation for the number of eigenvalues of $\BS_j = \BSigma_j \BSigma_j^H$ is also important, where $\BSigma_j$ is the rank-one perturbation of $\BSigma$. To this end, we need the following proposition.
\begin{proposition}
\label{Prop_uniform_x_ajQjaj}
Assume the same conditions as in Proposition \ref{Prop_bilinear_of_resolvent_uniform_x}. Let $v = n^{-\alpha / 228}$ and $\Bu \in \mathbb{C}^{p}$ be a deterministic vector with uniformly bounded Euclidean norm. Then, we have
\begin{align}
  \max_{0 \leq i \leq n} \E_i \left\{ v^{-l} \sup_{x \in \mathbb{R}} \abs{ \Bu^H\left[\BQ_j(z) -  \BTheta_j(z)\right] \Bu}^l  \right\} = o_{a.s.}(1), 
\end{align}
where $\BTheta_j(z)$ is defined as
\begin{equation}
  \BTheta_j = \left[-z \left(\BI_p + \sum_{l=1}^n \frac{\BOmega_l \widetilde{\delta}_l}{n} \right) + \BA_j(\BI_{n-1} + \BD_{\delta, j})^{-1} \BA_j^H \right]^{-1},
\end{equation}
with $\BD_{\delta, j} = \diag(\delta_1, \ldots, \delta_{j-1}, \delta_{j+1}, \ldots, \delta_n)$ and $\BA_j = [\Ba_1, \ldots, \Ba_{j-1}, \Ba_{j+1}, \ldots, \Ba_n]$.
\end{proposition}
\textit{Proof: }According to the identity $\BQ = \BQ_j + z\widetilde{q}_j \BQ_j \Bxi_j\Bxi_j^H\BQ_j$ and $-z\widetilde{q}_j = \eta_j + z \widetilde{q}_j \eta_j \Delta_j$, we have 
\begin{align*}
  \Ba_j^H\BQ\Ba_j &= \frac{(1 + \alpha_{j,j}) \Ba_j^H\BQ_j\Ba_j}{1 + \alpha_{j,j} + \Ba_j^H\BQ_j\Ba_j} + \varPi_j, 
\end{align*}
where 
\begin{equation}
    \varPi_j = - z\widetilde{q}_j\eta_j\Delta_j \Ba_j^H\BQ_j \Ba_j\Ba_j^H\BQ_j\Ba_j  + z\widetilde{q}_j \Ba_j^H\BQ_j \By_j\By_j^H\BQ_j\Ba_j + z \widetilde{q}_j \Ba_j^H\BQ_j \By_j\Ba_j^H\BQ_j\Ba_j + z\widetilde{q}_j\Ba_j^H\BQ_j \Ba_j\By_j^H\BQ_j\Ba_j.
\end{equation}
Hence, we can obtain 
\begin{align*}
  \Ba_j^H\BQ\Ba_j + \frac{\Ba_j^H\BQ\Ba_j\Ba_j^H\BQ_j\Ba_j}{1 + \alpha_{j,j}}  =  \Ba_j^H\BQ_j\Ba_j+ \frac{\varPi_j(1 + \alpha_{j,j} + \Ba_j^H\BQ_j\Ba_j)}{1 + \alpha_{j,j}} = \Ba_j^H\BQ_j\Ba_j + \overline{\varPi}_j. \numberthis \label{Eq_a_j_Q_a_j}
\end{align*}
By the Woodbury matrix identity, we can get a similar identity related to $\Ba_j^H \BTheta \Ba_j$
\begin{align*}
  % &\Ba_j^H \BTheta \Ba_j = \frac{(1 + \delta_j)\Ba_j^H\BTheta_j \Ba_j}{1 + \delta_j + \Ba_j^H \BTheta_j \Ba_j} \\
  \Ba_j^H \BTheta \Ba_j + \frac{\Ba_j^H \BTheta_j \Ba_j\Ba_j^H \BTheta \Ba_j}{1 + \delta_j} = \Ba_j^H\BTheta_j \Ba_j. \numberthis \label{Eq_a_j_Theta_a_j}
\end{align*}
By subtracting \eqref{Eq_a_j_Q_a_j} and \eqref{Eq_a_j_Theta_a_j}, we obtain
\begin{align}
  \Ba_j^H(\BQ - \BTheta)\Ba_j + \frac{\Ba_j^H\BQ\Ba_j\Ba_j^H\BQ_j\Ba_j}{1 + \alpha_{j,j}} - \frac{\Ba_j^H \BTheta_j \Ba_j\Ba_j^H \BTheta \Ba_j}{1 + \delta_j}  = \Ba_j^H(\BQ_j - \BTheta_j)\Ba_j + \overline{\varPi}_j
\end{align}
As a result, it holds that 
\begin{align*}
  % &-\overline{\varPi}_j + a_j^H(Q - \Theta)a_j + \frac{a_j^HQa_ja_j^HQ_ja_j}{1 + \alpha_{j,j}} - \frac{a_j^HQa_ja_j^HQ_ja_j}{1 + \delta_j} + \frac{a_j^HQa_ja_j^HQ_ja_j - a_j^H \Theta_j a_ja_j^H \Theta a_j}{1 + \delta_j}  = a_j^H(Q_j - \Theta_j)a_j \\
  & -\overline{\varPi}_j + \Ba_j^H(\BQ - \BTheta)\Ba_j + \frac{\Ba_j^H\BQ\Ba_j\Ba_j^H\BQ_j\Ba_j}{1 + \alpha_{j,j}} - \frac{\Ba_j^H\BQ\Ba_j\Ba_j^H\BQ_j\Ba_j}{1 + \delta_j} \\
  & + \frac{\Ba_j^H\BQ\Ba_j\Ba_j^H\BQ_j\Ba_j - \Ba_j^H\BTheta \Ba_j \Ba_j^H\BQ_j\Ba_j +\Ba_j^H \BTheta \Ba_j\Ba_j^H\BQ_j\Ba_j - \Ba_j^H \BTheta_j \Ba_j\Ba_j^H \BTheta \Ba_j}{1 + \delta_j}  = \Ba_j^H(\BQ_j - \BTheta_j)\Ba_j, \numberthis
\end{align*}
which yields
\begin{align*}
  \left(1 - \frac{\Ba_j^H \BTheta \Ba_j}{1 + \delta_j}\right)
  \Ba_j^H(\BQ_j - \BTheta_j)\Ba_j &=  -\overline{\varPi}_j + \Ba_j^H(\BQ - \BTheta)\Ba_j + \frac{\Ba_j^H\BQ\Ba_j\Ba_j^H\BQ_j\Ba_j}{1 + \alpha_{j,j}} \\
  &- \frac{\Ba_j^H\BQ\Ba_j\Ba_j^H\BQ_j\Ba_j}{1 + \delta_j} + \frac{\Ba_j^H(\BQ - \BTheta)\Ba_j\Ba_j^H\BQ_j\Ba_j }{1 + \delta_j}  = \underline{\varPi}_j. \numberthis
\end{align*}
Using \eqref{Eq_a_j_Theta_a_j}, the following is obtained 
\begin{align*}
    \abs{\Ba_j^H(\BQ_j - \BTheta_j)\Ba_j}   &= \abs{ \frac{\underline{\varPi}_j (1 + \delta_j + \Ba_j^H\BTheta_j \Ba_j)}{1 + \delta_j}} \leq K{v^{-1}(1 + v^{-1}) |\underline{\varPi}_j|} \\
    & \leq K \left(v^{-4}|{\varPi}_j| +  v^{-4}|\Ba_j^H (\BQ - \BTheta) \Ba_j| + v^{-6} |\alpha_{j, j} - \delta_j| \right). \numberthis
\end{align*}
Similarly, we have 
\begin{align*}
\abs{\Ba_j^H (\BQ_j - \BTheta_j) \Bu} &\leq K\left( v^{-2}\abs{\Ba_j^H (\BQ - \BTheta) \Bu} + v^{-4} |\alpha_{j, j} - \delta_j| + v^{-2}\abs{\Ba_j^H (\BQ_j - \BTheta_j) \Ba_j} + v^{-2} |\varPi_j^{u, 1}| \right), \numberthis \\
\abs{\Bu_j^H (\BQ_j - \BTheta_j) \Ba_j} &\leq K\left( v^{-2}\abs{\Bu_j^H (\BQ - \BTheta) \Ba_j} + v^{-4} |\alpha_{j, j} - \delta_j| + v^{-2}\abs{\Ba_j^H (\BQ_j - \BTheta_j) \Ba_j} + v^{-2} |\varPi_j^{u, 2}| \right), \numberthis \\
\abs{\Bu^H (\BQ_j - \BTheta_j) \Bu} &\leq K \Big( \abs{\Bu^H( \BQ - \BTheta) \Bu} + v^{-4} |\alpha_{j,j} - \delta_j| + v^{-4}|\Ba_j^H(\BTheta_j - \BQ_j) \Ba_j| \\
&+ v^{-2}\abs{\Ba_j^H(\BTheta_j - \BQ_j) \Bu} + v^{-2}\abs{\Bu^H(\BTheta_j - \BQ_j) \Ba_j}+ \abs{\varPi_j^{u, u}} \Big), \numberthis
\end{align*}
where $\varPi^{u,1}_j$, $\varPi^{u,2}_j$, and $\varPi^{u,u}_j$ have similar forms as $\varPi_j$.
Letting $r = 9$ in Proposition \ref{Prop_bilinear_of_resolvent_uniform_x}, we can get the probabilistic bound $\PP(\max_{x \in S_n} v^{-1}\abs{\Bu^H(\BQ_j - \BTheta_j)\Bu} \geq \epsilon) \leq K_d \epsilon^{-d} n^{-l}$, which implies 
\begin{align}
  \max_{0 \leq i \leq n} \E_i \left\{ v^{-l} \sup_{x \in \mathbb{R}} \abs{ \Bu^H\left[\BQ_j -  \BTheta_j\right] \Bu}^l  \right\} = o_{a.s.}(1).
\end{align}
Therefore, Proposition \ref{Prop_uniform_x_ajQjaj} is proved. \qed
\par
By the properties of the Stieltjes transform, it can be verified that $\BTheta_j$ is a matrix-valued Stieltjes transform and  $\Bu^H \BTheta_j(z) \Bu$ is the Stieltjes transform of a non-negative measure denoted by $F^{\Bu, nj}$. 
By Proposition \ref{Prop_same_support}, $[a', b']$ is outside the support of the underline measure of $\delta_j$.
As a result, for any $x_1, x_2 \in [a', b']$, we have 
\begin{equation}
   \sup_{x_1, x_2 \in [a', b']}\frac{\abs{ \delta_j(x_1) - \delta_j(x_2)}}{|x_1 - x_2|} \leq \sup_{x_1, x_2 \in [a', b']}\int_{\R^+} \frac{\mu_j(\dd \lambda)}{\abs{x_1 - \lambda}{\abs{x_2 - \lambda}}} \leq \frac{ \mu_j(\R^+)}{\dist^2([a', b'], \mathsf{Supp}(\mu_j))} \leq K,
\end{equation}
where $K$ is a constant that is independent of $j$ and $n$.
Similarly, we can also get $\sup_{x_1, x_2 \in [a', b']}|\widetilde{\delta}_j(x_1) - \widetilde{\delta}_j(x_2)| / |x_1 - x_2| \leq K$.
Therefore, we can assume $\liminf_n \inf_{x \in [a', b']} \min_{j \leq n} \{|\widetilde{\delta}_j(x)|, |1 + \delta_j(x)|\} > 0$ without loss of generality. Letting $z \in \mathbb{C}^+ \to x \in [a', b']$ and by the relation 
\begin{align*}
    \Bu^H \BTheta \Bu = \Bu^H \BTheta_j \Bu + z\widetilde{\delta}_j \Bu^H \BTheta_j \Ba_j \Ba_j^H \BTheta_j \Bu = \Bu^H \BTheta_j \Bu + \frac{\Bu^H \BTheta \Ba_j \Ba_j^H \BTheta \Bu}{z\widetilde{\delta}_j(1 + \delta_j)^2}, \numberthis
\end{align*}
% \begin{equation}
%     \Ba_j^H \BTheta_j(z) \Ba_j = \frac{-\Ba_j^H \BTheta(z) \Ba_j}{(1 + \delta_j(z))z\widetilde{\delta}_j(z)},
% \end{equation}
we have $\lim_{z \in \mathbb{C}^+ \to x} \Im(\Bu^H\BTheta_j(z)\Bu) = 0$ since $|1 + \delta_j(x)|, |\widetilde{\delta}_j(x)| > 0$, $\lim_{z \in \mathbb{C}^+ \to x} \BTheta(z) = \BTheta_0$ exists, and $\BTheta_0$ is Hermitian. In the above derivation, we also used the result that $x$ is outside the support of the underlying matrix-valued measure of $\BTheta(z)$, which is proved in Section \ref{Sec_A_rate_one_F}. As a result, it holds that $F^{\Bu, nj}([a', b']) = 0$.
Recalling that $F^{\Bu ,\BS_j ,n}$ is the underline measure of the Stieltjes transform $\Bu^H \BQ_j \Bu$,
Proposition \ref{Prop_uniform_x_ajQjaj} together with $F^{\Bu, nj}([a', b']) = 0$ implies
\begin{align}
  &\max_{ 0\leq i \leq n, 1 \leq j \leq n} \E_i \left\{ F^{\Bu, \BS_j, n}([a', b'])\right\}^l = o_{a.s.}\left(n^{-l\alpha  / 228}\right), ~~ \forall l \geq 1. \label{Eq_ESD_aj_Sj_rate}
\end{align}
\subsection{Convergence of $|m_{\BS_n} - \E m_{\BS, n}|$}
\label{Sec_Conv_random_part}
 In this section, we take $v = n^{-\alpha / 912}$ and prove that
\begin{equation}
  \sup_{x \in [a, b]} pv\abs{m_{\BS, n} - \E m_{\BS, n}} \xrightarrow{a.s.} 0. \label{Eq_pv_Diff_m_as_0}
\end{equation}
Following the same argument as Proposition \ref{Prop_bilinear_of_resolvent_uniform_x}, we note that it is sufficient to prove
\begin{equation}
  \sup_{x \in S_n} pv\abs{m_{\BS, n} - \E m_{\BS, n}} \xrightarrow{a.s.} 0,
\end{equation}
where $S_n = \{a + t(b - a) / (n^2 - 1), t = 0, \ldots, n^2-1\}$.
To proceed, we define the following quantities
\begin{align*}
  &\gamma_j = \frac{1}{1 + \frac{1}{n} \Tr \BOmega_j \E \BQ_j + \Ba_j^H \E \BQ_j \Ba_j}, ~~ \Xi_j = \frac{1}{n}\Tr \BOmega_j(\BQ_j - \E \BQ_j) + \Ba_j^H (\BQ_j - \E \BQ_j) \Ba_j \\
  & \underline{\Delta}_j = \Bxi_j^H\BQ_j^2 \Bxi_j - \frac{1}{n} \Tr \BOmega_j \BQ_j^2 - \Ba_j^H \BQ_j^2 \Ba_j. \numberthis \label{Eq_Def_2_gamma_Delta_xi}
\end{align*}
The following proposition shows the bound for the terms defined in \eqref{Label_Def_of_terms_bilinear} and \eqref{Eq_Def_2_gamma_Delta_xi} with $x \in [a, b]$. These bounds will be frequently used in the proof.
\begin{proposition}
\label{Prop_bound_ab_K_o}
For $z = x + \jmath v = x + \jmath n^{-\alpha / 912}$, there exists constant $K^o > 0$ such that for any $x \in [a, b]$ and sufficiently large $n$, the following holds.
\begin{equation}
  \max_{j \leq n}\left\{ |\delta_j|, | z\widetilde{\delta}_j|, |\kappa_j|,  |z\widetilde{\kappa}_j|, |\E\eta_j| , |\gamma_j|, \frac{1}{|1 + \delta_j|}, \frac{1}{|1 + \kappa_j|}, \frac{1}{|1 + \E \alpha_{j, j}|} \right\} \leq K^o.
\end{equation}
\end{proposition}
\textit{Proof:} By Proposition \ref{Prop_same_support}, we know the interval $[a', b']$ is outside the support of $\mu_j$, $\widetilde{\mu}_j$. Hence, we have
\begin{equation}
    |\delta_j(z)| \leq \int_{\R^+} \frac{\mu_j(\dd \lambda)}{|\lambda - z|} \leq \frac{\mu_j(\R)}{\underline{\epsilon}^2} \leq \frac{\Tr \BOmega_j}{n\underline{\epsilon}^2},
\end{equation}
where we recall that $\underline{\epsilon} = |a - a'| = |b - b'|$.
The inequality $|\widetilde{\delta}_j| \leq \underline{\epsilon}^{-2}$ can be shown similarly. By Lemma \ref{Lemm_bilinear}, we know 
$|\delta_j - \kappa_j |$, $|z\widetilde{\delta}_j - z\widetilde{\kappa}_j|$, $|z\widetilde{\delta}_j + \E \eta_j|$, and $|z\widetilde{\delta}_j + \gamma_j|$ are of order $o(1)$ uniformly for $x \in [a, b]$, which proves the bounds for these terms.
The bound for $1 / |1 + \delta_j|$ is due to the assumption that $\liminf_n \inf_j |1 + \delta_j| > 0$ for $x \in [a, b]$. Since $|\delta_j - \kappa_j|, |\kappa_j - \E\alpha_{j, j}| = o(1)$, $1 / |1 + \kappa_j|, 1 / |1 + \E \alpha_{j, j}| < \infty$ are proved.  \qed
\par
Writing $p\E m_{\BS, n} - pm_{\BS, n}$ as the sum of the martingale difference sequence and using $-z\widetilde{q}_j = \eta_j + z\widetilde{q}_j \Delta_j  \eta_j$, we have 
\begin{align*}
&p \E m_{\BS, n} - p m_{\BS, n}  = \sum_{j=1}^n [\E_{j} - \E_{j-1}] -z\widetilde{q}_j \Bxi_j^H \BQ_j^2 \Bxi_j \\
&= \sum_{j=1}^n [\E_{j} - \E_{j-1}] \eta_j \Bxi_j^H \BQ_j^2 \Bxi_j + [\E_{j} - \E_{j-1}] z\widetilde{q}_j \eta_j \Delta_j \Bxi_j^H \BQ_j^2 \Bxi_j \\
& = \sum_{j=1}^n [\E_{j} - \E_{j-1}] \eta_j \Bxi_j^H \BQ_j^2 \Bxi_j + [\E_{j} - \E_{j-1}] z\widetilde{q}_j \eta_j \Delta_j \left(\Bxi_j^H \BQ_j^2 \Bxi_j - \frac{1}{n} \Tr \BOmega_j \BQ_j^2 - \Ba_j^H\BQ_j^2 \Ba_j\right) \\
&+ [\E_{j} - \E_{j-1}]z\widetilde{q}_j \eta_j \Delta_j \left(\frac{1}{n}\Tr \BOmega_j \BQ_j^2 + \Ba_j^H\BQ_j^2\Ba_j\right)  = \sum_{j=1}^n \varGamma_{1j} + \varGamma_{2j} + \varGamma_{3j}. \numberthis
\end{align*}
Let $F^{\BS_j, n}$ be the ESD of $\BS_j = \BSigma_j \BSigma_j^H$. 
% By \eqref{Eq_ESD_aj_Sj_rate}, \eqref{Eq_ESD_S_rate} and \cite[Lemma 2.12]{bai1998no}, we get 
% $
%     \max_{i, j \leq n} \E_i\left\{ F^{\BS_j, n}([a', b'])\right\}^2 = o(v^8)
% $ and $$.
Define the indicator function 
\begin{equation}
  \mathscr{I}_j = \Ind_{\left\{\E_{j-1} \left[F^{\BS_j, n}([a', b']) + F^{\Ba_j, \BS_j, n}([a', b'])\right] \leq v^4\right\} \cap \left\{\E_{j-1} \left[F^{\BS_j, n}([a', b']) + F^{\Ba_j, \BS_j, n}([a', b'])\right]^2 \leq v^8\right\}}.
\end{equation}
Then, by \eqref{Eq_ESD_u_S_rate}, \eqref{Eq_ESD_S_rate}, \eqref{Eq_ESD_aj_Sj_rate}, and \cite[Lemma 2.12]{bai1998no}, we can obtain $\PP(\cup_{j=1}^n \{\mathscr{I}_j = 0\} \mathrm{~~ i.o.}) = 0$. As a result, we have, for any $\epsilon > 0$,
\begin{align*}
    \PP \left( v \max_{x \in S_n} \abs{\sum_{j=1}^n \varGamma_{1j}} > \epsilon \mathrm{~~ i.o.} \right) &\leq \PP \left( \left(\left[\max_{x \in S_n} v\abs{\sum_{j=1}^n \varGamma_{1j}}> \epsilon \right]\right) \bigcap_{j=1}^n [\mathscr{I}_j = 1]  \bigcup\left(\bigcup_{j=1}^n[\mathscr{I}_j = 0] \right)\mathrm{~~ i.o.} \right) \\
    &\leq \PP  \left(\max_{x \in S_n} v\abs{\sum_{j=1}^n \varGamma_{1j} \mathscr{I}_j}> \epsilon \mathrm{~~ i.o.}\right). \numberthis \label{Eq_Probabilistic_io_varGamma1}
\end{align*}
By the independence between $\By_j$ and $\mathscr{F}_{j-1}$, we can obtain 
\begin{equation}
  \varGamma_{1j}\mathscr{I}_j = \E_{j}\mathscr{I}_j\eta_j \left(\By_j^H \BQ_j^2 \By_j - \frac{1}{n}\Tr \BOmega_j \BQ_j^2\right) + \E_j \mathscr{I}_j\eta_j \Ba_j^H \BQ_j^2 \By_j + \E_j\mathscr{I}_j \eta_j \By_j^H \BQ_j^2 \Ba_j = \chi_{1j} + \chi_{2j} + \chi_{3j}.
\end{equation}
Since $\mathscr{I}_j \in\mathscr{F}_{j-1}$, $\{ \chi_{ij} \}_{1 \leq j\leq n}$ forms martingale difference sequences for $i=1, 2, 3$, Lemma \ref{Lemm_Burk_Ineq} shows
\begin{align*}
  \E \abs{\sum_{j=1}^nv\varGamma_{1j}\SI_j}^{2t}  \lesssim_t \sum_{i=1}^3\left[  \E \left( \sum_{j=1}^n \E_{j-1}|v\chi_{ij}|^2 \right)^{t} + \sum_{j=1}^n \E \abs{v\chi_{ij}}^{2t} \right], t \geq 3. \numberthis \label{Eq_Burk_varGamma1j_chi}
\end{align*}
Next, we evaluate the terms on the RHS of \eqref{Eq_Burk_varGamma1j_chi} that are related to $\chi_{1j}$. By writing $1 = \Ind_{\{|\eta_j| > 2K^o \}} + \Ind_{\{|\eta_j| \leq  2K^o \}}$, we have
\begin{align*}
  &v^2\sum_j \E_{j-1} |\chi_{1j}|^2 \leq v^2\sum_j \E_{j-1} \SI_j|\eta_j|^2 \abs{\By_j^H \BQ_j^2 \By_j - \frac{1}{n}\Tr \BOmega_j \BQ_j^2}^2 \lesssim \sum_j\frac{v^2}{n^2} \E_{j-1} \SI_j|\eta_j|^2 \Tr \BOmega_j \BQ_j^2 \BOmega_j \BQ_j^{H, 2} \\
  & \lesssim \sum_j\frac{v^2}{n^2} \E_{j-1} \SI_j \Tr \BOmega_j \BQ_j^2 \BOmega_j \BQ_j^{H, 2} + \frac{1}{nv^4} \E_{j-1} \Ind_{\{|\eta_j| > 2K^o \}}. \numberthis \label{Eq_I_chi_1j}
\end{align*}
Denoting $\lambda_{jl}$ as the $l$-th smallest eigenvalue of $\BS_j$, we have, on the set $\{ \SI_j = 1 \}$,
\begin{align*}
  &\E_{j-1}\Tr \BOmega_j \BQ_j^2 \BOmega_j \BQ_j^{H, 2} \leq (\omega^+)^2 \Tr \BQ_j^2\BQ_j^{H, 2} = (\omega^+)^2 \E_{j-1}\sum_{l=1}^p \frac{1}{((\lambda_{jl} - x)^2 + v^2)^2} \\
  &= (\omega^+)^2 \E_{j-1}\left[ \sum_{\lambda_{jl} \in [a', b']} \frac{1}{((\lambda_{jl} - x)^2 + v^2)^2} +  \sum_{\lambda_{jl} \notin [a', b']} \frac{1}{((\lambda_{jl} - x)^2 + v^2)^2} \right] \\
  &\leq  (\omega^+)^2 \left[p v^{-4} \E_{j-1}F^{\BS_{j}, n}([a', b']) + p \underline{\epsilon}^{-4} \right]  \leq Kn, \numberthis \label{Eq_I_j_TrQQ}
\end{align*}
for some constant $K$ and $\omega^+ := \sup_{j, n} \norm{\BOmega_j}$. Moreover, $|\eta_j| > 2K^o$ implies that $|\gamma_j^{-1} + \Xi_j| < 1 / (2K^o)$ and  thus the following holds by $|\gamma_j| \leq K^o$
\begin{equation}
    \Ind_{\{ |\eta_j| > 2K^o\}} \leq  \Ind_{\{ |\gamma_j^{-1} + \Xi_j| < 1 / (2K^o)\}} \leq \Ind_{\{|\Xi_j| \geq 1/(2K^o)\}}.
\end{equation}
As a result, by $(\sum_{j=1}^n |a_j|)^t \leq (\sum_{j=1}^n |a_j|^t)n^{t-1}$, we get
\begin{align*}
  &\E \left( v^2\sum_{j=1}^n \E_{j-1} |\chi_{1j}|^2 \right)^{t}  \lesssim v^{2t} + \frac{1}{n^{t}v^{4t}} \E \left(\sum_{j=1}^n \E_{j-1} \Ind_{\{|\Xi_j| \geq 1/(2K^o) \}} \right)^t \\
  & \leq v^{2t} + \frac{\sum_{j=1}^n \PP (|\Xi_j| \geq  1 / (2K^o))}{nv^{4t}} \leq  v^{2t} + \frac{(2K^o)^m\sum_{j=1}^n \E|\Xi_j|^m}{nv^{4t}} \overset{(a)}{\lesssim_t} v^{2t}, \numberthis \label{Eq_E_v2_E_j_x_1j2}
\end{align*}
where step $(a)$ follows by choosing a large $m$ and applying \eqref{Eq_Part_I}. By Lemma \ref{Lemm_trace}, we can obtain 
\begin{align*}
  \sum_{j=1}^n \E |\chi_{1j}|^{2t} \leq \sum_{j=1}^n \E \abs{\eta_j}^{2t}\abs{\By_j^H \BQ_j^2 \By_j - \frac{1}{n}\Tr \BOmega_j \BQ_j^2}^{2t} \lesssim_t \frac{1}{v^{2t}}\sum_{j=1}^n \frac{1}{n^{4\alpha t + 1} v^{4t}} = \frac{1}{n^{4\alpha t} v^{6t}}. \numberthis
\end{align*}
Next, we estimate the terms related to $\chi_{2j}$ on the RHS of \eqref{Eq_Burk_varGamma1j_chi}. Similar to \eqref{Eq_I_chi_1j} and using $\Ba^H \BC \Ba \leq \Ba^H \Ba \norm{\BC}$, we have
\begin{align*}
  &v^2 \sum_{j=1}^n \E_{j-1} |\chi_{2j}|^2 \leq \sum_{j=1}^n  v^2\E_{j-1} \mathscr{I}_j|\eta_j \Ba_j^H \BQ_j^2 \By_j|^2 \leq \sum_{j=1}^n \frac{v^2}{n}\E_{j-1} \mathscr{I}_j |\eta_j|^2 \Ba_j^H\BQ_j^2\BOmega_j\BQ_j^{2, H}\Ba_j \\
  &\lesssim \sum_{j=1}^n\frac{v^2}{n}\E_{j-1} \mathscr{I}_j |\eta_j|^2 {\Ba_j^H\BQ_j^2\BQ_j^{2, H}\Ba_j} \lesssim \sum_{j=1}^n \frac{v^2}{n} \E_{j-1} \mathscr{I}_j \Ba_j^H\BQ_j\BQ_j^{2, H} \Ba_j + \frac{1}{nv^4} \E_{j-1 } \Ind_{\{|\Xi_j| > 1/(2K_o)\}}. \numberthis \label{Eq_v2_E_j_chi2j2}
\end{align*}
Write the singular value decomposition of $\BS_j$ as $\BU_j \diag(\lambda_{jl}; 1 \leq l \leq p) \BU_j^H$ and  denote $\overline{\Ba}_j = \BU_j^H\Ba_j$. Here, on the set $\SI_j$, the following estimation holds by \eqref{Eq_ESD_aj_Sj_rate} 
\begin{align*}
  & \E_{j-1} \Ba_j^H\BQ_j\BQ_j^{2, H} \Ba_j = \E_{j-1} \left[ \sum_{\lambda_{jl} \in [a', b']} \frac{|[\overline{\Ba}_j]_l|^2}{((x - \lambda_{jl})^2 + v^2)^2} + \sum_{\lambda_{jl} \notin [a', b']} \frac{|[\overline{\Ba}_j]_l|^2}{((x - \lambda_{jl})^2 + v^2)^2}  \right]  \\
  & \leq  \underline{\epsilon}^{-4}{\norm{\Ba_j}^2} + v^{-4}\E_{j-1}F^{\Ba_j, \BS_j, n}([a', b']) \leq K, \numberthis \label{Eq_E_j_ajQjQjaj}
\end{align*}
with some constant $K$.  Therefore,  substituting \eqref{Eq_E_j_ajQjQjaj} into equation \eqref{Eq_v2_E_j_chi2j2}, and using the same method as in \eqref{Eq_E_v2_E_j_x_1j2}, we obtain
\begin{equation}
    \E \left( v^2 \sum_{j=1}^n \E_{j-1} |\chi_{2j}|^2 \right)^{2t} \lesssim_t v^{2t}. \label{Eq_burk_chi_2j_bound_1}
\end{equation}
The following can be derived by Lemma \ref{Lemm_trace}
\begin{equation}
    \sum_{j=1}^n \E |\chi_{2j}|^{2t} \lesssim_t \frac{1}{n^{2 \alpha t + 1} v^{6t}}. \label{Eq_burk_chi_2j_bound_2}
\end{equation}
The evaluation of the terms related to $\chi_{3j}$ is similar to \eqref{Eq_burk_chi_2j_bound_1} and \eqref{Eq_burk_chi_2j_bound_2}. Therefore, we get 
\begin{equation}
    \PP \left(\max_{x \in S_n}\abs{\sum_{j=1}^nv\varGamma_{1j}\SI_j} \geq \epsilon \right) \leq \frac{1}{\epsilon^{2t}}\sum_{x \in S_n}\E \abs{\sum_{j=1}^nv\varGamma_{1j}\SI_j}^{2t} \leq K_{t, \epsilon}n^{2-2t\alpha/912}, ~~ t \geq 3.
\end{equation}
By Borel-Cantelli Lemma and \eqref{Eq_Probabilistic_io_varGamma1}, we know  $\max_{x \in S_n}|\sum_{j}v\varGamma_{1j}| \to 0$ almost surely.
\par 
Next, we handle $\sum_j \varGamma_{2j}$. We first give the following probabilistic bound
\begin{equation}
  \PP \left( |z\widetilde{q}_j \eta_j \Delta_j| > n^{-\alpha / 3} \right) \leq \frac{\E\left[|z\widetilde{q}_j| |\eta_j| |\Delta_j|\right]^r}{n^{-\frac{ \alpha r}{3}}} \lesssim \frac{ \E^{\frac{1}{2}}{\norm{\BSigma}}^{4r} \E^\frac{1}{2} |\Delta_j|^{2r} }{v^{2r}n^{-\frac{ \alpha {r}}{3}}} \lesssim_r \frac{1}{n^{2 \alpha r / 3}v^{3r}},
\end{equation}
for large $r$. Hence, we have
\begin{align*}
  \sum_{j=1}^n \E_{j-1} \abs{v\varGamma_{2j}}^2 &\lesssim  v^2\sum_{j=1}^n \E_{j-1}\abs{z \widetilde{q}_j \eta_j \Delta_j \underline{\Delta}_j}^2 \lesssim n^{-{2\alpha}/3}\sum_{j=1}^n v^2\E_{j-1}\abs{ \underline{\Delta}_j}^2 + \sum_{j=1}^n \frac{n^8}{v^8}\E_{j-1} \Ind{\{|z\widetilde{q}_j \eta_j \Delta_j| > n^{-\alpha / 3}\}} \\
  & \lesssim \frac{1}{n^{{2\alpha}/3}v^{2}} + \frac{n^8}{v^8} \sum_{j=1}^n \E_{j-1}\Ind{\{|z\widetilde{q}_j \eta_j \Delta_j| > n^{- \alpha/3}\}}, \numberthis
\end{align*}
since $\norm{\Bxi_j} \leq \sqrt{n}$, $\norm{\BSigma} \leq n$, $|z \widetilde{q}_j| \lesssim n^2 / v$, $|\Delta_j| \lesssim n / v$, and $|\underline{\Delta}_j| \lesssim n / v^2$. As a result, we can obtain
\begin{align*}
  \E \abs{ \sum_{j=1}^n v \varGamma_{2j} }^{2t} \lesssim_t \sum_{j=1}^n  \E \left[\sum_{j=1}^n\E_{j-1}\abs{v\varGamma_{2j}}^2 \right]^t + \sum_{j=1}^n \E|v \varGamma_{2j}|^{2t} \lesssim \frac{1}{n^{2\alpha t / 3}v^{2t}} + \frac{n^{9t}}{v^{8t}v^{3r}n^{2\alpha r / 3}} + \frac{1}{v^{8t} n^{2 \alpha t}}. \numberthis \label{Eq_varGamma_2_bound}
\end{align*}
 By choosing sufficiently large $t$ and $r$ in the above inequality, we have $\max_{x \in S_n}|\sum_{j}\varGamma_{2j}|  = o_{a.s.}(v^{-1})$.
 \par
 In the following, we handle $\sum_j \varGamma_{3j}$. By the similar argument as in \eqref{Eq_Probabilistic_io_varGamma1}, it is sufficient to handle $\sum_j \varGamma_{3j} \SI_j$. Using the same method as in \eqref{Eq_I_chi_1j}, we get 
\begin{align*}
  \sum_{j=1}^n \E_{j-1}\abs{v\SI_j \varGamma_{3j}}^2 &\lesssim v^2\sum_{j=1}^n \SI_j \E_{j-1} \abs{z \widetilde{q}_j \eta_j}^2\abs{ \frac{1}{n} \Tr \BOmega_j \BQ_j^2 + \Ba_j^H\BQ_j^2\Ba_j}^2 |\Delta_j|^2 \\
  &\lesssim v^2 \sum_{j=1}^n \SI_j \E_{j-1}\abs{ \frac{1}{n} \Tr \BOmega_j \BQ_j^2 + \Ba_j^H\BQ_j^2\Ba_j}^2 |\Delta_j|^2 + \frac{n^6}{v^8}\E_{j-1} \Ind_{\{|\eta_j| \geq 2K^o\} \cup \{|z\widetilde{q}_j| \geq 2K^o\}}. \numberthis \label{Eq_varGamma_3j_split}
\end{align*}
By using Cauchy-Schwarz inequality $|\Tr \BC\BD^H|^2 \leq \Tr\BC\BC^H \Tr \BD \BD^H $ and $\Ba^H\Bb \leq \norm{\Ba}^2 \norm{\Bb}^2$, we have 
% \begin{align*}
%     \abs{ \frac{1}{n} \Tr \BOmega_j \BQ_j^2 + \Ba_j^H\BQ_j^2\Ba_j}^2 &\leq 2 \left(\frac{1}{n} \Tr \BOmega_j^{1/2} \BQ_j\BOmega_j^{1/2} \BQ_j^H\right)^2 + 2\left(\Ba_j^H \BQ_j \BQ_j^H \Ba_j \right)^2. \\
%     &\leq  2 (\omega^+)^2 \left(\frac{1}{n} \Tr \BQ_j \BQ_j^H\right)^2 + 2\left(\Ba_j^H \BQ_j \BQ_j^H \Ba_j \right)^2
% \end{align*}
\begin{align*}
    \abs{ \frac{1}{n} \Tr \BOmega_j \BQ_j^2 + \Ba_j^H\BQ_j^2\Ba_j}^2 &\leq 2 \abs{\frac{1}{n} \Tr \BOmega_j \BQ_j^2}^2 + 2\Ba_j^H \BQ_j^2 \Ba_j \Ba_j^H \BQ_j^{H, 2} \Ba_j  \\
    &\leq  2 \frac{\Tr \BOmega_j^2}{n^2} \Tr \BQ_j^2 \BQ_j^{H, 2} + 2\norm{\Ba_j}^2 \Ba_j^H \BQ_j^2 \BQ_j^{H, 2} \Ba_j. \numberthis
\end{align*}
Then, we can obtain, for $x \in [a, b]$,
\begin{align*}
    &v^2 \sum_{j=1}^n \SI_j \E_{j-1}\abs{ \frac{1}{n} \Tr \BOmega_j \BQ_j^2 + \Ba_j^H\BQ_j^2\Ba_j}^2 |\Delta_j|^2 \\
    &\lesssim \frac{v^2}{n} \sum_{j=1}^n \SI_j \E_{j-1}\left[ \frac{1}{p} \Tr \BQ_j^2\BQ_j^{H,2} + \Ba_j^H \BQ_j^2 \BQ_j^{H, 2} \Ba_j  \right] \left[\frac{1}{p} \Tr \BQ_j\BQ_j^H + \Ba_j^H \BQ_j\BQ_j^H \Ba_j \right] \\
    &\leq \frac{v^2}{n} \sum_{j=1}^n \E_{j-1} \SI_j \Big\{  \left[\left(1 + \norm{\Ba_j}^2 \right)\underline{\epsilon}^{-4} + v^{-4}F^{\BS_j, n}([a', b']) + v^{-4} F^{\Ba_j, \BS_j, n}([a', b']) \right] \\
   &~~~  \times \left[\left(1 + \norm{\Ba_j}^2 \right) \underline{\epsilon}^{-2} + v^{-2} F^{\BS_j, n}([a', b']) + v^{-2} F^{\Ba_j, \BS_j, n}([a', b']) \right] \Big\}  \lesssim v^2. \numberthis \label{Eq_F_aj_Sj_F_Sj_2}
\end{align*}
As a result, we have 
\begin{equation}
    \E\left( \sum_{j=1}^n \E_{j-1}\abs{v\SI_j \varGamma_{3j}}^2\right)^t \lesssim_t v^{2t} + {v^{-8t}n^{7t-1} \sum_{j=1}^n \left[ \PP(|\Xi_j| > 1/(2K^o)) + \PP(|\Delta_j + \Xi_j| > 1/(2K^o)) \right] } \lesssim_t v^{2t}.
\end{equation}
Therefore, the following bound holds
\begin{align*}
    \E\abs{\sum_{j=1}^n v \SI_j \varGamma_{3j}}^{2t} \lesssim_t v^{2t} + \sum_{j=1}^n \E |v\varGamma_{3j}|^{2t} \lesssim_t v^{2t} + \frac{\log^{4t}(n)}{v^{8t}n^{2 \alpha t}}, \numberthis
\end{align*}
which implies $\max_{x \in S_n} |\sum_j \varGamma_{3j}| = o_{a.s.}(v^{-1})$. Gathering the above analysis, we can prove that 
\begin{equation}
    \sup_{x \in S_n} \abs{ \sum_{j=1}^n \varGamma_{1j} + \varGamma_{2j} + \varGamma_{3j} } = o_{a.s.}(v^{-1}), 
\end{equation}
which completes the proof of \eqref{Eq_pv_Diff_m_as_0}.
\subsection{Convergence of $\E m_{\BS, n} - m_n$} In this section, for $z = x + \jmath v = x + \jmath n^{-\alpha / 912}$, we prove 
\begin{align*}
  \sup_{x \in [a, b]} p \abs{\E m_{\BS, n}(z) - m_n(z)}  \lesssim 1. \numberthis \label{Eq_Mean_Part}
\end{align*}
We begin by proving the following bound
\begin{equation}
    \sup_{x \in [a, b]} \norm{\BR(z)} \lesssim 1. \label{Eq_sup_ab_norm_R}
\end{equation}
According to discussion in Section \ref{Sec_A_rate_one_F}, we know $\Bmu([a, b]) = \mathbf{0}_{p \times p}$.  To get the bound for $\norm{\BTheta(z)}$, we choose a sequence of simple functions  \cite{rudin1987real} $\{f_k\}_{k \geq 1}$, where $f_k$ is a piece-wise constant function with $f_{k}(\lambda) = \sum_{l=1}^{N_k} a_{k, l} \Ind\{\lambda \in A_{k, l}\}$, $|f_k(\lambda)| \leq  |1 / (\lambda - z)|$, and $f_k(\lambda) \to 1 / (\lambda - z)$ as $k \to \infty$ for any $\lambda \in \R^+$. By the dominant convergence theorem, we have
\begin{align*}
    &\norm{\BTheta(z)} = \norm{\int_{\R^+}\frac{\Bmu(\dd \lambda)}{\lambda - z}} =\norm{\int_{\R^+}\lim_{k \to \infty} f_k \Bmu(\dd \lambda)}  = \norm{\lim_{k \to \infty} \sum_{l=1}^{N_k} a_{k, l}{\Bmu(A_{k, l})}} \\
    &\overset{(a)}{\leq} \norm{\lim_{k \to \infty} \sum_{l=1}^{N_k} |a_{k, l}|{\Bmu(A_{k, l})} } = \norm{\lim_{k \to \infty}\int_{\R^+} |f_k|{\Bmu(\dd \lambda)}} = \norm{\int_{\R^+} \frac{\Bmu(\dd \lambda)}{|\lambda - z|}} \leq \frac{1}{\underline{\epsilon}}, \numberthis \label{Eq_sup_ab_norm_Theta}
\end{align*}
where step $(a)$ follows by Lemma \ref{Lemm_Spnorm_Bound_of_sum_semidef_matrices}. Here, we used the fact that $\sup_{x \in [a, b]} |1 / (\lambda - z)| \leq 1 / \underline{\epsilon}$ since $[a, b] \subset [a', b']$.
By \eqref{Eq_Boud_norm_diff_R_Theta}, we can get $\sup_{x \in [a, b]} \norm{\BR - \BTheta} = o(1)$. Hence, \eqref{Eq_sup_ab_norm_R} holds. 
\begin{proposition}
\label{Prop_EQ_ER_EwQ_EwR}
    Let $\BD_T \in \mathbb{C}^{p \times p}$ be a diagonal matrix and $\BT \in \mathbb{C}^{p \times p}$ be a Hermitian non-negative matrix such that $\norm{\BD_T}, \norm{\BT} \leq T < \infty$. Then, we have  
\begin{align*}
  &\sup_{x \in [a, b]} \abs{ \Tr \BT \left( \E \BQ - \BR \right) } \lesssim 1, \numberthis \label{Eq_EQ_R}\\
  &\sup_{x \in [a, b]} \abs{z\Tr \BD_T(\E\widetilde{\BQ} - \widetilde{\BR}) } \lesssim 1. \numberthis \label{Eq_EwQ_wR}
\end{align*}
\end{proposition}
\textit{Proof: } The bound in \eqref{Eq_EQ_R} can be obtained by following the similar approach as in Section \ref{Para_Part_II}, i.e., by using the resolvent identity to break the target into a sum of terms with means close to 0 and finite moments. The difference here is that we need to bound the spectral norm of $\BR$ by \eqref{Eq_sup_ab_norm_R} for $x \in [a, b]$, as well as the following bound which can be obtained by the same argument as in \eqref{Eq_I_j_TrQQ}, \eqref{Eq_E_j_ajQjQjaj}, and \eqref{Eq_F_aj_Sj_F_Sj_2}
\begin{equation}
    \sup_{x \in [a, b]}\E \left[ \frac{1}{p} \Tr \BQ_j^{t_1} \BQ_j^{t_1, H} + \Ba_j^H \BQ_j^{t_1} \BQ_j^{t_1, H} \Ba_j \right]^{a_1}\left[ \frac{1}{p} \Tr \BQ_j^{t_2} \BQ_j^{t_2, H} + \Ba_j^H \BQ_j^{t_2} \BQ_j^{t_2, H} \Ba_j \right]^{a_2} \leq K, \forall a_i,  t_j \in \{0, 1, 2\}, \label{Eq_Trace_bilinear_Q_sup_ab_bound}
\end{equation}
where $K$ is a constant independent of $n$. By utilizing these bounds, we can eliminate the polynomial in the denominator of \eqref{Eq_Part_II}. The details are omitted for simplicity.
\par
Since the probabilities of $\{|z\widetilde{q}_j| > 2K^o\}, \{|\eta_j| > 2K^o\}$, and $\{|1 + \alpha_{j, j}|^{-1} > 2K^o\}$ are less than any polynomial of $n^{-1}$ for $x \in [a, b]$ by Proposition \ref{Prop_bound_ab_K_o}, we may assume $|z\widetilde{q}_j|, |\eta_j|, |1 + \alpha_{j, j}|^{-1} \leq K^o$. We write $\BD_T = \diag(t_1, \ldots, t_n)$ and $-z\widetilde{\kappa}_j - (-z[\widetilde{\BR}]_{j,j}) = \sum_{i = 1}^5 X_{ij}$ as in \eqref{Eq_kappa_j_wR_jj}. Using $-z \widetilde{q}_j = \eta_j + z\widetilde{q}_j \eta_j \Delta_j$, we have 
\begin{align*}
  & \sum_{j=1}^n t_j X_{1j} = \sum_{j=1}^n t_j \E z\widetilde{q}_j \eta_j \Delta_j = \sum_{j=1}^n t_j \E (-\eta_j - z\widetilde{q}_j \eta_j \Delta_j) \eta_j \Delta_j  \\
  &= \sum_{j=1}^n -t_j \E  \eta_j^2 \Delta_j - t_j \E z\widetilde{q}_j\eta_j^2 \Delta_j^2 = \sum_{j=1}^n t_j \E z\widetilde{q}_j\eta_j^2  \Delta_j^2. \numberthis
\end{align*}
Therefore, we can obtain 
\begin{align*}
    & \sup_{x \in [a, b]}\abs{\sum_{j=1}^n t_j X_{1j} } \leq T\sup_{x \in [a, b]} \sum_{j=1}^n  \E \abs{z\widetilde{q}_j\eta_j^2  \Delta_j^2}  \lesssim \sup_{x \in [a, b]} \sum_{j=1}^n \E |\Delta_j|^2 \\
    & \lesssim \sup_{x \in [a, b]} \frac{1}{n}\sum_{j=1}^n \E \left[ \frac{1}{p} \Tr \BQ_j \BQ_j^H + \Ba_j^H\BQ_j\BQ_j^H\Ba_j \right] \overset{(a)}{\lesssim} 1, \numberthis \label{Eq_sup_ab_tX1j}
\end{align*}
where step $(a)$ follows from the same argument as in \eqref{Eq_I_j_TrQQ} and \eqref{Eq_E_j_ajQjQjaj}. Next, we evaluate $\sup_{x \in [a, b]}| \sum_{j} t_j X_{2j}|$. By using $\BQ = \BQ_j + z\widetilde{q}_j \BQ_j \Bxi_j \Bxi_j^H \BQ_j$, the following holds by the similar argument as in \eqref{Eq_sup_ab_tX1j}
\begin{align*}
  &\sup_{x \in [a, b]}\abs{ \sum_{j=1}^n t_j X_{2j}} = \sup_{x \in [a, b]}\abs{- \sum_{j=1}^n t_j \E\left( \frac{\eta_j \Ba_j^H\BQ_j\BOmega_j\BQ_j
  \Ba_j}{n (1 + \alpha_{j,j})^2} \right) + \sum_{j=1}^n t_j X_{2j, 3}} \\
  &\lesssim \sup_{x \in [a, b]}\abs{\sum_{j=1}^n \frac{1}{n}\E{ \Ba_j^H\BQ_j\BOmega_j\BQ_j
  \Ba_j}} + \sup_{x \in [a, b]}  \sum_{j=1}^n \E |\Delta_j| \abs{\Ba_j^H\BQ_j\Bxi_j\Bxi_j^H\BQ_j
  \Ba_j - \frac{1}{n}\Ba_j^H \BQ_j \BOmega_j \BQ_j \Ba_j - \Ba_j^H \BQ_j \Ba_j\Ba_j^H \BQ_j \Ba_j} \\
  &+ \sup_{x \in [a, b]} \sum_{j=1}^n \E |\Delta_j|  \left(\frac{1}{n}\Ba_j^H \BQ_j \BOmega_j \BQ_j \Ba_j + \Ba_j^H \BQ_j \Ba_j\Ba_j^H \BQ_j \Ba_j\right) \lesssim 1. \numberthis
\end{align*}
The evaluation for $\sup_{x \in [a, b]} |\sum_j t_j X_{3j}|$ is given by
\begin{align*}
  &\sup_{x \in [a, b]} \abs{\sum_{j=1}^n t_j X_{3j}} = \sup_{x \in [a, b]}\abs{\sum_{j=1}^n \E \frac{t_j (\frac{1}{n} \Tr \BQ_j - \frac{1}{n} \Tr \BQ)}{(1 + \alpha_{j,j})(1 + \kappa_j)}} \\
  &= \sup_{x \in [a, b]} \abs{\sum_{j=1}^n \frac{1}{n} \E \frac{-zt_j\widetilde{q}_j \Bxi_j^H\BQ_j^2\Bxi_j}{(1 + \alpha_{j,j})(1 + \kappa_j)}} \lesssim \sup_{x \in [a, b]} \abs{\sum_{j=1}^n \frac{1}{n} \E \Bxi_j^H \BQ_j \BQ_j^H \Bxi_j}  \lesssim 1. \numberthis \\
\end{align*}
By writing $\BD_{\kappa} = \diag(1 / (1 + \kappa_j); j \leq n)$ and using \eqref{Eq_EQ_R}, we have 
\begin{align*}
  \sup_{x \in [a, b]} \abs{ \sum_{j=1}^n t_j X_{4j} } = \sup_{x \in [a, b]} \abs{ \sum_{j=1}^n \frac{t_j \Ba_j^H (\BR - \E \BQ) \Ba_j}{(1 + \kappa_j)^2}} = \sup_{x \in [a, b]} \abs{ \Tr (\BR - \E \BQ)\BA \BD_T \BD_{\kappa}^2 \BA^H} \lesssim 1. \numberthis
\end{align*}
The evaluation of $\sup_{x \in [a, b]}|\sum_j t_j X_{5j}|$ is quite similar to that of $\sup_{x \in [a, b]}|\sum_j t_j X_{3j}|$  and is omitted here. Therefore, we have completed the proof. \qed
\par
Next, we give a more tight evaluation for \eqref{Eq_I_minus_z_Phi_delta} and \eqref{Eq_I_minus_J_delta_inv} with $x \in [a, b]$. By the integration representation of $\delta_j$, we have 
\begin{equation}
    \sup_{x \in [a, b]} v^{-1}\Im(\delta_j(z)) = \sup_{x \in [a, b]} \int_{\mathbb{R}^+} \frac{\mu_j(\dd \lambda)}{(x - \lambda)^2 + v^2} \leq \underline{\epsilon}^{-2} \mu_j(\R^+) \lesssim 1. \numberthis
\end{equation}
% and 
% \begin{equation}
%     \sup_{x \in [a, b]} v^{-1}\Im(z\delta_j(z)) = \sup_{x \in [a, b]} \int_{\mathbb{R}^+} \frac{\lambda \mu_j(\dd \lambda)}{(x - \lambda)^2 + v^2} \leq \underline{\epsilon}^{-2} \int \lambda_{\mathbb{R}^+} \mu_j (\dd \lambda) \overset{(a)}{\lesssim} 1,
% \end{equation}
% where step $(a)$ follows form the tightness of the measure $\mu_j$. The same conclusion also holds for $v^{-1}\Im(\widetilde{\delta}_j)$ and $v^{-1}\Im(z\widetilde{\delta}_j)$. By \eqref{Eq_I_minus_z_Phi_delta}, we have 
% \begin{equation}
%    \sup_{x \in [a, b]} \norm{(\BI_n - |z|^2 \BPhi_{\delta})^{-1}}_{\infty} \leq \sup_{x \in [a, b]}\frac{\max_j \Im(z \widetilde{\delta}_j)}{v|z|^2} \frac{1}{\inf_{x \in [a, b]}\min_j |\widetilde{\delta}_j|} \lesssim 1,
% \end{equation}
% since $\liminf_n \min_j |\widetilde{\delta}_j| \geq K > 0$ for $x \in [a, b]$ according to the assumption. 
We then show that $\inf_{x \in [a, b], j \leq n} [\Bv_{\delta}]_j$ is bounded away from 0. By the definition of $\BTheta(z)$, we have 
\begin{equation}
    \BI_p = -z \BTheta - \sum_{l=1}^n \frac{\BTheta\BOmega_l }{n} z \widetilde{\delta}_l - z\BTheta\BA\BF \BA^H . \label{Eq_Theta_I_Identity}
\end{equation}
Taking the norm $\norm{\cdot}_{\Omega, j} := \sqrt{\Tr(\cdot)^H \BOmega_j (\cdot)}$ at both sides of \eqref{Eq_Theta_I_Identity} and using $|\Tr\BC\BD| \leq \norm{\BC} \Tr \BD$ for Hermitian non-negative matrix $\BD$ yields
\begin{align*}
  \sqrt{\Tr \BOmega_j} &\leq  |z|^2 (\Tr \BTheta^H \BOmega_j \BTheta)^{1/2} + \sum_{l=1}^n \frac{|z \widetilde{\delta}_l|}{n} (\Tr \BOmega_l \BTheta^H \BOmega_j\BTheta \BOmega_l)^{1/2} \\
  &+ |z|^2( \Tr \BA\BF^H\BA^H \BTheta^H \BOmega_j \BTheta \BA \BF\BA^H )^{1/2} \\
  &\leq (\Tr \BTheta^H\BOmega_j\BTheta)^{1/2} \left[|z|^2 + \sum_{l=1}^n \frac{|z\widetilde{\delta}_l| \norm{\BOmega_l}}{n} + \frac{\norm{\BA}^2}{\min_{l \leq n} |1 + \delta_l|} \right]. \numberthis
\end{align*}
Since $\liminf_n \min_j |1 + \delta_j| > 0$ for $x \in [a, b]$, we have $\inf_{x \in [a, b]} \frac{1}{n} \Tr \BOmega_j\BTheta\BTheta^H \geq K > 0$ for some constant $K$ independent of $n$.  By \eqref{Eq_I_minus_J_delta_inv}, we get
\begin{equation}
    \sup_{x \in [a, b]} \norm{(\BI_n - \BJ_{\delta} )^{-1}}_{\infty} \leq \frac{\sup_{x \in [a, b]}\max_j(\delta_j)}{v} \frac{1}{\inf_{x \in [a, b]} \min_j[\Bv_{\delta}]_j} \lesssim 1.
\end{equation}
Similarly, we know the max-row norms $\| (\BI_n - \BJ_{\kappa})^{-1} \|_{\infty}$, $\| (\BI_n - |z|^2\BPsi_{\delta})^{-1} \|_{\infty}$, and $\| (\BI_n - |z|^2\BPsi_{\kappa})^{-1} \|_{\infty}$ are of order $O(1)$ uniformly for $x \in [a, b]$. Rewrite the matrices as $\BPsi_{\delta} = \BPsi_{\delta}(z, z^*)$ and  $\BPsi_{\kappa} = \BPsi_{\kappa}(z, z^*)$ with 
\begin{equation}
    [\BPsi_{\delta}(z, z^*)]_{j, l} = \frac{1}{n} [\widetilde{\BF}(z)\BA^H\BTheta(z)\BOmega\BTheta(z^*) \BA \widetilde{\BF}(z^*)]_{l, l}, ~~ [\BPsi_{\kappa}(z, z^*)]_{j, l} = \frac{1}{n} [\widetilde{\BG}(z^*)\BA^H\BR(z^*)\BOmega\BR(z) \BA \widetilde{\BG}(z)]_{l, l}.
\end{equation}
By swapping the roles of $z$ and $z^*$, we can similarly prove that
\begin{equation}
    \sup_{x \in [a, b]} \norm{(\BI_n - |z|^2\BPsi_{\delta}(z^*, z))^{-1}}_{\infty
    }, \sup_{x \in [a, b]}\norm{(\BI_n - |z|^2\BPsi_{\kappa}(z^*, z))^{-1}}_{\infty
    } \lesssim 1.
\end{equation}
As a result, there holds
\begin{equation}
    \sup_{x \in [a, b]} \norm{(\BI_n - |z|^2\BPhi_{\delta})^{-1}}_{1
    }, \sup_{x \in [a, b]}\norm{(\BI_n - |z|^2\BPhi_{\kappa})^{-1}}_{1
    } \lesssim 1,
\end{equation}
due to the identities $\BF(z) \BA \widetilde{\BTheta}(z) = \BTheta(z) \BA \widetilde{\BF}(z)$, $\BG(z) \BA \widetilde{\BR}(z) = \BR(z) \BA \widetilde{\BG}(z)$, and $\BPhi_{\delta} = [\BPsi_{\delta}(z^*, z)]^T$, $\BPhi_{\kappa} = [\BPsi_{\kappa}(z^*, z)]^T$.
\par
In the following, we give a tight bound on $\norm{\Bkappa - \Bdelta}_{\infty}$. To this end, the following estimation is needed 
\begin{equation}
    \sup_{x \in [a, b]}\max_{i, j \leq n} \abs{[\BGamma_{\kappa, \delta} (\BI_n - z^2 \BPhi_{\kappa, \delta})^{-1}]_{i, j}} \lesssim n^{-1}. \label{Eq_max_ij_Gamma_Phi}
\end{equation}
In fact, we have 
\begin{align*}
    &\max_{i, j \leq n} \abs{[\BGamma_{\kappa, \delta} (\BI_n - z^2 \BPhi_{\kappa, \delta})^{-1}]_{i, j}} \leq \max_{i, j \leq n} \abs{ \frac{1}{n^2}\Tr \BOmega_i \BR \BOmega_j \BTheta } \norm{(\BI_n - z^2 \BPhi_{\kappa, \delta})^{-1}}_1 \\
    &\overset{(a)}{\leq} \frac{p(\max_j\norm{\BOmega_{j}})^2}{n^2} \norm{\BR}\norm{\BTheta} \sqrt{\norm{(\BI_n - |z|^2 \BPhi_{\kappa})^{-1}}_1\norm{(\BI_n - |z|^2 \BPhi_{\delta})^{-1}}_1}, \numberthis
\end{align*}
where step $(a)$ follows from Lemma \ref{Lemm_matrix_Spr_maxrow_Norm}. Taking $\sup_{x \in [a, b]}$ on both sides of the above inequality, \eqref{Eq_max_ij_Gamma_Phi} holds.
Hence, by Proposition \ref{Prop_EQ_ER_EwQ_EwR}, we have 
\begin{align*}
    \sup_{x \in [a, b]}\norm{\BGamma_{\kappa, \delta}(\BI_n - z^2\BPhi_{\kappa, \delta})^{-1} \left(z \widetilde{\Bd}\right) }_{\infty} = \sup_{x \in [a, b]} \max_{i \leq n} \abs{\sum_{j = 1}^n [\BGamma_{\kappa, \delta}(\BI_n - z^2\BPhi_{\kappa, \delta})^{-1}]_{i, j}\left(z\widetilde{\kappa}_j - z[\widetilde{\BR}]_{j, j}\right)} \lesssim n^{-1}. \numberthis
\end{align*}
Again, by Lemma \ref{Lemm_matrix_Spr_maxrow_Norm} and Proposition \ref{Prop_EQ_ER_EwQ_EwR}, we can obtain
\begin{align*}
    \sup_{x \in [a, b]}\norm{\Bkappa - \Bdelta}_{\infty} &= \sup_{x \in [a, b]}\norm{(\BI_n - \BJ_{\kappa, \delta})^{-1} \left[\Bd + \BGamma_{\kappa, \delta}(\BI_n - z^2\BPhi_{\kappa, \delta})^{-1} \left(z \widetilde{\Bd}\right) \right]}_{\infty}. \\
    &\leq \sup_{x \in [a, b]}\norm{(\BI_n - \BJ_{\kappa, \delta})^{-1}}_{\infty} \norm{\Bd}_{\infty} + \norm{(\BI_n - \BJ_{\kappa, \delta})^{-1}}_{\infty} \norm{\BGamma_{\kappa, \delta}(\BI_n - z^2\BPhi_{\kappa, \delta})^{-1} \left(z \widetilde{\Bd}\right) }_{\infty} \\
    & \lesssim \sup_{x \in [a, b]}\sqrt{\norm{(\BI_n - \BJ_{ \delta})^{-1}}_{\infty}\norm{(\BI_n - \BJ_{\kappa})^{-1}}_{\infty}} n^{-1} \lesssim \frac{1}{n}. \numberthis
\end{align*}
By \eqref{Eq_Bound_on_Diff_wkappa_wdelta} and the above estimation, we have $\sup_{x \in [a, b]}|\sum_{j=1}^n t_j (\widetilde{\kappa}_j - \widetilde{\delta}_j)| \lesssim 1$ with $\sup_{x \in [a, b], j \leq n} |t_j|$ uniformly bounded in $n$. By the resolvent identity, we have 
\begin{align}
   \Tr \BR - \Tr \BTheta = \sum_{j=1}^n \frac{\Tr \BR \BOmega_j \BTheta}{n}(z\widetilde{\kappa}_j - z\widetilde{\delta}_j) + \sum_{j=1}^n z^2 [\widetilde{\BF} \BA^H\BTheta \BR\BA \widetilde{\BG}]_{j, j}(\kappa_j - \delta_j).
\end{align}
Denoting $\frac{1}{n} \Tr \BR \BOmega_j \BTheta = t_j$, we have $\sup_{x \in [a, b]}\sup_j |zt_j| \lesssim 1$ and  
\begin{equation}
    \sup_{x \in [a, b]}\abs{\Tr (\BR - \BTheta)} \leq \sup_{x\in [a, b]} \abs{\sum_{j}t_j(\widetilde{\kappa}_j - \widetilde{\delta}_j)} +  \sup_{x \in [a, b]} \sup_{j \leq n} \frac{n \norm{\BA}^2 \norm{\BTheta} \norm{\BR}}{|1 + \delta_j||1 + \kappa_j|} \norm{\Bkappa - \Bdelta}_{\infty} \lesssim 1.
\end{equation}
By Proposition \ref{Prop_EQ_ER_EwQ_EwR} and the above inequality, \eqref{Eq_Mean_Part} holds. Therefore, we have completed the proof of 
\eqref{Eq_sup_diff_resolvent} \qed
\section{Proof of Corrolary \ref{Coro_min_eig_zero_forcing}}
\label{App_prof_of_coro_min_eig_zero_forcing}

The crucial step of the proof is to show that the boundedness of the terms listed in Proposition \ref{Prop_bound_ab_K_o} holds when $x = 0$. To that end, we need Theorem 1 from \cite{AblaNoEigen}, which we list below as a proposition. 
\begin{proposition}(\cite[Theorem 1]{AblaNoEigen})
\label{Prop_abla_noeigen}
    Assume the following holds
    \begin{itemize}
    \item [1)]  $0 < \liminf_n p / n \leq  \limsup_n p / n < 1$. 
    \item [2)] $ 0 < \liminf_n \min_j \lambda_{\min}(\BOmega_j) \leq \limsup_n \max_j \norm{\BOmega_j} < \infty$.
    \end{itemize}
    Then, there exists $\epsilon > 0$ suth that 
    \begin{equation}
        [0, \epsilon] \cap \mathsf{Supp}(F^n) = \emptyset
    \end{equation}
    for all large $n$.
\end{proposition}
\textit{Proof:} Following the same argument as in Proposition \ref{Prop_bound_ab_K_o}, it suffices to show that the following holds for $z = x+ \jmath v$, $v \to 0$ as $n \to \infty$
\begin{equation}
    \sup_{x \in [0, \epsilon']} \max_{j \leq n} \left\{ |z\widetilde{\delta}_j(z)|, |\delta_j(z)|, \frac{1}{|1 + \delta_j(z)|} \right\} \leq K, \label{ER_bounde_0e_qterms}
\end{equation}
for some $\epsilon' > 0$. Here, we choose $\epsilon' = \epsilon / 2$. By Proposition \ref{Prop_abla_noeigen} and the matrix valued integration representation \eqref{Eq_representation_Theta}, we can obtain that $[0, \epsilon]$ is outside the support of $F^n = \frac{1}{p} \Tr \Bmu$ and $[\Bmu]_{i, i}([0, \epsilon]) = 0$, for $i = 1, \ldots, p$. Since $\Bmu([0, \epsilon])$ is non-negative definite, $\Bmu([0, \epsilon]) = \Bzero_{p \times p}$ must holds, which implies $[0, \epsilon]$ is outside the support of the measures $\{\mu_j\}_{j \leq n}$. In particular, $\epsilon$ is less than or equal to the left endpoint of $\cup_j \mathsf{Supp}(\mu_j)$, which implies 
\begin{equation}
    \sup_{x \in [0, \epsilon']} |\delta_j(z)| \leq \int_{\R^+} \frac{\mu_j(\dd \lambda)}{|\lambda - z|} \leq \frac{2\Tr \BOmega_j}{n \epsilon},
\end{equation}
and 
\begin{equation}
    \inf_{x \in [0, \epsilon']} \delta_j(x) = \int_{\R^+} \frac{\mu_j(\dd \lambda)}{\lambda - x} \geq 0.
\end{equation}
Hence, for large $n$, it holds that 
\begin{equation}
    \frac{1}{|1 + \delta_j(z)|} \leq \frac{1}{|1 + \delta_j(x)| - |\delta_j(z) - \delta_j(x)|} \leq \frac{1}{1 - v \frac{4 \Tr \BOmega_j}{\epsilon^2n}} \leq 1 / 2.
\end{equation}
\par
Since $\BTheta(z)$ admits the representation $\int_{\R^+} \frac{\Bmu(\dd \lambda)}{\lambda - z}$ and $\Bmu([0, \epsilon]) = 0$, by the same argument as in \eqref{Eq_sup_ab_norm_Theta}, we have 
\begin{equation}
    \sup_{x \in [0, \epsilon']} \norm{\BTheta(z)} \leq \sup_{x \in [0, \epsilon']} \norm{\int_{\R^+} \frac{\Bmu(\dd \lambda)}{|\lambda - z|}} \leq \frac{2}{\epsilon}.
\end{equation}
Therefor, by the identity in \eqref{Eq_Res_identity_wdeltaj}, we can obtain
\begin{equation}
    \abs{z\widetilde{\delta}(z)} \leq \frac{1}{|1 + \delta_j(z)|} + \frac{\norm{\Ba_j}^2 \norm{\BTheta(z)}}{|1 + \delta_j(z)|^2},
\end{equation}
which is also uniformly upper bounded for $x \in [0, \epsilon']$.
\par
Hence, we have proved \eqref{ER_bounde_0e_qterms}. To prove that there are no eigenvalues of $\BS$ in the interval $[0, \epsilon']$, we only need to follow the same logic as in Appendix \ref{App_Proof_Thm_noeigen} by replacing the interval $[a, b]$ with $[0, \epsilon']$, and $[a', b']$ with $[- \epsilon'/2, 3\epsilon'/2]$. \qed

\bibliographystyle{IEEEtran}
\bibliography{reference}
\end{document}